\documentclass[a4paper,10pt]{amsart}
\usepackage{palatino, mathpazo}
\usepackage[plainpages=false]{hyperref}
\usepackage{amsfonts,latexsym,rawfonts,amsmath,amssymb,amsthm, mathrsfs, lscape}
\usepackage{verbatim}
\usepackage{stmaryrd}

\usepackage[all]{xy}
\usepackage{graphicx,psfrag}

\usepackage{array, tabularx}

\usepackage{setspace}

\usepackage{epsfig}
\usepackage{graphicx}
\usepackage{eufrak}
\usepackage{dsfont}
\usepackage[english]{babel}
\usepackage{color}

\usepackage{palatino, mathpazo}
\usepackage{amscd} % to be able to use "CD" environment
\usepackage{amssymb}
\usepackage{dsfont}
\usepackage{xypic} % to be able to use "diagram" environment of
\LaTeXdiagrams % xypic to typeset commutative diagrams
\newbox\mybox
\def\overtag#1#2#3{\setbox\mybox\hbox{$#1$}\hbox to
 0pt{\vbox to 0pt{\vglue-#3\vglue-\ht\mybox\hbox to \wd\mybox
 {\hss$\ss#2$\hss}\vss}\hss}\box\mybox}
\def\undertag#1#2#3{\setbox\mybox\hbox{$#1$}\hbox to 0pt{\vbox to
 0pt{\vglue#3\vglue\ht\mybox\hbox to \wd\mybox
 {\hss$\ss#2$\hss}\vss}\hss}\box\mybox}
\def\lefttag#1#2#3{\hbox to 0pt{\vbox to 0pt{\vss\hbox to
 0pt{\hss$\ss#2$\hskip#3}\vss}}#1}
\def\righttag#1#2#3{\hbox to 0pt{\vbox to 0pt{\vss\hbox to
 0pt{\hskip#3$\ss#2$\hss}\vss}}#1}
\let\ss\scriptstyle

\def\Dot{\lower.2pc\hbox to 2pt{\hss$\bullet$\hss}}
\def\Circ{\lower.2pc\hbox to 2pt{\hss$\circ$\hss}}
\def\Vdots{\raise5pt\hbox{$\vdots$}}
\newcommand\lineto{\ar@{-}}
\newcommand\dashto{\ar@{--}}
\newcommand\dotto{\ar@{.}}

\allowdisplaybreaks[3]

\newtheorem{prop}{Proposition}[section]
\newtheorem{lem}[prop]{Lemma}

\newtheorem{cor}[prop]{Corollary}
\newtheorem{thm}[prop]{Theorem}
\newtheorem{rmk}[prop]{Remark}
\newtheorem{example}{Example}

\newtheorem{question}[prop]{Question}

\newtheorem{defn}[prop]{Definition}

\newtheorem{con}[prop]{Conjecture}

 {\nolinebreak $\Box$ \end{trivlist}}

\newcommand{\noprint}[1]{}

\renewcommand{\tilde}{\widetilde}

\newcommand{\Ext}{\mbox{Ext}}

\newcommand{\pt}{\mathop{pt}}

\newcommand{\XX}{{\mathfrak X}}

\newcommand{\MM}{{\mathfrak M}}

\newcommand{\zz}{{\mathbb Z}}
\newcommand{\hh}{{\mathbb H}}

\newcommand{\aaa}{{\mathbb A}}
\newcommand{\ttt}{{\mathbb T}}

\renewcommand{\ll}{{\mathbb L}}
\newcommand{\qq}{{\mathbb Q}}
\newcommand{\pp}{{\mathbb P}}
\newcommand{\cc}{{\mathbb C}}
\newcommand{\rr}{{\mathbb R}}

\newcommand{\sT}{{\mathcal T}}

\newcommand{\sC}{{\mathcal C}}
\newcommand{\sE}{{\mathcal E}}

\newcommand{\sO}{{\mathcal O}}

\newcommand{\sX}{{\mathcal X}}
\newcommand{\sY}{{\mathcal Y}}

\newcommand{\sZ}{{\mathcal Z}}

\newcommand{\sA}{{\mathcal A}}

\DeclareMathOperator{\Sch}{Sch}

\DeclareMathOperator{\Aut}{Aut}

\DeclareMathOperator{\ind}{ind}

\DeclareMathOperator{\vir}{vir}

\DeclareMathOperator{\vd}{vd}

\DeclareMathOperator{\CM}{CM}

\DeclareMathOperator{\lci}{lci}

\DeclareMathOperator{\Bub}{Bub}
\DeclareMathOperator{\rt}{rt}

\newcommand{\pr}{\mathop{\rm pr}\nolimits}

\renewcommand{\top}{\mathop{\rm top}}

\renewcommand{\div}{\mathop{\rm div}}

\newcommand{\spec}{\mathop{\rm Spec}\nolimits}

\numberwithin{equation}{subsection}

\newcommand {\mat} [1] {\left(\begin{array}{#1}}
\newcommand {\rix} {\end{array}\right)}

\begin{document}

\title{A Bubble Compactification of  KSBA Moduli Space of Surfaces of General Type}

\author{Yunfeng Jiang}
\address{Department of Mathematics\\ University of Kansas\\ 405 Snow Hall 1460 Jayhawk Blvd\\Lawrence KS 66045 USA}
\email{y.jiang@ku.edu}
\date{\today}

\begin{abstract}
We construct a bubble-tree enhancement of the KSBA moduli space of surfaces of general type using index-one-cover Deligne--Mumford stacks. The construction is inspired by bubbling phenomena for K\"ahler--Einstein metrics. Non-klt semi-log-canonical singularities are replaced by admissible simple-normal-crossing configurations: the cusp and degenerate-cusp components are rational, while the simple-elliptic case also contains elliptic ruled components. A key feature is that the construction avoids lci covers.

We construct the moduli stack of   polarized bubble-tree index-one-cover Deligne-Mumford stacks. The stack is in general not separated, and we construct its underlying separated algebraic  space. 
We also construct the moduli stack of bubble-tree index-one-cover Deligne-Mumford stacks  with labeled sections that record algebraic bubbling centers and compare it with the pointed metric bubble tree. The labeled sections determine a numerical class whose stability proves the separatedness of the moduli stack. 

We establish a perfect obstruction theory and a virtual fundamental class on the bubble-tree moduli stacks. This provides a novel and geometric pathway to prove the existence of a virtual fundamental class on the original KSBA moduli space, thus offering a new approach to a conjecture of Donaldson.  We expect that parts of the construction extend to higher dimensions.
\end{abstract}

\maketitle

\section{Introduction}

We work over $\cc$.
The main purpose of this paper is to construct a new compactification of the KSBA moduli space of surfaces of general type. We define a polarized bubble-tree functor with a contraction to the KSBA moduli stack. We prove that this functor is represented by   a finite-type Deligne--Mumford stack, which carries a perfect obstruction theory and a  virtual fundamental class. Our approach does not use the lci covers introduced in \cite{Jiang_2022}; instead, it is inspired by bubble-tree structures in K\"ahler geometry. We begin with the main results.

\subsection{Main results}
The KSBA moduli functor $M_N:=\overline{M}_{K^2,\chi,N}$ parametrizes $\qq$-Gorenstein families $\sX\to T$ of stable semi-log-canonical (slc) surfaces $X$ for which $\omega_{\sX/T}^{[N]}=(\omega_{\sX/T}^{\otimes N})^{\vee\vee}$ is invertible, with fixed invariants $K_X^2=K^2$ and $\chi(\sO_X)=\chi$; see \cite{Kollar-Shepherd-Barron}.

Hacking introduced the canonical stack associated with an slc surface $X$ \cite{Hacking}; here we call it the index-one-cover Deligne--Mumford stack $\XX$. These stacks are Gorenstein and behave well in flat families.
The moduli stack $M_N^{\ind}:=\overline{M}^{\ind}_{K^2, \chi, N}$ of index-one covers, defined in \cite{AH,Jiang_2022}, parametrizes flat families of index-one-cover Deligne--Mumford stacks.
The stack $M_N^{\ind}$ is isomorphic to $M_N$.

In this paper, we define a moduli problem
 $M_N^{\Bub}$ of bubble-tree surfaces over the moduli stack $M_N^{\ind}$
of index-one covers. Similar bubbling phenomena appear in other moduli problems; see \cite{Feehan,MTT,Taubes,Uhlenbeck}. The bubble-tree construction avoids the lci covers used in \cite{Jiang_2022} and applies to the semi-log-canonical singularities considered here. In the cusp case, Inoue--Hirzebruch surfaces describe the analytic dual-cusp geometry, while the actual bubble components used in the algebraic moduli problem come from Type III models. The construction is motivated by rescaled limits of K\"ahler--Einstein metrics on surfaces of general type near non-klt log-canonical singularities.
The idea in this paper can also be applied to the compactification of K3 surfaces \cite{Jiang_K3}.

Let $f:\sX\to\Delta$ be a flat one-parameter smoothing of slc surfaces. Suppose that the central fiber $\sX_0$ contains simple elliptic, cusp, or degenerate cusp singularities $p_1,\ldots,p_n$ of embedding dimension at least six; these singularities may have nonzero higher obstruction spaces \cite{Jiang_2021}. 
We construct a bubble-tree stack $\XX_0^{\Bub}$ from the central fiber $\XX_0$ of the  index one cover Deligne-Mumford stack $f^{\ind}: \XX\to \Delta$ of $f: \sX\to \Delta$. 
Our first result is the following.

\begin{thm}\label{thm_smoothing_bubble_intro}(Theorem \ref{thm_bubble_tree_DM_smoothing})
Assume that the local models at the simple elliptic, cusp, and degenerate cusp singularities can be chosen to form an admissible bubble-tree datum in the sense of Definition \ref{defn_admissible_bubble_tree}. Then these singularities can be replaced by simple-normal-crossing bubble components. If, in addition, the $d$-semistability, global obstruction-vanishing, and equivariant compatibility hypotheses of Theorems \ref{thm_smoothing_bubble_tree} and \ref{thm_bubble_tree_DM_smoothing} hold, the resulting index-one-cover bubble-tree stack $\XX_0^{\Bub}$ admits a smoothing $f:\mathscr X^{\Bub}\to\Delta$ with central fiber $\XX_0^{\Bub}$.
\end{thm}

Within the local quotient types treated in this paper, every relevant singularity germ of $\XX_0$ is locally of the form $([Z/\mu_r],p_s)$, where $r\in\{1,2,3,4,6\}$ and $(Z,p_s)$ is a simple elliptic, cusp, or degenerate cusp singularity.
We summarize the construction as follows.
\begin{enumerate}
 \item For a simple elliptic singularity $(\sX_0,x)$, we
replace $x$ by a chain whose intermediate components are elliptic ruled surfaces and whose terminal component is a degree-$d$ del Pezzo surface $Y$ containing a smooth anticanonical elliptic curve $E$. The open surface $Y\setminus E$ occurs, after hyperk\"ahler rotation, as the bubble constructed in \cite{FHJ}.

\item For a cusp singularity $(\sX_0,x)$
that admits a smoothing, we use Looijenga's conjecture, proved in \cite{GHK15,Engel}, to replace the singularity by a Type III degeneration of log Calabi--Yau surfaces; see \S\ref{subsec_bubble_tree}. The correspondence between smoothable cusps and Looijenga pairs carrying the dual cusp cycle is therefore central to the construction.

Our construction also has a differential-geometric interpretation. The sequence of K\"ahler--Einstein metrics is expected to collapse near the cusp singularity. In \cite[Corollary 1.11]{DFS}, the authors prove that a complete K\"ahler--Einstein metric near the cusp
$x$ is asymptotic to $(\hh\times \hh)/\Gamma$, where
$\hh$
is the upper half-plane and
$\Gamma$ is a discrete parabolic subgroup of $\Aut(\hh\times\hh)$. The space
$(\hh\times \hh)/\Gamma$ is a noncompact manifold with a cusp end. We conjecture that the Type III model records suitable compactifications of the rescaled K\"ahler--Einstein bubbles.
\item For degenerate cusp singularities $(\sX_0,x)$,
suppose that they admit a one-parameter smoothing. We use a crepant semistable model to form the bubble components. The replacement is an open Kulikov surface, whose mirror-symmetry geometry is used in \cite{AAB2024}
to construct KSBA families of log Calabi--Yau surfaces. We attach this open Kulikov surface along the resolution cycle of the degenerate cusp.

It is expected that this open Kulikov surface is the bubbling limit structure in \cite{Sun}.
The algebraic-geometric study of the bubble limits can be found in \cite{dBS24} and \cite{Odaka}.
\end{enumerate}

The naive moduli problem has two defects: crepant modifications make it nonseparated, and unmarked log Calabi--Yau bubble components may have positive-dimensional automorphism groups.  We therefore impose a polarization together with a stable auxiliary divisor.  Fix
\[
\tau=(K^2,\chi,N,P,m,\epsilon,\mathcal I,P_c),
\]
where $P$ is the Hilbert polynomial of a relatively very ample line bundle, $N$ bounds the Cartier and stabilizer indices, and $\mathcal I$ is a finite set of decorated bubble types.  The coefficient $0<\epsilon\leq1$, the integer $m>0$, and the Hilbert polynomial $P_c$ of the contraction graph are also fixed. The enhanced stack structures are assumed to admit a finite-type presentation.

The main moduli functor is
\[
M_\tau^{\Bub,\mathrm{pol}}:\Sch_{\cc}^{\mathrm{op}}
\longrightarrow\mathrm{Groupoids}.
\]
For a scheme $T$, its objects are tuples
\[
(f:\mathscr X^{\Bub}\to T,\ c:\mathscr X^{\Bub}\to\mathscr X,
\ \mathscr L,\ \mathscr H,\ \xi).
\]
Here $f$ is a flat family of admissible index-one-cover bubble-tree Deligne--Mumford stacks; $c$ contracts the designated bubbles to a KSBA family $\overline{f}: \mathscr{X}\to T$; $\mathscr L$ descends to a relatively very ample line bundle with Hilbert polynomial $P$; and $\mathscr H\in|\mathscr L^{\otimes m}|$ descends to a divisor $H$ for which $(X_t^{\Bub},\epsilon H_t)$ is a stable pair.  Finally, $\xi$ is a first-order pre-smoothing direction that is nonzero along the strata required to smooth.   Isomorphisms preserve all five pieces of data, including the labels in $\mathcal I$.

This definition separates the ambient moduli problem from its smoothing locus.  The ambient stack contains obstructed pre-smoothings; the actual smoothing locus consists of those points for which every successive obstruction to extending $\xi$ vanishes and the resulting formal deformation is effective.  The stable pair condition, rather than the polarization alone, makes automorphism groups finite.  The fixed Hilbert polynomial, bounded index, and finite decorated set give boundedness.

The bubble-tree model should be viewed as a refinement of semistable reduction for an index-one-cover family.  Its existence is governed by the admissibility and smoothing conditions in Section \ref{sec_flat_family_bubble_tree}.  Different relative minimal models can be related by the birational  modifications of Types $(M0)$, $(M1)$, and $(M2)$ described in Section \ref{subsec_birational_modification}.  We call two such models $\mathbb{S}$-equivalent if they extend the same punctured one-parameter family and are related by the allowed crepant modifications.  The stack parametrizes actual polarized objects up to isomorphism; it does not quotient its groupoids by $\mathbb{S}$-equivalence.  That relation is instead reserved for a prospective algebraic space or a canonical closed representative.  Its expected relation to metric bubble trees is stated in Conjecture \ref{con_uniqueness}.

The first moduli theorem takes the following form; the precise hypotheses are collected in Theorem \ref{thm_moduli_stack}.

\begin{thm}\label{thm_moduli_Bub_intro}(Theorem \ref{thm_moduli_stack})
For fixed numerical data $\tau$, the stability condition on the auxiliary pair gives finite automorphism groups, and the fixed Hilbert polynomial and decorated local types give boundedness. The polarized functor $M_\tau^{\Bub,\mathrm{pol}}$ is represented by a finite-type Deligne--Mumford stack. 

The moduli stack $M_\tau^{\Bub,\mathrm{pol}}$ is, in general, nonseparated.  But there is a  contraction morphism 
$\mathbf{c}: M_\tau^{\Bub,\mathrm{pol}}\to M_N$.
\end{thm}

In the  moduli functor
$M_\tau^{\Bub,\mathrm{pol}}$, 
for a scheme $T$, and the family $f:\mathscr X^{\Bub}\to T$, even after we take the relative minimal model $f:\mathscr X_{\min}^{\Bub}\to T$,  different relative minimal models $f:\mathscr X_{\min}^{\Bub}\to T$ may be still related by type  $(M0), (M1)$ and $(M2)$ birational modifications. This causes the ``nonseparatedness" of the moduli. 

It is interesting to find some stability conditions of bubble-tree index-one-cover Deligne-Mumford stacks so that under the stability condition all the objects in the $\mathbb{S}$-equivalence class are isomorphic or the objects can be separated. 

Out first method is to  
use $\mathbb{S}$-equivalent relation on $M_\tau^{\Bub,\mathrm{pol}}$ to define a ``separated" algebraic space. 
We have 

\begin{thm}\label{thm_good_model}(Theorem \ref{thm_stack_good_moduli})
There exists a separated algebraic sapce $$\pi: M_\tau^{\Bub,\mathrm{pol}}\to M_N^{\mathrm{Bub}}.$$
If $N$ is sufficiently large and divisible, $M_N^{\mathrm{Bub}}$ is compact. 

The universal family 
$p^{\mathrm{\Bub}}: \mathscr{X}^{\mathrm{Bub}}\to M_\tau^{\Bub,\mathrm{pol}}$ induces a family 
    $$
    p^{\mathrm{alg}}: \mathscr{X}^{\mathrm{alg}}\to M_N^{\mathrm{Bub}}.
$$
 Then 
there exists a perfect obstruction theory 
\begin{equation}\label{eqn_OT_univ_intro_canonical}
\phi^{\mathrm{Bub}}: E^{\bullet}_{M_N^{\mathrm{Bub}}}\to \ll_{M_N^{\mathrm{Bub}}}^{\bullet}
\end{equation}
where 
$$E^{\bullet}_{M_N^{\mathrm{Bub}}}=Rp^{\mathrm{alg}}_{*}(\ll^{\bullet}_{\mathscr{X}^{\mathrm{alg}}/M^{\mathrm{Bub}}}\otimes \omega^{\mathrm{alg}})[-1].$$
\end{thm}

In the construction of the bubble-tree surface $X^{\Bub}$,
the number of bubbling surfaces depends on the smoothing of the singularity, the embedding dimension, and the self-intersection number of the resolution cycles. This is motivated by the behavior of bubbling near the singularities. The construction of the Type III degeneration for cusp singularities gives the direction for the bubbling limits of collapsing K\"ahler--Einstein metrics in differential geometry.

An interesting point is that we use the smoothing components of Type III degeneration of log Calabi--Yau surfaces—motivated by the smoothing of the cusp singularity—to construct smoothing of the bubble-tree surfaces in the cusp singularity case.

This construction is also inspired by mirror symmetry: the mirror of the smoothing of a cusp singularity corresponds to a vanishing sphere $\mathbb{S}^2$ in the dual SYZ fibration. This topological sphere $\mathbb{S}^2$
can be viewed as an integral affine sphere with singularities, where the number of singularities is determined by the cusp; see \cite{Aroux}.

\subsection{Connection to other moduli stacks}

\begin{con}
 The moduli stack $M_\tau^{\Bub,\mathrm{pol}}$ is obtained from the KSBA stack $M_N$ by a finite sequence of blow-ups along loci parametrizing surfaces with simple elliptic or cusp singularities of embedding dimension at least six.

 Moreover, there is a birational map
 $M_{N}^{\Bub}\dashrightarrow M^{\lci}_N$ to the moduli stack of lci covers, and
 a birational map
 $M_{N}^{\Bub}\dashrightarrow M^{\mathrm{alg}}_N$ to the algebraic space of the moduli stack of birational $\lci$ enhancements
 defined in \cite{Jiang_2022}. All stacks admit contraction morphisms to $M_N$, giving the diagram
 \[
 \xymatrix{
 M_{N}^{\Bub}\ar@{-->}[rr]\ar[drr]_{\mathbf{c}}\ar@{-->}[d]&& M_N^{\lci}\ar[d]^{f^{\lci}}\\
  M^{\mathrm{alg}}_N\ar[rr]_{f^{\mathrm{alg}}} && M_N.
 }
 \]
\end{con}

In \cite{ABJ}, we will show that the wall crossing between the moduli stack of twisted stable maps to the moduli stack $\overline{M}_{1,1}$ of stable genus-one curves with one marked point and the moduli stack of stable quasimaps to $\overline{M}_{1,1}$ provides evidence for the conjecture above. These are moduli stacks of lci covers over the KSBA moduli stack of elliptic fibered surfaces.

\subsection{Uniqueness of the bubble limit structure}\label{subsec_uniqueness_bubble}

S. Sun proves finite termination of the bubble tree in the non-collapsing K\"ahler--Einstein setting \cite{Sun}; this does not imply uniqueness of the entire bubble tree. Motivated by that result, we formulate a separate conjecture for smoothings of non-klt surface singularities.

In the cusp case, Type III degenerations of log Calabi--Yau surfaces are related by birational modifications; analogous nonuniqueness occurs for open Kulikov models of degenerate cusps. Thus two distinct central fibers may extend the same family over a punctured disk. This is the reason for introducing $\mathbb{S}$-equivalence.

For a flat family $f:\mathscr{X}^{\Bub}\to\Delta$, let
$\hat f:\mathscr{X}_{\min}^{\Bub}\to\Delta$ denote a relative minimal model obtained by the birational modifications allowed in Definition \ref{defn_M012_flop}. To discuss pointed bubble limits one must also fix a section $s:\Delta\to\mathscr X$ through the region in which curvature concentrates. Both the marking and the choice of rescaling sequence are essential data.

\begin{con}[Algebraic determination of marked bubble trees]\label{con_uniqueness}
Let $f:\mathscr X\to\Delta$ be a $\qq$-Gorenstein smoothing of a stable surface with a simple elliptic, cusp, or degenerate cusp singularity $p$, and fix a section $s$ specializing to $p$. Assume that the negative K\"ahler--Einstein metrics on the smooth fibers admit a finite pointed bubble tree along $s$. Then the following hold.
\begin{enumerate}
\item The marked metric bubble tree is expected to be compatible with a relative minimal bubble-tree model of $(\mathscr X,s)$. Recovering its metrics requires period, asymptotic, and scale data in addition to an algebraic $\mathbb{S}$-equivalence class; no bijection from unadorned central fibers is asserted.
\item For a simple elliptic singularity, the non-conical vertex adjacent to the main component is the $ALH^*$ bubble described by the corresponding del Pezzo pair $(Y,E)$.
\item For a cusp singularity, after hyperk\"ahler rotation and compactification of each non-conical Ricci-flat vertex, the resulting Looijenga pairs assemble, up to Type III modifications, into the Type III bubble-tree model of the dual cusp.
\item For a degenerate cusp, the analogous compactified vertices assemble into the relative minimal open Kulikov model.
\end{enumerate}
Every compactified vertex $(Y,D;L)$ is expected to be weakly open K-polystable, in the sense of \cite[Definition 2.3]{Odaka2}, for the polarization induced by the degeneration.
\end{con}

\begin{rmk}
The last sentence uses weak \emph{open} K-polystability, not ordinary log K-polystability. Theorem 1.3 of \cite{Odaka2} gives log K-semistability for polarized dlt log Calabi--Yau pairs; it does not give the stronger polystability assertion. The latter is therefore retained as part of the conjecture and is stated with the polarization included.
\end{rmk}

\begin{rmk}
 For a simple elliptic singularity, the replacement is a chain whose end components are degree-$d$ del Pezzo surfaces and whose intermediate components are elliptic ruled surfaces. The identification with the metric bubble should be understood in the precise sense of Conjecture \ref{con_uniqueness}; the analytic model of the first non-conical bubble is provided by \cite{FHJ}.
\end{rmk}

The pointed theory requires an additional moduli construction, which motivates our second method to deal with ``nonseparatedness".  A section through a smoothing records the center at which the metric is rescaled; after semistable reduction, its specialization may lie on a distinguished bubble component.  The algebraic section does not, by itself, record the numerical rescaling sequence or prove that this component is a metric bubble, but these algebraic sections  really encodes the degeneration/smoothing direction that the underlying special fiber by itself forgets.  In Section \ref{sec_marked_bubble_moduli} we therefore construct a separated moduli  stack of polarized bubble-tree Deligne--Mumford surfaces with labeled bubbling sections and compare it with $M_\tau^{\Bub,\mathrm{pol}}$.  
An interesting result is that there exists a algebraic numerical class 
$$A(\mathbf{s})$$ for the 
the smoothing family $f: \mathscr{X}^{\Bub}\to \Delta$ of polarized bubble-tree Deligne--Mumford stacks with labeled bubbling sections $\mathbf{s}=(s_1,\cdots,s_r)$, such that its stability with respect to the flopping curves determines, up to isomorphism, a unique flat family. 
This proves the ``separatedness" of the moduli stack with sections, see \S \ref{subsec_moduli_stack_bubbling_section}. 
The proposed identification with the pointed bubbles of \cite{Sun} is isolated as a conditional comparison statement.

We have 

\begin{thm}\label{prop_marked_stack_construction_intro}(Theorem \ref{prop_marked_stack_construction})
Let $\mathfrak{M}^{\Bub,\mathrm{sec}}_{\tau,\sigma, A_{\mathbf{s}}}$ be 
the fiber category of  polarized bubble-tree Deligne--Mumford surfaces with labeled bubbling sections $\sigma$.  Then $\mathfrak{M}^{\Bub,\mathrm{sec}}_{\tau,\sigma, A_{\mathbf{s}}}$ is represented by  a separated finite-type Deligne--Mumford stack.  The morphism forgetting the sections,
\[
F:\mathfrak{M}^{\Bub,\mathrm{sec}}_{\tau,\sigma, A_{\mathbf{s}}}
\longrightarrow M^{\Bub,\mathrm{pol}}_\tau,
\]
is representable. 
Moreover, there exists a proper morphism 
$$
\mathbf{c}_{\mathrm{sec}}: \mathfrak{M}^{\Bub,\mathrm{sec}}_{\tau,\sigma, A_{\mathbf{s}}}\to M_N
$$
if $N$ is sufficiently large and divisible. 
\end{thm}

\subsection{Virtual fundamental class}\label{subsec_pot_tautological_intro}

We now describe that the polarized moduli stack carries a perfect obstruction theory.

For the moduli stack $M_\tau^{\Bub,\mathrm{pol}}$, let $p:\mathscr X^{\Bub}\to M_\tau^{\Bub,\mathrm{pol}}$ be the universal family. It is flat and projective; we additionally assume that it is relatively Gorenstein. Let $\ll^{\bullet}_{\mathscr{X}^{\Bub}/M_\tau^{\Bub,\mathrm{pol}}}$ be its relative cotangent complex and set
$\omega^{\Bub}:=\omega_{\mathscr{X}^{\Bub}/M_\tau^{\Bub,\mathrm{pol}}}[2]$.

Let
$$E^{\bullet}_{M_\tau^{\Bub,\mathrm{pol}}}:=Rp_{*}\left(\ll^{\bullet}_{\mathscr{X}^{\Bub}/M_\tau^{\Bub,\mathrm{pol}}}\otimes\omega^{\Bub}\right)[-1],$$
where the relative dualizing sheaf satisfies
$$\omega_{\mathscr{X}^{\Bub}/M_\tau^{\Bub,\mathrm{pol}}}|_{p^{-1}(t)}\cong \omega_{\mathscr{X}^{\Bub}_t}.$$
Then the dualizing sheaf $\omega_{\mathscr{X}^{\Bub}_t}$ of the index-one-cover bubble-tree Deligne--Mumford stack $\mathscr{X}^{\Bub}_t$, which is locally given by $\omega_{X_t}^{[r]}$ at a singularity germ ($r$ is the index of the singular germ),
is invertible.

The Kodaira--Spencer map $\ll^{\bullet}_{\mathscr{X}^{\Bub}/M_\tau^{\Bub,\mathrm{pol}}}\to p^{*}\ll_{M_\tau^{\Bub,\mathrm{pol}}}^{\bullet}[1]$
induces a morphism
$$\phi:E^{\bullet}_{M_\tau^{\Bub,\mathrm{pol}}}\to \ll_{M_\tau^{\Bub,\mathrm{pol}}}^{\bullet}.$$
The lci condition makes the relative cotangent complex perfect, hence making the complex $E^{\bullet}_{M_\tau^{\Bub,\mathrm{pol}}}$ perfect of amplitude $[-1,0]$, Thus,  $\phi$ satisfies the Behrend--Fantechi perfect obstruction theory axioms and  is a perfect obstruction theory in the sense of \cite{BF,LT}. The obstruction sheaf may be nonzero; its nonvanishing records precisely the obstructed smoothings retained in our moduli problem.

Once we have a perfect obstruction theory, general theory of the intrinsic normal cone in \cite{BF}, and \cite{Jiang_2022}
implies that there is a virtual fundamental class
$$
[M_\tau^{\Bub,\mathrm{pol}}]^{\vir}\in A_{\vd}(M_\tau^{\Bub,\mathrm{pol}}),
$$
where $\vd$ is the virtual dimension.

\begin{thm}\label{thm_virtual_class_MN}
Let $\mathbf c:M_\tau^{\Bub,\mathrm{pol}}\to M_N$ be the contraction map. If the perfect obstruction theory admits a virtual push-forward along $\mathbf c$ of the expected degree, define
$$
[M_N]^{\vir}:=\mathbf c_*[M_\tau^{\Bub,\mathrm{pol}}]^{\vir}\in A_{\vd}(M_N).
$$
\end{thm}

\subsection{Virtual class on algebraic space $M_N^{\Bub}$}

For the universal  family 
    $$
    p^{\mathrm{alg}}: \mathscr{X}^{\mathrm{alg}}\to M_N^{\mathrm{Bub}}.
$$
 in Theorem \ref{thm_good_model}, the above analysis gives 
 a perfect obstruction theory 
$$
\phi^{\mathrm{Bub}}: E^{\bullet}_{M_N^{\mathrm{Bub}}}\to \ll_{M_N^{\mathrm{Bub}}}^{\bullet}.
$$
in (\ref{eqn_OT_univ_intro_canonical}).
Hence a virtual fundamental class $[M_N^{\Bub}]^{\vir}\in H_{\vd}(M_N^{\Bub})$.  Applying the contraction map $c: M_N^{\Bub}\to M_N$ we also get a virtual fundamental class in the KSBA moduli space $M_N$. 

The issue here is that for the family $p^{\mathrm{alg}}: \mathscr{X}^{\mathrm{alg}}\to M_N^{\mathrm{Bub}}$, and each point $P\in  M_N^{\mathrm{Bub}}$, the fiber $(p^{\mathrm{alg}})^{-1}(P)$ is a representative in the $\mathbb{S}$-equivalence class of bubble-tree index-one-cover Deligne-Mumford stacks.  Thus, we have to show that the perfect obstruction theory (\ref{eqn_OT_univ_intro_canonical}) is independent to the choices of the representatives in the $\mathbb{S}$-equivalence class. This can be similarly proved in \cite[Appendix A]{Jiang_2022}.

Thus, provided that the virtual push-forward is independent of the chosen bubble-tree presentation, the construction gives a virtual fundamental class on the KSBA moduli space of surfaces of general type and an approach to Donaldson's conjecture \cite{Donaldson}.

\subsection{Tautological invariants}\label{subsec_tautological_invariants}

Assume that $M_\tau^{\Bub,\mathrm{pol}}$ and $M_N$ are proper and that the virtual push-forward of Theorem \ref{thm_virtual_class_MN} is defined. We may then use the resulting class on $M_N$ to define tautological invariants.
For the universal family $p: \mathscr{X}\to M_N$ of the KSBA moduli stack, the CM (Chow-Mumford) line bundle $L_{\CM}$ was proved to be ample in \cite{PXu15}. We define the invariant
$$
I(M_N):=\int_{[M_N]^{\vir}}(c_1(L_{\CM}))^{\vd}
$$
as the tautological invariant.
Let $K_{\mathscr{X}/M_N}$ be the relative canonical class. Similar to \cite{Alexeev_2023}, we can define
$$
\kappa_i:=p_{*}(c_1(K_{\mathscr{X}/M_N})^{i+2})\in A_i(M_N).
$$
Such $\kappa_i$ are called Kappa classes, which are very important tautological classes for the moduli space $M_N$. The first Kappa class $\kappa_1$ is the first Chern class of the CM line bundle $L_{\CM}$. Thus,
\begin{defn}\label{defn_tautological_invariants}
Let $i_1,\ldots,i_m\in \zz_{>0}$ satisfy $i_1+\cdots+i_m=\vd$. Then
we define the tautological invariants as
 $$
I_{i_1, \cdots, i_m}(M_N):=\int_{[M_N]^{\vir}}\left(\prod_{j=1}^{m}\kappa_{i_j}\right).
$$
\end{defn}

The Kappa classes were calculated for some smooth KSBA moduli space of surfaces of general type, see \cite{Alexeev_2023}.
In \cite{Donaldson}, Donaldson conjectured that the Miller-Morita-Mumford classes (characteristic classes) of the diffeomorphism group of general type surfaces can be extended to the KSBA moduli space.
It is thus interesting to study more general tautological classes like Hodge classes using the virtual fundamental class $[M_N]^{\vir}$.

More general tautological classes can also be defined for the KSBA moduli space of stable log general type surfaces \((X, D)\), where
\[
D = \sum_{i=1}^{n} a_i D_i
\]
is a \(\mathbb{Q}\)-divisor. The KSBA moduli space of stable pairs was defined in \cite{Alexeev2}, but the perfect obstruction theory on this moduli space is quite subtle. We hope to return to the construction of the virtual fundamental class for the moduli of stable pairs in future work. In cases where the KSBA moduli space of stable pairs is irreducible and has the expected dimension, the virtual fundamental class coincides with the ordinary fundamental class. This allows one to define counting invariants and study tautological invariants; see \cite{Alexeev_Ann} for a concrete construction of such moduli spaces.

\subsection{Convention}

In this paper, a {\em simple elliptic singularity} $(X,x)$ is a normal Gorenstein surface singularity such that the exceptional divisor of the
minimal resolution is a smooth elliptic curve $E$. Let $E^2=-d$ for $d\in \zz_{>0}$.
A normal Gorenstein surface singularity $(X,x)$ is called a {\em cusp} if the exceptional divisor $D$ of the minimal resolution is a cycle of smooth rational curves or a rational nodal curve.
Let $D=\cup_i D_i$, and the self-intersection sequence is given by $(-d_1, \cdots, -d_s)$, where $d_i\ge 2$, and some $d_j\ge 3$.
The cusp $D$ is an lci cusp if $\sum_{i}(d_i-2)\le 4$.
A {\em degenerate cusp} $(X,x)$ is a non-normal Gorenstein surface singularity. If $f: Y\to X$ is a minimal semi-resolution, then the exceptional divisor $D$ is a cycle of smooth rational curves or a rational nodal curve.
The $Y$ has no pinch points and the irreducible components of $Y$ have cyclic quotient singularities.

For any Deligne--Mumford stack $M$, the Chow ring $A_*(M)$ is always taking $\qq$ as coefficient.

\subsection{Outline}

Here is the outline of the paper.
In Section \ref{sec_background} we review the construction of KSBA moduli space of general type surfaces. We introduce the index-one-cover for local singularity germs, and define the index-one-cover Deligne--Mumford stacks for any slc surface $X$. We recall the moduli stack $M_N^{\ind}$ of index-one-covers and the obstruction theory $E^{\bullet}_{M_N^{\ind}}$. We explain the properties of the simple elliptic, cusp singularities that cause higher obstruction spaces. In Section \ref{sec_deg_cusp} we recall the degenerate cusp singularities and the open Kulikov models. 

We construct bubble-tree surfaces in
Section \ref{sec_flat_family_bubble_tree}. We define the relevant decorated trees and construct bubble-tree surfaces by gluing their surface components. We distinguish carefully between an admissible central fiber, a first-order pre-smoothing, and an effective one-parameter smoothing: the last exists only after the successive global obstruction classes vanish. We then introduce index-one-cover bubble-tree Deligne--Mumford stacks and study their deformation theory.

In Section \ref{sec_moduli _Bubble} we define flat families of bubble-tree index-one-cover Deligne--Mumford stacks and the polarized functor $M_\tau^{\Bub,\mathrm{pol}}$. Its objects include a stable rigidifying divisor and a first-order pre-smoothing, whose obstruction is allowed to be nonzero. We prove finite automorphisms and boundedness from the fixed data $\tau$, and state explicitly the additional effectivity and valuative hypotheses needed for a finite-type Deligne--Mumford stack and a proper contraction to the KSBA stack.

In Section \ref{sec_marked_bubble_moduli} we construct the moduli stack with labeled bubbling sections. We compare its automorphisms, boundedness, deformation theory, and virtual dimension with those of the unmarked stack, and separate this algebraic construction from the conjectural metric identification.

In Section \ref{sec_discussion} we discuss the metric tangent cones, and bubble limits for K\"ahler--Einstein metrics. We propose questions for general type surfaces.

\subsection*{Acknowledgments}
Y. J. would like to thank Professor Simon Donaldson for suggesting this beautiful idea using bubble-tree surfaces to the compactification of KSBA moduli space of general type surfaces, valuable zoom discussion, and his continuous encouragement. Y.J. thanks Professor Richard Thomas for his interest in this work, continuous discussion on the lci covering Deligne--Mumford stacks, and valuable inputs on the proof of equivalence of perfect obstruction theories.

Y. J. thanks Hans-Joachim Hein, Yuji Odaka, Song Sun and Junsheng Zhang for the correspondence and valuable discussion on
semi-log-canonical singularities, metric tangent cones and bubble limits of K\"ahler--Einstein metrics. Y. J. thanks Dori Bejleri for valuable discussion on smoothing of cusp singularities.
This work is partially supported by NSF DMS-2401484, a Simon Collaboration Grant, and the College Black-Babcock Professorship Funds from University of Kansas.

\section{Background}\label{sec_background}

In this section, we review the construction of the moduli stack of index-one-covers. This stack is isomorphic to the KSBA moduli stack of surfaces of general type.

\subsection{KSBA moduli space of general type surfaces}\label{subsec_KSBA}

We first
fix some notations. Let $X$ be a projective surface. Let $K_X$ be the canonical class of $X$, which is a Weil divisor class, and let $\omega_X$ be the dualizing sheaf. For any integer $N>0$ we set
$$\omega_X^{[N]}:=\sO_X(NK_X)=(\omega_X^{\otimes N})^{\vee\vee}.$$
\cite[Appendix to \S 1, Theorem 7]{Reid} showed that $\omega_X$ is a torsion-free sheaf of rank one. If $X$ is normal, $\omega_X$ is a divisorial sheaf which satisfies the equivalent conditions in \cite[Appendix to \S 1, Proposition 2]{Reid}.
In particular, $\omega_X$ is reflexive if $X$ is normal.

We say that $X$ has slc singularities if the following conditions hold:
\begin{enumerate}
\item the surface $X$ is reduced, Cohen-Macaulay, and has only double normal crossing singularities $(xy=0)\subset \aaa_{\mathbf{k}}^3$ away from a finite set of points;
\item Let the pair $(X^{\nu}, \Delta^{\nu})$ be the normalization of $X$ with the inverse image of the double curve. Then $(X^{\nu}, \Delta^{\nu})$ has log canonical singularities;
\item for some $N>0$ the $N$-th reflexive tensor power
$\omega_X^{[N]}$ for the dualizing sheaf $\omega_X$ is invertible.
\end{enumerate}
A stable surface $X$ is an slc surface such that $\omega_X$ is ample.

We define the KSBA moduli functor of slc surfaces by fixing $K^2, \chi, N\in\zz_{>0}$ from \cite{Kollar-Shepherd-Barron}.
Define
$$M_N:=\overline{M}_{K^2, \chi, N}: \Sch_{\cc}\to \text{Groupoids}$$
to be the moduli functor
sending
\begin{equation}\label{eqn_condition_functor_M}
T\mapsto \left\{(f:\sX\to T) \left| \begin{array}{l}
 \text{$\bullet \sX\stackrel{f}{\rightarrow} T$ is a $\qq$-Gorenstein deformation family} \\
 \text{of stable s.l.c. surfaces;} \\
 \text{$\bullet$ Conditions (1)-(5) hold for each geometric fiber;} \\
 \text{$\bullet$ For each geometric point $t\in T$, we have}\\
 \omega^{[N]}_{\sX/T}\otimes \cc(t)\to \omega_{\sX_t}^{[N]}\text{~ is an isomorphism, where}\\
 \omega^{[N]}_{\sX/T}=j_*(\omega^{\otimes N}_{\sX^0/T}), \text{and~} j: \sX^0\to \sX
 \text{~ is the inclusion}\\
 \text{of the locus where $f$ is Gorenstein.}
 \end{array} \right\}\right.
\end{equation}
modulo equivalence. The Conditions (1)-(5) are given by
 \begin{enumerate}
 \item each fiber of the family $f: \sX\to T$ is a reduced projective surface;
 \item each $\sX_t$ is connected with only s.l.c. singularities;
 \item the sheaf $\omega_{\sX_t}^{[N]}$ which is defined by $\omega_{\sX_t}^{[N]}=j_*(\omega_{(\sX_t)^0}^{\otimes N})$ and
 $j: (\sX_t)^0\to \sX_t$ is the inclusion of Gorenstein locus of $\sX_t$, is an ample line bundle;
 \item $K_{\sX_t}^{2}=\frac{1}{N^2}(\omega_{\sX_t}^{[N]}\cdot \omega_{\sX_t}^{[N]})=K^2$ for any $t\in T$;
 \item $\chi(\sO_{\sX_t})=\chi$ for $t\in T$.
 \end{enumerate}

\cite{Kollar-Shepherd-Barron}, \cite{KP15} proved that
for the fixed $K^2, \chi, N\in\zz_{>0}$, the functor $M_N$ is represented by a Deligne--Mumford stack of finite type. Suppose that $N>0$ is sufficiently large and divisible, then the stack
 $M_N$ is a proper Deligne--Mumford stack with projective coarse moduli space.

\subsection{Moduli stack of index-one-covers}\label{subsec_index_one_DM}
We introduce the index-one-cover Deligne--Mumford stack for an slc surface $X$.

Let $(X,x)$ be an slc singularity germ. The local index of $x\in X$ is the least integer $r$ such that $\omega_X^{[r]}$ is invertible around $x$. Fix an isomorphism
$\theta: \omega_X^{[r]}\to \sO_X$, we define
$$Z:=\spec_{\sO_X}\left(\sO_S\oplus\omega^{[1]}_X\oplus\cdots\oplus \omega_X^{[r-1]}\right),$$
where the multiplication on $\sO_Z$ is defined by the isomorphism $\theta$. Then $\pi:Z\to X$ is a cyclic cover of degree $r$, called the index-one cover. The inverse image of $x$ is a single point; $\pi$ is \'etale over $X\setminus\{x\}$; and $Z$ is Gorenstein, meaning that it is Cohen--Macaulay and its dualizing sheaf $\omega_Z$ is invertible. The germ $(Z,x)$ is also slc. This construction is unique locally in the \'etale topology.

A deformation $(x\in \sX)/(0\in T)$ is $\qq$-Gorenstein if it is induced by an equivariant deformation of the index-one-cover of $(x\in X)$. This means that there exists a $\mu_r$-equivariant deformation $\sZ/T$ of $Z$ whose quotient is $\sX/T$.

Since an slc surface $X$ has only finitely many singularities of local index greater than one, the local index-one-cover stacks $[Z/\mu_r]$ glue to an index-one-cover Deligne--Mumford stack $\XX\to X$ whose coarse moduli space is $X$.

Let $\sX/T$ be a flat family of slc surfaces, and let $\omega_{\sX/T}$ be the relative dualizing sheaf.
From \cite[\S 5.4]{Kollar-Shepherd-Barron}, \cite[Appdedix to \S 1]{Reid} and \cite[\S 3.1]{Hacking}, we have
$\omega_{\sX/T}^{[N]}:=(\omega_{\sX/T}^{\otimes N})^{\vee\vee}=i_*(\omega_{\sX^0/T}^{\otimes N})$,
where $i: \sX^{0}\hookrightarrow \sX$ is the inclusion of the Gorenstein locus.

Let
$\sZ/(0\in T)$ be a $\mu_r$-equivariant deformation of $Z$ inducing a $\qq$-Gorenstein deformation $\sX/(0\in T)$ of $X$,
then we have
$$\sZ=\spec_{\sO_{\sX}}(\sO_{\sX}\oplus \omega^{[1]}_{\sX/T}\oplus\cdots\oplus \omega_{\sX/T}^{[r-1]}),$$
where the multiplication of $\sO_{\sZ}$ is given by fixing a trivialization of $\omega_{\sX/T}^{[r]}$.

The index-one cover of the slc germ $(\sX,x)$ is uniquely determined in the \'etale topology. These local index-one-cover data on $\sX/T$ glue to define a Deligne--Mumford stack
$\XX/T$ which we call the index-one-cover Deligne--Mumford stack associated with $\sX/T$. The dualizing sheaf $\omega_{\XX/T}$ is invertible.

There is a one-to-one correspondence between the isomorphism classes $\{f:\sX\to T\}$ of $\qq$-Gorenstein flat families of slc surfaces and the isomorphism classes $\{f:\XX\to T\}$ of index-one-cover Deligne--Mumford stacks. In \cite{AH}, \cite[\S 5]{Jiang_2022}, the moduli functor
$$
M_N^{\ind}:=\overline{M}^{\ind}_{K^2, \chi,N}: \Sch_{\cc}\to \text{Groupoids}
$$
is defined by sending
$$
T\mapsto \{\XX\stackrel{f}{\rightarrow} T\}
$$
such that the coarse moduli space $f: \sX\to T$ of
$f: \XX\to T$ is a $\qq$-Gorenstein deformation family satisfying the condition in (\ref{eqn_condition_functor_M}).

\begin{thm}\label{thm_index_one_DM_Moduli_intro}(\cite{AH}, \cite[Theorem 5.1]{Jiang_2022})
The moduli functor $M_N^{\ind}$ of index-one-covers is represented by a Deligne--Mumford stack $M_N^{\ind}$ of finite type, and it is isomorphic to the KSBA moduli stack $M_N$.
\end{thm}

Comparing with the definition (\ref{eqn_condition_functor_M}) of Koll\'ar-Shepherd-Barron, the moduli stack $M_N^{\ind}$ is a new formulation for the KSBA moduli stack.

The advantage of the stack $M_N^{\ind}$ lies in the obstruction theory constructed in \cite{Jiang_2022}. Let
$$
p: \mathscr{X}\to M_N^{\ind}
$$
be the universal family, which is projective, flat, and relatively Gorenstein, since $M_N^{\ind}$ is a fine moduli stack.
Let $\ll^{\bullet}_{\mathscr{X}/M_N^{\ind}}$ be the relative cotangent complex of $p$ and
$\omega^{\ind}:=\omega_{\mathscr{X}/M_N^{\ind}}[2]$.
Let
$$E^{\bullet}_{M_N^{\ind}}:=Rp_{*}\left(\ll^{\bullet}_{\mathscr{X}/M_N^{\ind}}\otimes\omega^{\ind}\right)[-1],$$
where the relative dualizing sheaf $\omega_{\mathscr{X}/M_N^{\ind}}$ satisfies the property
$$\omega_{\mathscr{X}/M_N^{\ind}}|_{(p)^{-1}(t)}\cong \omega_{\XX_t}.$$
Here the dualizing sheaf $\omega_{\XX_t}$ of the index-one-cover Deligne--Mumford stack $\XX_t\to X_t$, which is locally given by $\omega_{X_t}^{[r]}$ at a singularity germ ($r$ is the index of the singular germ),
is invertible. Then the Kodaira-Spencer map $\ll^{\bullet}_{\mathscr{X}/M_N^{\ind}}\to (p)^{*}\ll_{M_N^{\ind}}^{\bullet}[1]$
induces a morphism
$$\phi: E^{\bullet}_{M_N^{\ind}}\to \ll_{M_N^{\ind}}^{\bullet}.$$
\cite[Theorem 5.6]{Jiang_2022} proved that $\phi$ is an obstruction theory in the sense of \cite{BF}.

This obstruction theory $\phi: E^{\bullet}_{M_N^{\ind}}\to \ll_{M_N^{\ind}}^{\bullet}$ is not $perfect$ in general since there exists higher obstruction spaces for the bad singularities on slc surfaces $X$ making the complex $E^{\bullet}_{M_N^{\ind}}$ not perfect.
We will review the slc singularities in the next section and classify the bad slc singularities that cause higher obstructions.

\subsection{Semi-log-canonical singularities}\label{subsec_slc_singularity}

In this section we recall the classification of Koll\'ar-Shepherd-Barron on slc singularities in \cite[Theorem 4.24]{Kollar-Shepherd-Barron}.
The slc surface singularities are exactly as follows:
\begin{enumerate}
\item the semi-log-terminal singularities;
\item the Gorenstein singularities, where every Gorenstein surface $S$ is either semi-canonical (which is smooth, normal crossing, a pinch point or a Du Val singularity),
 or has simple
elliptic singularities, cusp, or degenerate cusp singularities;
\item the $\mu_2, \mu_3, \mu_4, \mu_6$ quotients of simple elliptic singularities;
\item the $\mu_2$ quotient of cusps and degenerate cusps.
\end{enumerate}
Here the semi-log-terminal (slt) singularities are exactly:
\begin{enumerate}
\item the quotient of $\cc^2$ by Brieskorn \cite{Brieskorn};
\item normal crossing or pinch points;
\item $(xy=0)$ modulo the group action given by
$x\mapsto \zeta^a x$, $y\mapsto \zeta^b y$, and $z\mapsto \zeta z$, where $\zeta$ is a primitive
$r$-th root of unity and $(a,r)=1, (b, r)=1$;
\item $(xy=0)$ modulo the group action
$x\mapsto \zeta^a y$, $y\mapsto x$, and $z\mapsto \zeta z$, where $\zeta$ is a primitive
$r$-th root of unity and $4| r$, $(a,r)=2$;
\item $x^2=zy^2$ modulo the group action given by $x\mapsto \zeta^{1+a} x$, $y\mapsto \zeta^a y$, and $z\mapsto \zeta^2 z$, where $\zeta$ is a primitive
$r$-th root of unity and $r$ odd, and $(a,r)=1$.
\end{enumerate}

The preceding classification shows that, apart from quotient, simple elliptic, cusp, and degenerate cusp singularities and their cyclic quotients, the singularities under consideration are lci. Throughout this paper the relevant moduli functors parametrize $\qq$-Gorenstein deformations or smoothings. Thus \cite[Proposition 3.10]{Kollar-Shepherd-Barron} implies that the quotient singularities of an slc surface $X$ that admit such smoothings are class $T$ singularities
which are either rational double points or cyclic quotient singularities of the form
$$\cc^2/\mu_{r^2 s},$$
where $\mu_{r^2 s}=\langle\alpha\rangle$ and there exists a primitive $r^2 s$-th root of unity
$\eta$ such that the action is given by:
$\alpha(x,y)=(\eta x, \eta^{dsr-1}y)$
and $(d, r)=1$. When $s=1$, these are called Wahl singularities. The index-one-cover of this singularity
is $[Z/\mu_r]$, where $Z=\cc^2/\mu_{rs}$ is a type $A_{rs-1}$-singularity. Thus, $Z$ is lci and Gorenstein.
We have the following result from the above argument:

\begin{prop}\label{prop_index_one_singularity}
Let $X$ be an slc surface. Apart from its lci singularities, the index-one-cover Deligne--Mumford stack $\pi:\XX\to X$ has only simple elliptic, cusp, or degenerate cusp singularities.
\end{prop}

Let $(X,x)$ be a simple elliptic, cusp, or degenerate cusp singularity. It is known from \cite{Laufer} that if the embedding dimension $\le 4$, then it is lci. The complete local equations for simple elliptic singularities can be found in \cite[\S 6]{Looijenga_Top}. The lci cusp case can be found in \cite{Nakamura}. If the embedding dimension is $5$, although it is not lci, its higher obstruction spaces to deform this singularity vanish; see \cite{Jiang_2021}.

If the embedding dimension of $(X,x)$ is $\ge 6$, it is never lci and the higher obstruction spaces do not vanish \cite{Jiang_2021}. These are the singularities of $X\in M_N$, and the corresponding index-one-cover Deligne--Mumford stack $\XX\in M_N^{\ind}$ that cause the obstruction complex
$E^{\bullet}_{M_N^{\ind}}$ not perfect, i.e., not a two-term complex. We will construct a bubble-tree surface or index-one-cover Deligne--Mumford stack such that these bad singularities are replaced by finite bubbling rational surfaces attached together. The new bubble-tree index-one-cover Deligne--Mumford stack $\XX^{\Bub}$ only has lci singularities.

\subsection{Smoothing of simple elliptic, cusp singularities and mirror symmetry}\label{subsec_mirror_symmetry}

Simple elliptic, cusp, and degenerate cusp singularities are non-klt slc surface singularities. Their smoothing has been studied extensively. We now summarize the aspects of their smoothing theory needed here.

\subsubsection*{Simple elliptic singularity}

Let $(X, x)$ be a normal simple elliptic singularity of degree $d$. \cite{Pinkham} proved that $(X, x)$ admits a smoothing if and only if $1 \leq d \leq 9$. Now, let $Y$ be a degree $d$ del Pezzo surface and define the affine scheme $\sY := \spec\left(\bigoplus_{i \ge 0} H^0(Y, (-K_Y)^{\otimes i})\right)$.
Then $\sY$ is an affine cone and there is a smoothing
$$
f: \sY\to \aaa^1
$$
such that the central fiber is an elliptic cone which contains a simple elliptic singularity of degree $d$. Therefore, del Pezzo surfaces are closely related to simple elliptic singularities.

\subsubsection*{Cusp singularity}

Cusp singularities typically exist in dual pairs $(X, x)$ and $(X^\prime, x^\prime)$, which arise as the two singularities of an Inoue-Hirzebruch surface
$(\overline{V}, x, x^\prime)$
\cite{Inoue}, also see \cite{Looijenga}. Let $D$ and $D^\prime$ be the resolution cycles of these two cusp singularities. These cycles $D$ and $D^\prime$ consist of rational curves with sequences of negative self-intersections $d=(d_1, \cdots, d_r)$ and $d^\prime=(d_1^\prime, \cdots, d_s^\prime)$, respectively.
For a singularity to be a cusp, each sequence must satisfy $d_i \ge 2$ for all $i$, and at least one $d_j \ge 3$ to ensure the intersection matrix of $D$ is negative-definite. Furthermore, the sequence $d$ can be recovered from $d^\prime$, and vice versa. This is a consequence of the fact that the monodromy matrix of the cusp $D$ is conjugate to that of $D^\prime$.

Recall that a pair $(Y, D)$ is called a Looijenga pair (or an anticanonical pair) if $Y$ is a smooth rational surface and $D \in |-K_Y|$ is an anticanonical divisor; see \cite{Looijenga}. The smoothing of cusp singularities is intimately related to Looijenga pairs. This duality is captured by Looijenga's conjecture, which is now a theorem.

\begin{thm}\label{thm_Looijenga_conjecture}(\cite{Looijenga}, \cite{GHK15}, \cite{Engel})
 A cusp singularity $(X,x)$ admits a smoothing if and only if the dual cycle $D^\prime$ of the dual cusp $(X^\prime, x^\prime)$ lies in a rational surface $Y$ as an anticanonical divisor.
\end{thm}

Thus, the smoothing of cusp singularities is canonically related to the theory of Looijenga pairs. We now briefly review the proof.

\cite[(2.10)]{Looijenga} established the necessary condition for Theorem \ref{thm_Looijenga_conjecture}. Specifically, Looijenga proved that the Inoue surface $(\overline{V}, x, x^\prime)$ admits a smoothing for the cusp $x$, and then, by taking a simultaneous resolution, one obtains the Looijenga pair $(Y, D^\prime)$.
The sufficient condition was proved in \cite{GHK15} using log Gromov-Witten invariants for the Looijenga pair $(Y, D^\prime)$. They used genus zero relative Gromov-Witten invariants of $(Y, D^\prime)$ to construct an algebra such that the spectrum of this algebra yields a smoothing of the cusp singularity associated with $D$.

In \cite{Engel}, Engel employed an approach involving Type III degeneration of Looijenga pairs, analogous to the Type III degeneration used for K3 surfaces. From the pair $(Y, D^\prime)$, Engel constructed a Type III degeneration of Looijenga pairs, which is a normal crossing variety
\begin{equation}\label{eqn_XX_0}
 \mathscr{X}_0=\bigcup_{i=0}^{m-1}V_i
\end{equation}
where $(V_0, D, D^\prime)$ is
the Inoue surface with the resolution cycle $D, D^\prime$ of the two dual cusps; and for $i>0$, $V_i$ are rational surfaces such that $(V_i, D_i)$ is a Looijenga pair. Here we have
$$
D=D_0=\bigcup_{i}D_{0i},
$$
where $D_{0i}=V_0\cap V_i$, and for $i>0$,
$D_i=\cup D_{ij}$, $D_{ij}=V_i\cap V_j$ are double curves.
Moreover, the dual complex $\Gamma(\mathscr{X}_0)$ is a triangulation of sphere. \cite[Theorem 5.2]{Engel} proved that $\mathscr{X}_0$ admits a smoothing $f: \mathscr{X}\to \Delta$ which gives the smoothing of the cusp $D$. The generic fiber of $f$ is the Looijenga pair $(Y,D^\prime)$.
More precisely, the smoothing of the cusp $D$ is given by the smoothing $f: \mathscr{X}\to \Delta$ as follows. From \cite{Shepherd-Barron}, all the rational surfaces $V_i$ for $i>0$ in the central fiber
$f^{-1}(0)=\mathscr{X}_0$ are contractible which gives the cusp $D$.

Finally, even if a cusp $(X,x)$ is smoothable, its own resolution cycle $D$ need not be an anticanonical cycle on a rational surface; the theorem concerns the cycle of the dual cusp. See \cite{FM} for examples.

\begin{example}\label{example1}
 Let $(X,x)$ be a cusp whose resolution cycle $D$ has the negative self-intersection sequence $d=(-12, -2, -4)$, then the dual cycle $D^\prime$ has the negative self-intersection sequence
 $d^\prime=(-4, \underbrace{-2, \cdots, -2}_{9}, -3, -2)$. The cusp $D^\prime$ is a hypersurface cusp which admits a smoothing, but $D^\prime$ cannot be an anticanonical divisor of a Looijenga pair since from \cite[Theorem (4.5)]{FM}, $d=(-12, -2, -4)$ does not have a rational dual.
\end{example}

\begin{example}\label{example2}
As Professor Friedman pointed out,
 the resolution cycles for ``most" smoothable cusp singularities are not realized as anticanonical cycles on a rational surface. For instance, any cusp singularity of multiplicity at most 3 is a hypersurface singularity, hence smoothable. Let the
 self-intersection sequence be $(-3, \underbrace{-2, -2, \cdots, -2)}_{r}$ with $r$ at least $11$, then this resolution cycle is not realized as an anticanonical cycle on a rational surface.
\end{example}

\section{Degenerate cusp singularities}\label{sec_deg_cusp}

Let $(X,x)$ be a degenerate cusp singularity with resolution cycle $D$. In this case $X$ is non-normal, and if
$$
\sigma: \widetilde{X}\to X
$$
is the normalization, then $\widetilde{X}$ is a disjoint union of surfaces $(X_i, E_i)$ where $E_i$ are double curves and
$(X_i, x)$ have cyclic quotient $\cc^2/\mu_{r_i}$-singularities which are log-terminal singularities.
It is well-known that degenerate cusp singularities always admit a smoothing.
The smoothing of degenerate cusp singularities also has mirror symmetry properties.

\subsection{Open Kulikov surfaces}\label{subsec_open_Kulikov}

We review the open Kulikov surfaces introduced in \cite{AAB2024}.

An log Calabi-Yau surface $(Y,D)$ is a log canonical pair $(Y,D)$, where $D$ is a reduced effective divisor with $K_Y+D$ Cartier and trivial. We always assume $D$ is non-empty and singular. From \cite[Proposition 2.1]{AAB2024}, away from $D$, the surface $Y$ has at worst du Val (ADE) singularities, and the germ of $(Y,D)$ at any point $p\in D$ has cyclic quotient singularities.

Any log Calabi-Yau surface $(Y,D)$ has a toric model $(\overline{Y}, \overline{D})$, where $\overline{Y}$ is a toric surface and $\overline{D}$ is the toric divisor. The log Calabi-Yau surface $(Y,D)$ can be obtained from $(\overline{Y}, \overline{D})$ by a sequence of interior blow-ups and corner blow-ups.

\begin{defn}\label{defn_open_logCY}
 An {\em open} log Calabi-Yau surface $(Y,D)$ is a smooth complex surface $Y$, together with a non-empty reduced normal crossing divisor $D\subset Y$, such that there exists a smooth log Calabi-Yau surface $(\widetilde{Y},\widetilde{D})$ and a non-empty connected curve $E$ in $\widetilde{Y}$, given by a union of irreducible components of $\widetilde{D}$, such that $Y=\widetilde{Y}\setminus E$, and $D=\widetilde{D}\setminus E$.
\end{defn}

\begin{rmk}
 The divisor $\widetilde{D}$ has at least two irreducible components, which is a cycle of projective lines. Therefore, $D=D_0+\cdots +D_{s}$, where
 $D_0=D_s=\aaa^1$, and $D_i\cong \pp^1$ for $1\le i\le s-1$.

 Any open log Calabi-Yau surface $(Y,D)$ can be completed into a smooth log Calabi-Yau surface.
\end{rmk}

\begin{defn}\label{defn_0_surface}
 A {\em $0$-surface} is an open log Calabi-Yau surface $(Y,D)$ obtained by interior and corner blow-ups of an affine toric surface $(\overline{Y}, \overline{D})$, where the interior blow-ups are taken along the torus fixed points of $\overline{Y}$.
\end{defn}

\begin{defn}\label{defn_1_surface}
 A {\em $1$-surface} is an open log Calabi-Yau surface $(Y,D)$ such that there exists a projective surjective morphism $Y\to \aaa^1$, which contracts the compact irreducible components of $D$ to $0\in \aaa^1$.
\end{defn}

\begin{defn}\label{defn_2_surface}
 A {\em $2$-surface} is a smooth log Calabi-Yau surface $(Y,D)$.
\end{defn}

An open Kulikov surface is a reduced normal crossing surface
\begin{equation}\label{eqn_XX_open_Kulikov}
 \XX_0=\bigcup_{i=0}^{m-1}X_i
\end{equation}
satisfying the following conditions:
\begin{enumerate}
\item The dual complex of $\XX_0$ is a disk with a triangulation $\sT$.
 \item The irreducible components of $\XX_0$ are labeled by $I=\{0, 1,\cdots, m-1\}$. For each $X_i$, the double curve
 $$D_i:=\bigcup_{j\in I\setminus \{i\}}X_i\cap X_j.$$
 Let $\sigma: \widetilde{X}_i\to X_i$ be the normalization, and $\widetilde{D}_i=\sigma^{-1}(D_i)$. Then for each interior vertex $v_i\in \text{Int}(\sT)$, the pair
 $(\widetilde{X}_i, \widetilde{D}_i)$ is a $2$-surface, and for each boundary vertex $v_j\in \partial \sT$, the pair $(\widetilde{X}_i, \widetilde{D}_i)$ is either a $0$-surface, or a $1$-surface.
 \item For every irreducible component $D_{ij}=X_i\cap X_j$,
 $$
 (D_{ij}|_{X_i})^2+ (D_{ij}|_{X_j})^2=-2+2g(D_{ij}),
 $$
 where $d(D_{ij})$ is the arithmetic genus of $D_{ij}$.
\end{enumerate}

\subsection{Smoothing of degenerate cusp singularities}\label{subsec_deg_cusp_smoothing}

Let $(X,0)$ be a degenerate cusp singularity. Then it is known that it always admits a smoothing $\overline{f}: \overline{\sX}\to \Delta$.

\begin{lem}
 There exists a crepant resolution $f: \sX\to \Delta$ for $\overline{f}: \overline{\sX}\to \Delta$.
\end{lem}
\begin{proof}
 The proof is similar to \cite[Theorem (2.2)]{Reid}. Let
 $\pi: \tilde{\sX}\to \overline{\sX}$ be some resolution, then Hironaka's resolution theorem implies that every crepant resolution $f: \sX\to \overline{\sX}$ can be obtained from the diagram
 \[
 \xymatrix{
 &\widetilde{\widetilde{\sX}}\ar[dl]\ar[dr]&\\
 \widetilde{\sX}\ar[dr]&& \sX\ar[dl]\\
 &\overline{\sX}&
 }
 \]
 where $\widetilde{\widetilde{\sX}}\to \widetilde{\sX}$ is a blow-up. The crepant prime divisors of $\sX$ are mapped injectively to those of the fixed $\widetilde{\sX}$.
\end{proof}
%\marginpar{Fix me. The proof of lemma}

Suppose that the central fiber of $f:\sX\to\Delta$ is an open Kulikov surface. We recall one of the main results of \cite{AAB2024}.

\begin{figure}[htbp]
 \centering
 \IfFileExists{open.pdf}{\includegraphics[width=35mm]{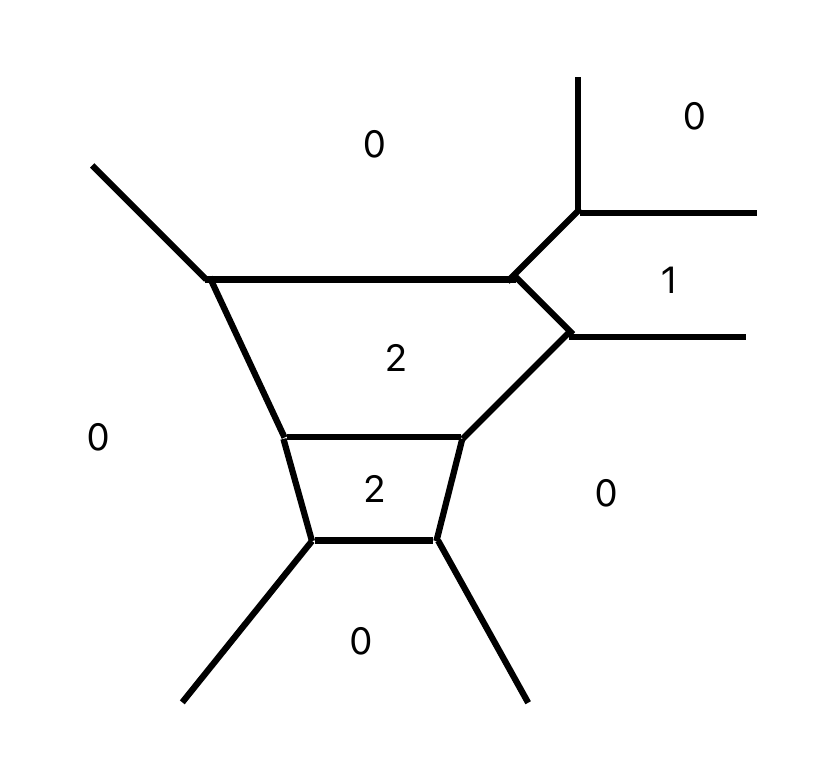}}{\IfFileExists{open.png}{\includegraphics[width=35mm]{open.png}}{\fbox{\parbox{35mm}{\centering Original figure \\ \texttt{open} not supplied}}}}
 \caption{Open Kulikov surface $\sX_0$}
 \label{fig-open_Kulikov}
 \centering
 \IfFileExists{central.pdf}{\includegraphics[width=30mm]{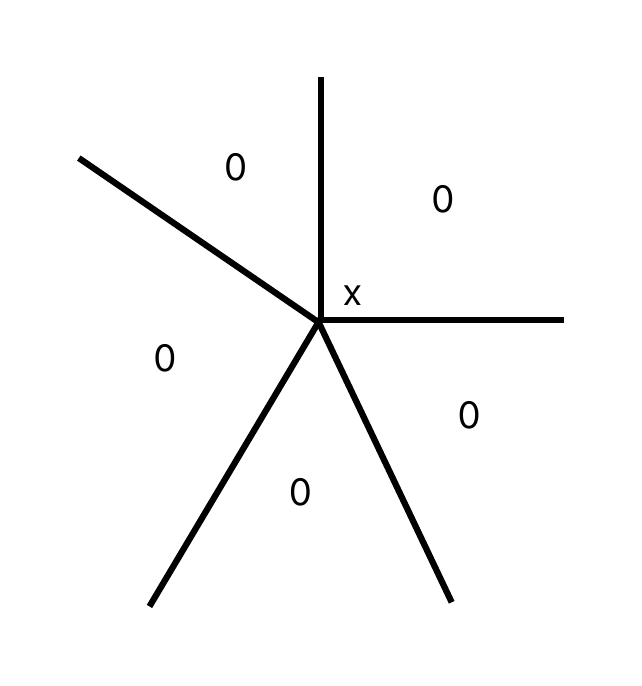}}{\IfFileExists{central.png}{\includegraphics[width=30mm]{central.png}}{\fbox{\parbox{30mm}{\centering Original figure \\ \texttt{central} not supplied}}}}
 \caption{degenerate cusp singularity $(\overline{\sX}_0,x)$}
 \label{fig-open_Kulikov_2}
\end{figure}

\begin{thm}\label{thm_semistable_degeneration}
Let $(Y,D,L)$ be a polarized log Calabi-Yau surface, where $L$ is an ample line bundle. Then there is a semistable degeneration
$$
f: \sX\to \Delta
$$
such that the central fiber $f^{-1}(0)=\sX_0$ is an open Kulikov surface.
Moreover, its affinization
$$
\overline{f}: \overline{\sX}\to \Delta
$$
is a smoothing of a degenerate cusp singularity determined by $D$.
\end{thm}
\begin{proof}
 The proof of this statement is from \cite[\S 4.2]{AAB2024}. We briefly review it here.

 Let $(\overline{Y}, \overline{D}, \overline{L})$ be a toric model of $(Y,D,L)$, and $\overline{P}$ be the toric momentum polytope. We perform the Symington process for the polytope $\overline{P}$ to get the Symington polytope $P$.
 Let $\mathscr{P}$ be the polyhedral decomposition induced by the polyhedral decomposition $\overline{\mathscr{P}}$ of $\overline{P}$.

Let $M\cong \zz^2$ be the character lattice for $\overline{Y}$, so that $\overline{P}\subseteq M_{\rr}$ where $M_{\rr}=M\otimes_{\zz}\rr$. Let
$$
\overline{f}: \sX_{\overline{\mathscr{P}}}\to \aaa^1
$$
be the Mumford toric degeneration, where $\sX_{\overline{\mathscr{P}}}$ is the $3$-dimensional toric variety with fan $C(\overline{\mathscr{P}})$ in
$M_{\rr}\times \rr_{\ge 0}$ consisting of the cones
$$
C(\sigma)=\rr_{\ge 0}\cdot (\sigma\times \{1\})\subset M_{\rr}\times \rr_{\ge 0}
$$
over the cells $\sigma\in \overline{\mathscr{P}}$. The central fiber $\sX_{\overline{\mathscr{P}},0}:=\overline{f}^{-1}(0)$ is a union of irreducible components
$X_{\overline{\mathscr{P}},v}$, indexed by the vertices
$v\in \overline{\mathscr{P}}^{[0]}$. We have
$$\sX_{\overline{\mathscr{P}},0}=\bigcup_{v\in \overline{\mathscr{P}}^{[0]}}X_{\overline{\mathscr{P}},v}$$
where all the $X_{\overline{\mathscr{P}},v}$ are toric surfaces with fan the star of $v$ in $\overline{\mathscr{P}}$.

The open Kulikov surface $\sX_{\mathscr{P},0}$ is constructed from a deformation of $\sX_{\overline{\mathscr{P}},0}$. We use $X_{\mathscr{P},v}$ indexed by $v\in \mathscr{P}^{[0]}$ to represent the irreducible components of $\sX_{\mathscr{P},0}$, which are constructed as follow.
Let
$$
\pr: \overline{\mathscr{P}}^{[0]}\to \mathscr{P}^{[0]}
$$
be the natural projection from the Symington process, since the Symington polytope $P$ is given by the quotient
$\overline{P}/\cup_{ij}\text{Int}(\Delta_{ij})$.
For any vertex $v\in \mathscr{P}^{[0]}$, there is a toric surface $\widetilde{X}_{\mathscr{P},v}$ which is a deformation of the union
$$
\bigcup_{v^\prime\in \pr^{-1}(v)}X_{\overline{\mathscr{P}},v},
$$
where each toric surface is an irreducible component in $\sX_{\overline{\mathscr{P}},0}$. From the toric pair $(\widetilde{X}_{\mathscr{P},v}, \partial \widetilde{X}_{\mathscr{P},v})$, \cite[\S 4.2]{AAB2024} obtained a deformation pair $(X_{\mathscr{P},v}, \partial X_{\mathscr{P},v})$ which is a log Calabi-Yau surface.

Then gluing the log Calabi-Yau surfaces $X_{\mathscr{P},v}$ for $v\in \mathscr{P}^{[0]}$ along the divisors corresponding to the edges in $\mathscr{P}$, and we get the open Kulikov surface $\sX_{\mathscr{P},0}$.

Since the dual complex of $\sX_{\mathscr{P},0}$ is a disk (the polytope $P$), the polyhedral decomposition implies that if the vertex $v\in \mathscr{P}^{[0]}$ is an interior vertex, then $(X_{\mathscr{P},v}, \partial X_{\mathscr{P},v})$ is a $2$-surface; and if the vertex $v\in \mathscr{P}^{[0]}$ is a boundary vertex, $(X_{\mathscr{P},v}, \partial X_{\mathscr{P},v})$ is a $0$, or $1$-surface.

The open Kulikov surface $\sX_{\mathscr{P},0}$ is quasi-projective, see \cite[Lemma 4.7]{AAB2024}, and admits a smoothing
$$
f: \sX_{\mathscr{P}}\to \Delta.
$$
This smoothing induces a smoothing
$$
\overline{f}: \overline{\sX}_{\mathscr{P}}\to \Delta
$$
so that $\overline{\sX}_{\mathscr{P}}$ is the affinization of
$\sX_{\mathscr{P}}$. The $2$-surfaces and the boundary of the $1$-surfaces contract to the singular point $x\in \overline{\sX}_{\mathscr{P}}$ so that the central fiber
$(\overline{\sX}_{\mathscr{P},0},x)$ is a degenerate cusp singularity.
We explain this contraction in Figure \ref{fig-open_Kulikov}
 and Figure \ref{fig-open_Kulikov_2}.
\end{proof}

\begin{rmk}
 From \cite{AAB2024}, the semistable degeneration $f: \sX_{\mathscr{P}}\to \Delta$ is a crepant resolution of the degenerate cusp singularity $(\overline{\sX}_{\mathscr{P},0},x)$. Two crepant resolutions $f_1: \sX^1_{\mathscr{P}}\to \Delta$
 and $f_2: \sX^2_{\mathscr{P}}\to \Delta$ are related by a birational modification.
\end{rmk}

\section{Smoothing of bubble-tree surfaces}\label{sec_flat_family_bubble_tree}

In this section, we define bubble-tree surfaces and study their smoothing. Before giving the definition, we first review the theory of integral affine surfaces. These surfaces will form the dual complexes of the attached bubble-tree surfaces for cusp singularities.

\subsection{Integral affine surfaces}\label{subsec_integral_affine}

A basis triangle in $\rr^2$ is a triangle of area $1/2$ whose vertices lie in $\zz^2$. At each vertex, the two primitive integral edge vectors form a basis of $\zz^2$.

\begin{defn}\label{defn_triangle_integral-affine}(\cite[Definition 3.1]{Engel})
A triangulated integral-affine surface, possibly with singularities, is a triangulated real surface $S$, possibly with boundary, such that:
\begin{enumerate}
\item the complement of the vertices $\{\alpha_i\}\subset S$ admits an atlas of charts into $\rr^2$ whose transition functions lie in
$SL_2(\zz)\ltimes \zz^2$.
\item the interior of every triangle admits a chart to a basis triangle.
\end{enumerate}
\end{defn}

An integral-affine surface has a canonical orientation induced from the standard orientation on $\rr^2$. Let $e_{ij}$ be the edge $\alpha_i-\alpha_j$ in the triangulation of $S$. Let $f_{ijk}$ be
the triangle whose counterclockwise vertices are $\alpha_i, \alpha_j, \alpha_k$. In this chart we can write $e_{ij}=\alpha_j-\alpha_i$.

Let $(Y,D^\prime)$ be a Looijenga pair. The pseudo-fan of $(Y,D^\prime)$ is a triangulated integral-affine surface whose underlying surface $S_{(Y,D^\prime)}$ is the cone over the dual complex of $D^\prime$, which is a polygon when $D^\prime$ is a cycle of rational curves.

Let $e_i$ be the edge from the cone point to the vertex corresponding to $D^\prime_i$. Then the negative self-intersection of $e_i$ is:
$$
d^{\prime}_i=
\begin{cases}
-D^{\prime 2}_i, & n>1;\\
2-D^{\prime 2}_i, & n=1.
\end{cases}
$$

Starting from a Looijenga pair $(Y, D^\prime)$, and its integral affine surface $S_{(Y, D^\prime)}$, \cite{Engel} constructed a Type III canonical degeneration $\XX_0$ in (\ref{eqn_XX_0}), such that the dual complex $\Gamma(\XX_0)$ is a triangulation of the sphere $\mathbb{S}^2$.
The vertices $\{\alpha_i\}$ correspond to the rational components $V_i$, the directed edges $e_{ij}$ correspond to double curves $D^{\prime}_{ij}$, and triangular faces $f_{ijk}$ correspond to
triple points such that
\begin{equation}\label{eqn_triple_point}
d^{\prime}_{ij}:=
\begin{cases}
-D^{\prime 2}_{ij}, & \ell(D^{\prime}_i)\ge 2;\\
2-D^{\prime 2}_{ij}, & \ell(D^{\prime}_i)=1
\end{cases}
\end{equation}
where $d^{\prime}_{ij}$ is the negative self-intersection of $e_{ij}$.

Actually the smoothing $f: \XX\to \Delta$ of the Type III degeneration $\XX_0$ in (\ref{eqn_XX_0}) is proved using the triple point formula (\ref{eqn_triple_point}).
The dual complex $\Gamma(\XX_0)$ carries the singular integral-affine structure produced by Engel's construction.

\subsection{Bubble-tree surfaces}\label{subsec_bubble_tree}

We begin with the combinatorial notion of a rooted tree.

\begin{defn}\label{defn_tree}
 A tree $\sT$ is a finite graph, oriented by a partial order $<$ such that
 \begin{enumerate}
 \item there exists a unique vertex $\alpha_{\rt}\in\sT$, called the root of $\sT$.
 \item for any vertex $\alpha\in \sT$, $\alpha\neq \alpha_{\rt}$, there exists a unique maximal $\beta<\alpha$, the predecessor of $\alpha$, which is denoted by $\alpha^{-}$.
 \item let $\alpha^{+}:=\{\beta\in \sT| \beta^{-}=\alpha\}$ be the set of direct successors of $\alpha\in \sT$. Also let $\sT_{\top}$ be the vertices of $\sT$ without successors (i.e., the top vertices of $\sT$).
 \end{enumerate}
\end{defn}

\begin{rmk}
Because $\sT$ is finite, it always has a terminal vertex. We allow the trivial tree consisting only of its root; otherwise every nonterminal vertex has a nonempty set of successors.
\end{rmk}

We define the bubble-tree surface using trees $\sT$. We fix a surface $X$ with singularities
$$\underbrace{p_1, \cdots, p_s}_{\text{simple elliptic}}, \underbrace{p_{s+1},\cdots, p_{t}}_{\text{cusp}}, \underbrace{p_{t+1},\cdots, p_n}_{\text{deg cusp}}.$$
We assume that these singularities are not lci. Some simple elliptic singularities may carry an action of $\mu_r$, with $r\in\{2,3,4,6\}$, whose quotient is rational. A cusp or degenerate cusp may likewise carry an involution with rational quotient.
We construct the attached log Calabi-Yau surfaces case by case.

\subsection*{Step 1:}

For the simple elliptic singularities $p_i$ of degree $d_i$ for $1\le i\le s$,
we take $\sT_i$ as a chain tree, which means $|\beta^+|=1$ for any $\beta\in \sT_i$.
Let $\alpha_i=\alpha_{\rt}$ be the unique root in the tree $\sT_i$.
For each nonterminal vertex $\beta$, let $X_\beta$ be an elliptic ruled surface, with the two distinguished sections used as its incoming and outgoing double curves. For the unique terminal vertex $\beta_{\top}$, let $X_{\beta_{\top}}=Y_i$ be a smooth del Pezzo surface of degree $d_i$, and choose a smooth anticanonical elliptic curve $E_i\in|-K_{Y_i}|$. Notice that the intermediate ruled surfaces are not rational.

We let
$$
X_{\sT_i}=\bigcup_{\beta\in \sT_i}X_{\beta}.
$$
be the surface obtained by gluing the above surfaces $X_{\beta}$ along the elliptic curves. The surface is a normal crossing surface.

Let $\pi_i:\widetilde X_i\to X$ be the minimal resolution of $p_i$, with exceptional elliptic curve $E_i^{\mathrm{exc}}$. The normal bundle has degree $(E_i^{\mathrm{exc}})^2=-d_i$. Choose the incoming section of the first ruled component with the dual normal bundle, and glue it to $E_i^{\mathrm{exc}}$. Successive components are glued along their distinguished sections, with dual normal bundles, and the final ruled component is glued to the anticanonical curve in $Y_i$. We obtain the simple-normal-crossing surface
\begin{equation}\label{eqn_simple_i}
 X^{\Bub}_{\sT_i}:=\widetilde X_i\mathop{\cup}_{E_i^{\mathrm{exc}}}X_{\sT_i}.
\end{equation}

This resembles a Type II degeneration of K3 surfaces, although the intermediate components here are elliptic ruled surfaces rather than rational surfaces.
From a construction in \cite{FHJ}, the affine surface $Y_i\setminus E_i$ is the bubble limit of the rescaling K\"ahler--Einstein metrics
for simple elliptic singularities.

If the simple elliptic singularity $p_i$ admits an action of a cyclic group $\mu_r$, for $r\in\{2,3,4,6\}$, an equivariant construction requires compatible lifts to every ruled and del Pezzo component and to all gluing data. Subject to this condition, the corresponding index-one-cover model is the quotient stack of the chain by $\mu_r$.

\subsection*{Step 2:}\label{step2}

For a cusp $p_j$, with $s+1\le j\le t$, write the self-intersection sequence of its resolution cycle $D_j$ as $(-d_j^1,\ldots,-d_j^{r_j})$. Write the corresponding sequence for the dual cycle $D_j^\prime$ as $(-d_j^{\prime 1},\ldots,-d_j^{\prime s_j})$. The two cycles are dual in the sense of \cite{Nakamura}.

For a cusp singularity $p_j$ that admits a smoothing $f: \XX \to \Delta$, the resolution cycle $D_j$ may not be an anticanonical divisor of a rational surface $Y_j$, as shown in Example \ref{example1} and Example \ref{example2}. Therefore there are two cases.
\begin{enumerate}
 \item Case one is that the resolution cycle $D_j$ of the cusp singularity $p_j$ is not an anticanonical divisor of a rational surface $Y_j$.
 \item Case two is that the resolution cycle $D_j$ of the cusp singularity $p_j$ is an anticanonical divisor of a rational surface. In this case, the existence of the Looijenga pair $(Y, D_j)$ implies, by the Looijenga conjecture, that for the dual cusp $p_j^\prime$, not only is the resolution cycle $D_j^\prime$ an anticanonical divisor of a rational surface, but the singularity $p_j^\prime$ also admits a smoothing.
\end{enumerate}

In both these cases, we construct a new tree $\mathcal{T}_j$ whose root $\alpha_{\rt}$ is $\alpha_j$; and the construction assumes that the cusp $p_j$ admits a smoothing.
By Theorem \ref{thm_Looijenga_conjecture}, the dual cycle $D_j^\prime$ is an anticanonical cycle of a rational surface $Y_j$, so that $(Y_j, D_j^\prime)$ forms a Looijenga pair.

For the tree $\mathcal{T}_j$, we allow the vertices $\beta_k$ to form triangles. We label the vertices as
$\{0, 1, \cdots, m-1\}$ so that $\alpha_j=\alpha_0$ and other vertices are $\beta_1, \cdots, \beta_{m-1}$.
Starting from the Looijenga pair $(Y_j, D_j^\prime)$, the construction yields a triangulated integral affine sphere $S_j \cong \mathbb{S}^2$ associated with the tree $\mathcal{T}_j$, as described in equation (\ref{eqn_XX_0}).

The method for constructing a Type III degeneration, as in (\ref{eqn_XX_0}), involves the Symington process applied to the polytope associated with the toric model of the Looijenga pair $(Y_j, D_j^\prime)$. This polytope is then completed into an integral affine sphere via triangulation. Further details can be found in \cite[Section 5]{Engel}.

Thus we have a Type III central fiber
$$
X_{\sT_j}:=\bigcup_{k=0}^{m-1}X_{k}
$$
where for the root $\alpha_j=\alpha_{0}$, we have $X_{\alpha_0}=(V_j, D_j, D_j^\prime)$ is an Inoue surface in Case one (in Case two we may choose $X_{\alpha_0}=Y$). The other $X_{k}$ are rational surfaces for other $\beta_k\in \sT_j$ so that $(X_{k}, D_{k})$ are log Calabi-Yau surfaces.
Here we have
$$
D_j=D_{\alpha_0}=\cup_{s}D_{0s}:=\cup_{s=1}^{m-1} X_{\alpha_0}\cap D_{s}
$$
and
$$D_{k}=D_{k0}\cup_{s} D_{ks}=(X_{k}\cap X_{\alpha_0})\cup \cup_{s}(X_{k}\cap X_{s})$$
for $1\le k\le m-1$. The dual complex of $X_{\sT_j}$ is the integral affine sphere $\mathbb{S}^2$.

For the cusp singularity $(X, p_j)$, we let $\pi: \widetilde{X}\to X$ be the minimal resolution of the cusp singularity $p_j$ with resolution cycle $D_j$. Let $U_j\subset \widetilde{X}$ be an open neighborhood of $D_j\in \widetilde{X}$ which is the same as the open neighborhood of $D_j$ in the Inoue surface $V_j$.
Let $X_{\sT_j}^\prime:=\bigcup_{k=1}^{m-1}X_k$. Its boundary contains a cycle canonically identified, with its normal-bundle data, with $D_j$. Replace the analytic Inoue component $X_{\alpha_0}$ by $\widetilde X$ and glue along this cycle. This is a pushout along the divisor $D_j$, not a gluing of punctured neighborhoods.
This construction is similar to the Type III degeneration of K3 surfaces. Thus, we have
\begin{equation}\label{eqn_cusp_j}
 X^{\Bub}_{\sT_j}:=\widetilde X\mathop{\cup}_{D_j}X^\prime_{\sT_j}.
\end{equation}

If $p_j$ carries an involution with rational quotient, an equivariant bubble-tree construction additionally requires the involution to extend to the chosen Looijenga pair and its Type III model. Such an extension is a hypothesis, not an automatic consequence of the nonequivariant construction. Related equivariant Type III degenerations of K3 surfaces are studied in \cite{AET23}.

\begin{rmk}
 In the above construction, we have to use the analytic neighborhood $U_j$ of the cusp singularity in the Inoue surface $V_j$. It is interesting to have a pure algebraic construction.
\end{rmk}

\begin{example}
 Let $(X, p_j)$ be a cusp singularity with resolution cycle $D_j$ consisting of two rational curves and
 $d_j=(-6,-9)$. The dual cycle $D_j^\prime$ is
 $d_j^\prime=(-3,-2,-2,-2,-3,-2,-2,-2,-2,-2,-2)$. Assuming that $(X,p_j)$ is smoothable, Looijenga's theorem gives a Looijenga pair $(Y,D_j^\prime)$. According to \cite[Example 5.9]{Engel}, there are two deformation families of such pairs, associated with two toric models. One Type III central fiber has the form
 $$
 X_{\sT_j}=V_0\cup V_1\cup V_2
 $$
 where $(V_0, D_j, D_j^\prime)$ is the Inoue surface, where
 $D_j=D_{01}\cup D_{02}$, whose self-intersection sequence is given by $(-6,-9)$. Here $D_1=D_{10}\cup D_{12}$ with $D_{10}^2=4$, $D_{12}^2=1$; and
 $D_2=D_{20}\cup D_{21}$ with $D_{20}^2=7$, $D_{21}^2=-3$.

 So if there is a smoothing $f: \sX\to \Delta$ of the singularity $(X, p_j)$, we take minimal resolution
 $\widetilde{X}\to X$ with resolution cycle $D_j$ whose self-intersection sequence is $(-6,-9)$. Then we glue along this $D_j$ and $D_1, D_2$ with $V_1\cup V_2$.
\end{example}

\begin{rmk}
 It is interesting to find concrete examples for the moduli space $M_N^{\Bub}$ so that the cusp singularities are replaced by a Type III degeneration of Looijenga pairs. Since the Type III degeneration of log Calabi-Yau surfaces is similar to the Type III degeneration of K3 surfaces in \cite{AE23}, it is interesting to compare the moduli compactification of K3 surfaces to the bubble-tree surfaces defined below.

 In \cite{ABJ}, the bubble moduli stack $M_N^{\Bub}$ and the moduli stack of lci covers for fibered surfaces as in \cite{Jiang_2022} are related by wall crossings.
\end{rmk}

\subsection*{Step 3:}

For the degenerate cusp singularities $(X,p_k)$ with $t+1 \le k \le n$, let $D_k$ be the resolution cycle, with its sequence of negative self-intersections given by $d_k = (-d_k^1, \dots, -d_k^{r_k})$.
By the definition of a degenerate cusp singularity, the normalization
$$
\sigma: \widetilde{X}\to (X,p_k)
$$
of this singularity germ has only cyclic quotient singularities.

We assume that $d_k^s\ge2$ for every $s$ and that $d_k^t\ge3$ for at least one index $t$. Any assertion about the quotient types on the normalization or the vanishing of higher obstruction spaces must be checked from the local classification; it does not follow from these inequalities alone.

Suppose that there is a smoothing
$$
\overline{f}: \overline{\sX}\to \Delta
$$
of this degenerate cusp singularity. Assume that, after a finite base change if necessary, it admits a crepant semistable model
\[
f:\sX\to\Delta
\]
from which the bubble surface is constructed.

We now apply Theorem \ref{thm_semistable_degeneration} to this degenerate cusp.
The construction in \S \ref{sec_deg_cusp}, and \cite{AAB2024} imply that there is a dual cycle $D_k^\prime$ of the cycle $D_k$ such that $D^\prime_k$ is an anticanonical divisor of a smooth surface $Y$. Thus, we get a log Calabi-Yau surface $(Y,D_k^\prime,L_k)$ after fixing a polarization $L_k$. Then from Theorem \ref{thm_semistable_degeneration}, there is a smoothing of an open Kulikov surface
$$
f: \sX\to \Delta
$$
such that its affinization is $\overline{f}: \overline{\sX}\to \Delta$.

We write the open Kulikov surface as
$$
\sX_0=
\bigcup_{v\in\sT_k^{[0]}}X_{v},
$$
where $\sT_k$ is the triangulated complex of the disk.
This is the polyhedral decomposition $\mathscr{P}$ in \S \ref{sec_deg_cusp}.
The number of $0$-surfaces in the open Kulikov surface
$f^{-1}(0)=\sX_0$ is the same as the number of components in the normalization $\widetilde{X}$.

To obtain a genuine normal-crossing pushout, one must specify an identification between the conductor curves on the semiresolution of $(X,p_k)$ and the boundary curves of the $0$-surfaces, together with compatible normal-bundle data. Assuming these data satisfy Definition \ref{defn_admissible_bubble_tree}, we define
\begin{equation}\label{eqn_deg_cusp_k}
 X^{\Bub}_{\sT_k}:=\widetilde X\mathop{\cup}_{D_k}\sX_0,
\end{equation}
where the notation denotes this conductor pushout. Boundary components not meeting the semiresolution are completed by log Calabi--Yau components as in Definition \ref{defn_open_logCY}.

We conjecture in Section \ref{sec_discussion} that this algebraic model records the metric bubble tree. Sun's non-collapsing finiteness theorem \cite{Sun} motivates that expectation but does not prove it in the present collapsing, non-normal setting.

For an equivariant version, we additionally require a lift of the $\mu_2$-action to the chosen open Kulikov model and to the conductor gluing. When such a lift exists, it acts compatibly on the components of $D_k=\bigcup_{s=1}^{r_k}D_k^s$ and hence on $X^{\Bub}_{\sT_k}$.

\begin{rmk}\label{rmk_bubble_sun}
Corollary 1.13 of \cite{Sun} proves finiteness of the bubble tree under its non-collapsing hypotheses. We do not use that theorem to construct the algebraic surface above; rather, it motivates the expectation that a finite combinatorial model should exist. Relating the two constructions is part of the conjectural differential-geometric picture; see also \cite{Odaka}.
\end{rmk}

Finally, from (\ref{eqn_simple_i}), (\ref{eqn_cusp_j}) and (\ref{eqn_deg_cusp_k}), we define the bubble-tree surfaces.

\begin{defn}\label{defn_admissible_bubble_tree}
A bubble-tree datum is \emph{admissible} if the following conditions hold.
\begin{enumerate}
\item Every component is smooth and the boundary on each component has simple normal crossings; each irreducible double curve is smooth; the resulting pushout is a simple-normal-crossing surface (with the usual triple-point local model where three components meet).
\item For every component $V$, the reduced double divisor $D_V$ satisfies $K_V+D_V\sim0$ on a bubble component. On the resolved main component it is the log pullback of the canonical class of $X$.
\item For every double curve $C=V_i\cap V_j$, the infinitesimal normal bundle satisfies the triple-point formula
\[
N_{C/V_i}\otimes N_{C/V_j}\otimes\sO_C(T_C)\cong\sO_C,
\]
where $T_C$ is the divisor of triple points on $C$. In addition, these trivializations must agree at triple points and define a global trivialization of $T^1$ on the singular locus. Componentwise degree identities alone do not imply this global condition.
\item The identifications of double curves satisfy the cocycle condition at triple points. If a group acts, all gluing data and normal-bundle identifications are equivariant.
\item A polarization or boundary marking is chosen so that the automorphism group of the resulting object is finite.
\end{enumerate}
\end{defn}

\begin{defn}\label{defn_bubble_tree_surface}
 Let $X$ be a reduced surface (not necessarily projective) with singularities $p_1,\ldots,p_n$ of the types considered above, and choose an admissible bubble-tree datum at each $p_i$.
 The associated bubble-tree surface $X^{\Bub}$ is the pushout obtained by replacing the germ at each $p_i$ with the corresponding resolved main component and attached bubble components. The gluing is the divisor gluing specified in (\ref{eqn_simple_i}), (\ref{eqn_cusp_j}), or (\ref{eqn_deg_cusp_k}); it is not a set-theoretic union with $X$.
\end{defn}

\begin{prop}
Let $X$ be a stable slc surface and let $X^{\Bub}$ arise from an admissible bubble-tree datum. Then $X^{\Bub}$ is slc and, for each bubble component $V$ with reduced double divisor $D_V$,
\[
\omega_{X^{\Bub}}|_V\cong\omega_V(D_V)\cong\sO_V.
\]
The automorphism group preserving the polarization or boundary marking in Definition \ref{defn_admissible_bubble_tree} is finite.
\end{prop}
\begin{proof}
The local models in Definition \ref{defn_admissible_bubble_tree} are simple normal crossings, hence demi-normal, and the log-canonical condition follows componentwise from the normalization criterion for slc surfaces. The adjunction formula for a normal-crossing union gives $\omega_{X^{\Bub}}|_V\cong\omega_V(D_V)$, which is trivial by admissibility. Finiteness is part of the stability condition imposed by the polarization or boundary marking; it would be false for arbitrary unmarked log Calabi--Yau components.
\end{proof}

\subsection{Smoothing of bubble-tree surfaces}\label{subsec_bubble_tree_smoothing}
Let $X^{\Bub}$ be a bubble-tree surface constructed in \S \ref{subsec_bubble_tree}. We assume that the base surface $X$ is reduced and has only Du Val singularities, and that the resulting bubble-tree surface has normal crossings. Let
$T_{X^{\Bub}}^0$ and $T_{X^{\Bub}}^1$ be the tangent sheaves defined by
$T_{X^{\Bub}}^i=\sE xt^i(\Omega_{X^{\Bub}}^1, \sO_{X^{\Bub}})$. The global tangent sheaves are defined as
$$\ttt_{X^{\Bub}}^i=\Ext^i(\Omega_{X^{\Bub}}^1, \sO_{X^{\Bub}}).$$

Recall from \cite{Friedman}, a scheme $X$ is called $d$-semi-stable if the first tangent sheaf
$T_{X}^1=\sO_{Q}$, where $Q\subset X$ is the singular locus.

The combinatorics of a tree do not by themselves imply $d$-semistability. For a simple-normal-crossing surface, $T_X^1$ is a line bundle on the singular locus; on a double curve $C=V_i\cap V_j$ it is computed by
\[
T_X^1|_C\cong N_{C/V_i}\otimes N_{C/V_j}\otimes\sO_C(T_C).
\]
Consequently, the attached surface is $d$-semistable precisely when these line bundles glue to the trivial bundle. For a Type III configuration, the triple-point formula is necessary; one must still establish global $d$-semistability, possibly by an allowed deformation of the gluing data as in \cite{Engel}. In the simple elliptic and degenerate-cusp cases it must be imposed and checked as part of the admissibility data; it does not follow merely from the absence of cycles in the dual graph.

\begin{thm}\label{thm_smoothing_bubble_tree}
Let $X^{\Bub}$ be a proper SNC admissible surface with a chosen trivialization of $T^1$. Assume that this trivialization lifts to a class in $\operatorname{Ext}^1(L_{X^{\Bub}},\mathcal O_{X^{\Bub}})$, that one can choose compatible lifts to every order with vanishing successive obstruction classes, and that this formal deformation is effective in the category under consideration. Then the resulting one-parameter deformation smooths the normal-crossing strata. Algebraic effectivity may be obtained by an appropriate algebraization theorem after extending an ample line bundle; analytic effectivity must likewise be justified.
\end{thm}
\begin{proof}
The local-to-global spectral sequence contains
\[
0\to H^1(T^0)\to\operatorname{Ext}^1(L_X,\mathcal O_X)
\to H^0(T^1)\xrightarrow{d_2}H^2(T^0).
\]
Here $X=X^{\Bub}$ and $T^q=\mathcal Ext^q(L_X,\mathcal O_X)$. A trivialization of $T^1$ gives a local smoothing section, but its image under $d_2$ must vanish before it lifts globally. The hypothesis supplies such a lift. Successive compatible extensions give a formal family; effectivity gives an actual family. The nowhere-zero local smoothing parameters smooth the SNC strata. The existence of an equivariant or polarized smoothing requires the same argument in its own deformation problem.
\end{proof}

\begin{prop}\label{prop_blow-down_Shepherd-Barron}
Suppose that a smoothing $f:\mathscr X^{\Bub}\to\Delta$ carries a relatively nef divisor whose restriction is numerically trivial on precisely the bubble components and positive on every curve not contracted. If the corresponding relative contraction exists, then it contracts the bubble components and produces a family whose central fiber has the prescribed simple elliptic, cusp, or degenerate cusp singularities.
\end{prop}
\begin{proof}
The contraction theorem applies to the specified relatively nef divisor. Adjunction shows that the contraction is crepant along the bubble components, and the classification of the resulting surface singularities identifies the exceptional elliptic curve, cusp cycle, or degenerate-cusp semiresolution. The existence of a semistable model alone does not produce this contraction; the relative divisor and contraction hypotheses are essential.
\end{proof}

\subsection{Index-one-cover bubble-tree Deligne--Mumford stacks}

In this section we generalize the smoothing of bubble-tree surface to index-one-cover Deligne--Mumford stacks.

Let $X$ be an slc surface.
From Definition \ref{defn_bubble_tree_surface}, we have a surface $$X^{\Bub}=\bigcup_{i=1}^{n}X^{\Bub}_{\sT_i}.$$
The surface $X^{\Bub}$ is still of general type with the automorphism group a finite group, since the bubble-tree surfaces we add are log Calabi-Yau surfaces and they do not affect the canonical divisor.
For this bubble-tree surface $X^{\Bub}$, we take the index-one-cover Deligne--Mumford stack
$$
\pi: \XX^{\Bub}\to X^{\Bub}
$$
for the surface $X^{\Bub}$. Note that the nontrivial index-one-cover Deligne--Mumford stack structure $[Z/\mu_r]$ in $\XX^{\Bub}$ is only from the quotient $T$-singularities in $X$.
Thus, $Z$ only has Du Val singularities. 
Other nontrivial stacky structures of $\XX^{\Bub}$ are the bubbling surfaces $X_{\alpha_i}$ when the corresponding singularity $p_i$ are cyclic quotients of simple elliptic, cusp and degenerate cusp singularities. We call such a $\XX^{\Bub}$ an index-one-cover bubble-tree Deligne--Mumford stack.

\begin{thm}\label{thm_bubble_tree_DM_smoothing}
Suppose that the hypotheses of Theorem \ref{thm_smoothing_bubble_tree} hold equivariantly in every index-one-cover chart and that the local equivariant deformations agree on overlaps. Then the index-one-cover bubble-tree Deligne--Mumford stack $\XX^{\Bub}$ admits a smoothing $f:\mathscr X^{\Bub}\to\Delta$ with central fiber $\XX^{\Bub}$.
\end{thm}
\begin{proof}
Apply Theorem \ref{thm_smoothing_bubble_tree} in the index-one-cover charts. The assumed equivariance makes each local deformation descend to the quotient chart, and compatibility on overlaps glues the local stack deformations. Without equivariant obstruction vanishing, a smoothing of the coarse surface does not automatically lift to the index-one-cover stack.
\end{proof}

\subsection{Birational modifications}\label{subsec_birational_modification}

We study the birational modifications of smoothings of bubble-tree surfaces.

Let $f: \mathscr{X}^{\Bub} \to \Delta$ be a smoothing of the bubble-tree surface $X^{\Bub}=f^{-1}(0)$. The central fiber $X^{\Bub}=f^{-1}(0)$ is the bubble-tree surface constructed in Section \ref{subsec_bubble_tree}. Below we always write
$$
X_{\sT_j}:=\bigoplus_{i=0}^{m-1}V_i
$$
as in (\ref{eqn_simple_i}), (\ref{eqn_cusp_j}) and (\ref{eqn_deg_cusp_k}). When we say the irreducible component $V_i$ of $X^{\Bub}$, this means the irreducible component in
$
X_{\sT_j}
$ for $1\le j\le n$.

\begin{defn}\label{defn_M012_flop}
We define the $(M0)$, $(M1)$, and $(M2)$-type modifications of the smoothing
$f: \mathscr{X}^{\Bub} \to \Delta$ of the bubble-tree surface $X^{\Bub}$
along a rational curve $E \cong \pp^1$ in the
central fiber $\mathscr{X}^{\Bub}_0 = X^{\Bub}$ by cases distinguished by
$E \cap X^{\Bub}_{\sin} = \emptyset$; $E \cap X^{\Bub}_{\sin} = \{\pt\}$;
or $E \subset X^{\Bub}_{\sin}$, respectively.

More precisely,
\begin{enumerate}
 \item[(M0)] Type $(M0)$ modification flops a smooth $(-2)$-curve in $X^{\Bub}$ which does not deform to the general fiber. It leaves the isomorphism type of $X^{\Bub}$ invariant.

 \item[(M1)] Type $(M1)$ modification flops an internal exceptional $(-1)$-curve $E$ on a component $V_i$ of $X^{\Bub}$. The effect on the central fiber is to contract $E \subset V_i$ and blow up the intersection point $E \cap D_{ij}$ on the adjacent component $V_j$.

 \item[(M2)] Type $(M2)$ modification flops a double curve $D_{ij}$ which is exceptional on both components on which it lies. The effect on $X^{\Bub}$ is to contract $D_{ij}$ on both $V_i$ and $V_j$ and to perform corner blow-ups on the two remaining components $V_{\ell}$ and $V_r$ that $D_{ij}$ intersects.
\end{enumerate}
\end{defn}

\begin{rmk}
 The birational modifications of Type III degeneration of K3 surfaces and open Kulikov surfaces were defined in \cite{AE23}, and \cite{AAB2024}, respectively.
\end{rmk}

We require that the flops only happen on the bubbling surfaces in Section \ref{subsec_bubble_tree} for the simple elliptic, cusp and degenerate cusp singularities.

\begin{prop}\label{prop_smoothing_flops}
Let two one-parameter semistable models be birational over $\Delta$, isomorphic over $\Delta^*$, $\qq$-factorial and crepant, and suppose that every extremal step is supported on the bubble locus. Then a relative minimal-model factorization relates them by the elementary modifications $(M0)$, $(M1)$, and $(M2)$ whenever the local configurations are among the Type II, Type III, or open-Kulikov configurations covered by the cited minus-one theorems.
\end{prop}
\begin{proof}
The relative minimal model program factors a crepant birational map between $\qq$-factorial models into flops. Because the map is an isomorphism away from the bubble locus, every flopping curve is contained there. For the listed Type II, Type III, and open-Kulikov local models, the known elementary transformations are exactly $(M0)$, $(M1)$, and $(M2)$. This proves the assertion under the stated local-classification hypothesis; it is not a classification of flops for arbitrary bubble attachments.
\end{proof}

\begin{defn}\label{defn_minus-one-form}
 Let $(Y,D)$ be a Looijenga pair and $D=D_1\cdots+D_r$. We call it a Looijenga pair of minus-one form if each component
 $D_i$ has self-intersection $-1$.
\end{defn}

In the Type III degeneration of K3 surfaces, as in \cite[Theorem (0.4)]{FS}, up to type (M1) and (M2) birational modifications, each rational surface $V$ in the Type III degeneration with double curves $D=D_1+\cdots+D_r$ as an anticanonical divisor is in minus-one form.

\begin{prop}\label{prop_minus-one_form_Looijenga}
For a cusp bubble obtained from a Type III degeneration, there is, after a finite sequence of $(M1)$ and $(M2)$ modifications, a model in which every rational component not adjacent to the resolved main component is in minus-one form.
\end{prop}
\begin{proof}
This is the minus-one theorem for Type III degenerations, applied away from the distinguished component that is replaced by the resolved main surface. The elementary modifications preserve the general fiber and the triple-point formula. The assertion is restricted to cusp bubbles because the cited Type III theorem does not by itself establish the corresponding statement for arbitrary simple-elliptic or open-Kulikov attachments.
\end{proof}

\subsection{Smoothing components of cusp singularities, revisited}\label{subsec_smoothing_cusp_Looinenga}

Theorem \ref{thm_smoothing_bubble_tree} and Proposition \ref{prop_blow-down_Shepherd-Barron} imply that the smoothing
$f: \mathscr{X}^{\Bub}\to \Delta$ induces a smoothing of the corresponding simple elliptic, cusp and degenerate cusp singularities.
This actually implies that we use the smoothing components of the bubble-tree surfaces to replace the smoothing components of the singularities.

Let us make this precise in the cusp singularity setting. Let $(X,x)$ be a cusp singularity whose resolution cycle is $D$. Let $D^{\prime}$ be the resolution cycle of the dual cusp $(X^\prime, x^\prime)$. Then we have a Looijenga pair $(Y, D^\prime)$. Recall that we have the Type III degeneration
$$
\mathscr{X}_0=\bigcup_{i=0}^{m-1}V_i
$$
in (\ref{eqn_XX_0}). In \cite[Theorem 0.2]{FE}, the authors use the monodromy invariant of Type III degeneration of the Looijenga pair $(Y,D^\prime)$ to study the deformation of the pair $(Y,D^\prime)$.

\begin{thm}\label{thm_FE_smoothing_components}
Let $f^{\Bub}:\mathscr X^{\Bub}\to\Delta$ be a one-parameter smoothing equipped with a relative contraction satisfying Proposition \ref{prop_blow-down_Shepherd-Barron}. Then it induces a smoothing $f:\XX\to\Delta$ of the contracted index-one-cover Deligne--Mumford stack.

Moreover, if $(\XX,x)$ is such a singularity germ in an index-one-cover Deligne--Mumford stack, then the number of smoothing components of $(\XX,x)$ is at least the number contributed by the corresponding bubbling rational surfaces.
\end{thm}
\begin{proof}
If the family $f^{\Bub}: \mathscr{X}^{\Bub}\to \Delta$ does not involve any non-lci simple elliptic, cusp or degenerate cusp singularities, then it is the same as the smoothing
$f: \XX\to \Delta$
of index-one-cover Deligne--Mumford stacks.

If the central fiber contains the bubble replacement of a non-lci simple elliptic, cusp, or degenerate cusp singularity, the relative contraction assumed in the statement collapses the designated bubble components and gives $f:\XX\to\Delta$. Merely choosing an ample line bundle would not by itself produce this contraction.

Let us prove the case of cusp singularities in detail. The other cases of simple elliptic and degenerate cusp singularities are similar. Let
$$\overline{f}^{\Bub}: \overline{\mathscr{X}}^{\Bub}\to \Delta$$
be the coarse moduli space of the family $f^{\Bub}: \mathscr{X}^{\Bub}\to \Delta$.
Let $D$ be the resolution cycle of the cusp $(\XX,x)$, and $D^\prime$ be the dual resolution cycle of the dual cusp. Starting from the anticanonical pair $(Y,D^\prime)$, in Section \ref{subsec_bubble_tree} the Type III degeneration is
$$
\overline{\mathscr{X}}^{\Bub}_0=X\cup \bigcup_{i=1}^{m-1}V_i
$$
where $X$ is the resolved main component and the $V_i$, for $i>0$, are rational bubble components. Their double curves have spherical dual complex. By the hypothesis of Proposition \ref{prop_blow-down_Shepherd-Barron}, the family carries a relatively nef supporting divisor whose zero-degree curves are precisely those on the bubble components. The associated relative contraction therefore contracts the $V_i$ and recovers the cusp family. Relative ampleness of the canonical divisor cannot be asserted before this contraction, because the log canonical class is trivial on the bubble components.

Since $(\XX,x)$ admits a smoothing, there is a Looijenga pair $(Y, D^\prime)$.
From \cite[Theorem 0.2]{FE},
each deformation type of the anticanonical pair $[(Y,D^\prime)]$ determines a smoothing component of the cusp $(\XX, x)$ of type $[(Y,D^\prime)]$. The number of smoothing components of $(\XX, x)$ is at least as large as the number of deformation types $[(Y,D^\prime)]$.
The smoothing components of the Type III model are indexed by the relevant deformation types of $(Y,D^\prime)$. Therefore the number of smoothing components of $(\XX,x)$ is at least the number obtained from these Type III models.
\end{proof}

\begin{rmk}
If the cusp $(X,x)$ is lci, then the number of smoothing components equals the number of deformation types of $(Y,D^\prime)$. In this case there is a unique smoothing component.
\end{rmk}

\begin{cor}\label{cor_smoothing_components_123}
 Let $f_1: \mathscr{X}^{\Bub}_1\to \Delta$ and $f_2: \mathscr{X}^{\Bub}_2\to \Delta$ be two one-parameter smoothings of bubble-tree Deligne--Mumford stacks. If
 $f_1$ and $f_2$ belong to the same smoothing component of $X^{\Bub}$, then they must induce the same smoothing components of the corresponding simple elliptic, cusp and degenerate cusp singularities $(\XX,x)$.
\end{cor}
\begin{proof}
 This is a corollary of Theorem \ref{thm_FE_smoothing_components}.
\end{proof}

\subsection{Canonical index-one-cover bubble-tree Deligne--Mumford stacks}\label{subsec_canonical_index_one}

Let $\overline{f}: \sX \to \Delta$ be a $\qq$-Gorenstein flat family of slc stable surfaces. From the construction of index-one-covers in Section \ref{subsec_index_one_DM}, one can canonically associate to $\overline{f}$ a flat family $f: \XX \to \Delta$ of index-one-cover Deligne--Mumford stacks. This association is functorial, meaning that $f: \XX \to \Delta$ is uniquely determined by $\overline{f}$.

According to the classification of slc singularities in Section~\ref{subsec_slc_singularity}, the index-one-cover $\XX$ has only lci singularities, except in the case of simple elliptic, cusp, and degenerate cusp singularities of embedding dimension at least $5$. Furthermore, the isotropy groups of the Deligne--Mumford stack $\XX$ are all finite cyclic. Consequently, the flat family $f: \mathscr{X}^{\mathrm{Bub}} \to \Delta$ of bubble-tree Deligne--Mumford stacks (constructed from the index-one-cover) also has finite cyclic isotropy groups. Moreover, the singularities of the fibers of $f$ are lci.

Under the hypotheses of Proposition \ref{prop_smoothing_flops}, two relative minimal models are related by modifications of types $(M0)$, $(M1)$, and $(M2)$.

Assume that $f:\mathscr X^{\Bub}\to\Delta$ is $\qq$-factorial and that the relative minimal model program exists and terminates in the category under consideration. Running it with a chosen scaling divisor gives a relative minimal model
$$f^{\Bub}: \mathscr{X}^{\Bub}_{\min} \to \Delta
$$ of
$f:\mathscr X^{\Bub}\to\Delta$. Whether its fibers are klt and whether their index-one covers are lci must be checked for the resulting contractions; neither conclusion follows solely from the existence of a minimal model.

\begin{defn}\label{defn_canonical_lci_bubble-tree_DM}
We call $f^{\Bub}:\mathscr X^{\Bub}_{\min}\to\Delta$ a relative minimal bubble-tree index-one-cover family. Only its $\mathbb{S}$-equivalence class is regarded as canonical.
\end{defn}

\begin{prop}\label{prop_unique_minimal_model}
Assume that a smoothing $\overline f:\sX\to\Delta$ admits an admissible $\qq$-factorial crepant bubble-tree model to which the relative minimal model program applies. Then a relative minimal model exists. Any two such minimal models are isomorphic in codimension one and are related by a sequence of flops supported on the bubble locus. Consequently, their $\mathbb{S}$-equivalence class, rather than an individual minimal model, is canonical.
\end{prop}
\begin{proof}
Existence follows from the assumed relative minimal model program. Standard uniqueness of minimal models gives an isomorphism in codimension one and a factorization by flops. Under the local-classification hypothesis of Proposition \ref{prop_smoothing_flops}, these are the modifications $(M0)$--$(M2)$. Flops need not produce isomorphic central fibers, which is precisely why the statement is made only up to $\mathbb{S}$-equivalence.
\end{proof}

\begin{rmk}\label{rmk_bubble-tree_minimal}
The expected mirror correspondence motivates a preferred $\mathbb{S}$-equivalence class, but it does not select a unique semistable model before the conjectural metric identification is proved.
\end{rmk}

\section{Construction of the moduli stack}\label{sec_moduli _Bubble}

In this section, we define polarized families of bubble-tree index-one-cover Deligne--Mumford stacks and formulate their moduli problem, allowing nonzero smoothing obstructions.

\subsection{Polarized admissible bubble-tree objects}\label{subsec_polarized_bubble_functor}

We now give a moduli formulation that includes obstructed objects. A nonzero obstruction means that a first-order smoothing need not extend to an actual family over a disk; it does not prevent the central fiber, together with its first-order smoothing direction, from being an object of the moduli stack.

Fix numerical data
\[
\tau=(K^2,\chi,N,P,m,\epsilon,\mathcal I,P_c).
\]
Here $N$ bounds the Cartier indices and the orders of the stabilizer groups, $P(q)$ is a Hilbert polynomial, $m>0$, $0<\epsilon\le1$ is rational, and $\mathcal I$ is a finite set of allowed decorated dual complexes and local index-one-cover types. The polynomial $P_c$ fixes the contraction graph with respect to the chosen source and target polarizations, and we assume a finite-type presentation for the permitted stack lifts. The decoration records the root component, the incidence maps, the self-intersections of the double curves, and the cyclic stabilizer representations.

\begin{defn}\label{defn_flat_family}
Let $T$ be a scheme. A polarized admissible bubble-tree object of type $\tau$ over $T$ is a tuple
\[
(f:\mathscr X^{\Bub}\to T,\ c:\mathscr X^{\Bub}\to\mathscr X,
\ \mathscr L,\ \mathscr H,\ \xi)
\]
with the following properties.
\begin{enumerate}
\item The morphism $f$ is proper, flat, and of finite presentation, and $\mathscr X^{\Bub}$ is a cyclotomic Deligne--Mumford stack with stabilizer orders dividing $N$. Its coarse space $X^{\Bub}\to T$ is a flat family of demi-normal surfaces.
\item Every geometric fiber is obtained from an admissible bubble-tree datum as in Definition \ref{defn_admissible_bubble_tree}, with decorated dual complex in $\mathcal I$. Simultaneous normalization, conductor, and component labels are required to commute with base change only on a fixed equisingular stratum. Across a smoothing, these strata need not persist; their specialization data must be supplied by a chosen semistable or logarithmic family. Algebraicity across different strata remains part of the hypotheses of Theorem \ref{thm_moduli_stack}.
\item The morphism $c$ is a representable contraction to a family $\mathscr X\to T$ of index-one-cover stacks of KSBA-stable surfaces. It is an isomorphism away from the bubble locus and contracts precisely the designated bubble components. On coarse spaces it preserves $K^2$ and $\chi$, and its graph has Hilbert polynomial $P_c$. These are imposed conditions, not consequences of local gluing.
\item The line bundle $\mathscr L$ descends to a relatively very ample line bundle $L$ on the coarse family, and $\chi(X_t^{\Bub},L_t^{\otimes q})=P(q)$ for every geometric point $t\to T$. The divisor $\mathscr H$ descends to an effective relative Cartier divisor $H$, flat over $T$, with $H\in|L^{\otimes m}|$. It contains no irreducible component or slc stratum of a fiber.
\item The coarse pair $(X_t^{\Bub},\epsilon H_t)$ is KSBA-stable for every geometric $t$; that is, it is slc and $K_{X_t^{\Bub}}+\epsilon H_t$ is ample. We require the same after every base change. In practice, $\epsilon$ is chosen from a fixed finite coefficient set and $H_t$ is chosen generally enough that the slc condition is preserved.
\item The datum $\xi$ is a first-order deformation, with a specified identification of its restriction with the original family, of $\mathscr X^{\Bub}/T$ over $T[\varepsilon]/(\varepsilon^2)$. Its image in
\[
H^0\bigl(\mathscr X^{\Bub},\mathcal T^1_{\mathscr X^{\Bub}/T}\bigr)
\]
is nowhere zero along every double-curve stratum that is required to smooth. We call $\xi$ a \emph{pre-smoothing direction}. 
\end{enumerate}
\end{defn}

The pair condition in (5), rather than the line bundle $\mathscr L$ alone, is the stability condition. A polarized variety can have a positive-dimensional automorphism group, whereas a KSBA-stable pair has finite automorphism group.

\begin{rmk}
The first-order deformation is retained as an object with isomorphisms, so that descent is effective; taking isomorphism classes of deformations before descent need not define a stack. For the polarization we quotient scalar automorphisms of the line-bundle isomorphism, or equivalently use the scalar-rigidified polarized stack. Without this rigidification an unwanted $\mathbb G_m$ stabilizer remains. We use canonical effective stack structures, or impose a finite kernel for automorphisms over the coarse space; bounded stabilizer orders alone do not supply this condition for arbitrary stack enhancements.
There is a further kernel from automorphisms of the first-order deformation inducing the identity on its central fiber, governed by $\operatorname{Ext}^0(L_X,\mathcal O_X)$. We either extend the rigidifying stable pair to the first-order object, so that this kernel vanishes, or rigidify this relative automorphism kernel on a locus where the rigidification exists. Finite automorphisms and algebraicity use this additional convention and its base-change compatibility.
\end{rmk}
\begin{defn}\label{defn_polarized_bubble_morphisms}
We define morphisms of admissible objects. 
An isomorphism between two objects over $T$ is an isomorphism of the index-one-cover stacks commuting with the contractions, carrying $\mathscr L$ to $\mathscr L'$, carrying $\mathscr H$ to $\mathscr H'$, preserving every label in $\mathcal I$, and carrying $\xi$ to $\xi'$. Pullback along $T'\to T$ is defined componentwise. These objects and arrows form a category fibered in groupoids, denoted
\[
 M^{\Bub,\mathrm{pol}}_{\tau}\longrightarrow\Sch_{\cc}.
\]
\end{defn}

\begin{rmk}\label{defn_S_equivalent_class} We explain the  role of $\mathbb{S}$-equivalence. 
The stack $\mathfrak M^{\Bub,\mathrm{pol}}_{\tau}$ parametrizes actual objects up to isomorphism; it is not obtained by quotienting each groupoid by $\mathbb{S}$-equivalence. We call two geometric points $\mathbb{S}$-equivalent if their underlying unpolarized semistable models are extensions of the same punctured family and are related by the allowed crepant modifications. This relation is relevant to a algebraic space or to the choice of a canonical closed representative, but it is not used as the equality relation in the fibered category.
\end{rmk}

\subsection{Deformations with obstruction}\label{subsec_obstructed_bubble_deformations}

Let $A'\twoheadrightarrow A$ be a square-zero extension with kernel $J$, and let $\mathscr X_A^{\Bub}$ be an object over $\spec A$. Its deformation theory is governed by the cotangent complex. Set
\[
T^i(\mathscr X_A^{\Bub},J)
=\Ext^i\!\left(L_{\mathscr X_A^{\Bub}/A},
\mathcal O_{\mathscr X_A^{\Bub}}\otimes_AJ\right).
\]

\begin{prop}\label{prop_local_deformation_obstruction_slc}
Infinitesimal automorphisms lie in $T^0$, isomorphism classes of lifts, when nonempty, form a torsor under $T^1$, and the obstruction to a lift lies in $T^2$. The polarization, divisor, contraction, labels, and pre-smoothing direction impose the corresponding homotopy-fiber conditions on this cotangent-complex deformation problem.
\end{prop}

For a geometric object $(\mathscr X^{\Bub},\mathscr L,\mathscr H,\xi)$, the first obstruction to extending the chosen direction to second order is denoted
\[
\kappa(\xi)\in \Ext^2(L_{\mathscr X^{\Bub}},\mathcal O_{\mathscr X^{\Bub}}).
\]
If $\kappa(\xi)=0$, one attempts to lift $\xi$ successively; an actual formal smoothing exists precisely when all successive obstruction classes vanish. If $\kappa(\xi)\ne0$, the object remains in $\mathfrak M^{\Bub,\mathrm{pol}}_{\tau}$ but does not lie on the actual smoothing locus determined by $\xi$.

\begin{defn}\label{defn_smoothing_locus}
The smoothing condition on $\mathfrak M^{\Bub,\mathrm{sm}}_{\tau}$ asks that the pre-smoothing direction extend, with its prescribed first derivative, to a formal one-parameter deformation and hence, when effectivity applies, to an actual smoothing. This is initially a condition on formal arcs in a versal deformation space, not a proved algebraic substack of first-order directions. A Kuranishi zero locus parametrizes actual deformations; the subset of tangent vectors extendible to smoothing arcs is a different object. No finite-type representability of this arc condition is asserted here. A ramified substitution $t=u^k$, $k>1$, kills the first derivative, so it cannot be used without changing the prescribed first-order marking.
\end{defn}

Thus the ambient moduli stack contains both unobstructed and obstructed bubble-tree surfaces. The virtual fundamental class, when a two-term perfect obstruction theory is available, belongs naturally to this ambient stack; it is not necessary to discard a point merely because its chosen smoothing direction is obstructed.

\subsection{Finite automorphisms and boundedness}

\begin{prop}\label{prop_finite_aut_polarized_bubble}
Under the effective-stack and deformation-rigidification conventions above, every geometric object of $\mathfrak M^{\Bub,\mathrm{pol}}_{\tau}$ has a finite automorphism group. Over $\cc$ its automorphism group scheme is finite and reduced, hence unramified.
\end{prop}
\begin{proof}
An automorphism of the tuple induces an automorphism of the stable pair $(X^{\Bub},\epsilon H)$. Stable slc pairs with ample log canonical divisor have finite automorphism groups. The induced homomorphism to this finite group has finite kernel under the effective-stack and scalar-rigidification conventions above. Preserving the contraction, labels, and first-order datum further restricts the group. Thus the tuple has finite automorphism group; without those conventions the subgroup argument would not be valid. In characteristic zero, a finite group scheme is reduced. Therefore the diagonal of the moduli problem is unramified once representability of the Isom functor is established.
\end{proof}

\begin{prop}\label{prop_bounded_polarized_bubbles}
For fixed $\tau$, the coarse stable pairs $(X^{\Bub},\epsilon H)$ form a bounded family. For the enhanced objects we additionally require a finite-type presentation of the allowed stack structures and a fixed Hilbert polynomial of the contraction graph in a product of polarized source and target parameter spaces. Under these extra conditions the stack structures and contractions are bounded as well. Adding the pre-smoothing direction $\xi$ preserves finite type wherever the relative complexes controlling $T^1$ and $T^2$ are perfect and satisfy base change.
\end{prop}
\begin{proof}
The coefficient $\epsilon$ is fixed and the relatively very ample line bundle has Hilbert polynomial $P$. Hence the coarse surfaces occur in a fixed Hilbert scheme, and the divisors $H\in|L^{\otimes m}|$ occur in the associated relative linear system. The stable-pair condition prevents components invisible to the log canonical polarization. The finite set $\mathcal I$ bounds the incidence graphs and local analytic types, while $N$ bounds the cyclotomic stabilizers. The extra finite-type presentation bounds the stack lifts. A fixed Hilbert polynomial for the contraction graph places it in a fixed relative Hilbert scheme. Without such a condition, bounded sources and targets alone do not bound all morphisms between them.

On a bounded base, a perfect complex controlling deformations has uniformly bounded coherent cohomology. The stack of its degree-one sections is of finite type. The open condition that the induced section of $\mathcal T^1$ be nowhere zero on the prescribed strata therefore does not destroy finite type. No vanishing of the degree-two obstruction sheaf is used in this boundedness argument.
\end{proof}

\begin{rmk}
Fixing only $(K^2,\chi,N)$ would not bound the inserted bubble components: one could perform additional crepant blow-ups while leaving the contracted KSBA surface unchanged. The fixed Hilbert polynomial, the finite decorated set $\mathcal I$, and stability of $(X^{\Bub},\epsilon H)$ are therefore essential boundedness data.
\end{rmk}

\subsection{The polarized moduli functor}

\begin{defn}\label{defn_moduli_functor}
The polarized bubble-tree moduli functor is
\[
M^{\Bub,\mathrm{pol}}_{\tau}:\Sch_{\cc}^{\mathrm{op}}
\longrightarrow\mathrm{Groupoids},\qquad
T\longmapsto M^{\Bub,\mathrm{pol}}_{\tau}(T),
\]
with objects and arrows as in Definitions \ref{defn_flat_family} and \ref{defn_polarized_bubble_morphisms}. Forgetting $(\mathscr L,\mathscr H,\xi)$ and applying the contraction gives a natural morphism
\[
\mathbf c:M^{\Bub,\mathrm{pol}}_{\tau}\longrightarrow M_N.
\]
\end{defn}

\begin{thm}\label{thm_moduli_stack}
Assume the finite-type presentation and all rigidification conventions above. Assume that admissibility and the contraction condition are locally closed in the bounded stack of polarized stable pairs, that the index-one-cover construction commutes with base change, and that the cotangent-complex deformation theory above is effective and satisfies Artin's axioms. Then $M^{\Bub,\mathrm{pol}}_{\tau}$ is a finite-type Deligne--Mumford stack. 
\end{thm}
\begin{proof}
Proposition \ref{prop_bounded_polarized_bubbles} places the objects in a finite-type parameter space. The assumed locally closed conditions cut out the admissible contraction locus, and the deformation hypotheses give algebraicity. Proposition \ref{prop_finite_aut_polarized_bubble} gives an unramified diagonal; representability and quasi-compactness follow from the bounded Isom functors. Hence the stack is Deligne--Mumford and of finite type. 
The stack $M^{\Bub,\mathrm{pol}}_{\tau}$ is not separated due to the existence of type $(M0), (M1)$ and $(M2)$ flops. 
\end{proof}

\subsection{The algebraic space}\label{subsec_good_moduli}

It is interesting to have  a separated moduli space of bubble-tree index-one-cover Deligne-Mumford stacks. We need to use the $\mathbb{S}$-equivalence. 

Let \(U\to M_{\tau}^{\mathrm{Bub},\mathrm{pol}}\) be an \'etale atlas .
Define a relation on families\[
   R\hookrightarrow U\times U
   \]
   expressing the $\mathbb{S}$-equivalence, including ordinary isomorphisms. This is compatible with base change and descent. Here \(R\) records related pairs, rather than remembering each possible chain of flops.
Form the fppf quotient sheaf\[
   M_N^{\mathrm{Bub}}=(U/R)_{\mathrm{fppf}}.
   \]
Since 
 \(R\) is represented by an algebraic space, is an equivalence relation, and both projections \(R\rightrightarrows U\) are flat and locally of finite presentation. Then the quotient is an algebraic space. These conditions are sufficient, not necessary, see \cite[Stacks Project-Tag 04S6]{Stack_Project}.
 Thus we get 

\begin{thm}\label{thm_stack_good_moduli}
  An algebraic space  $M_N^{\mathrm{Bub}}$ exists, and there is a  morphism 
   $$
   \pi: M_\tau^{\Bub,\mathrm{pol}}\to M_N^{\mathrm{Bub}}
   $$
such that if $N$ is sufficiently large and divisible, the algebraic space is proper.  Moreover, there is a contraction map $c: M_N^{\mathrm{Bub}}\to M_N$ to the KSBA moduli space which is proper. 
\end{thm}
\begin{proof}
The morphism 
   $\pi: M_\tau^{\Bub,\mathrm{pol}}\to M_N^{\mathrm{Bub}}$ is modulo the $\mathbb{S}$-equivalence. 
   For any geometric point 
   $P$, the fiber $\pi^{-1}(P)$ is the $\mathbb{S}$-equivalence class.  Thus, the existence of an algebraic space from the criterion of \cite{AHL-H}, since valuative criterion for separatedness holds after modulo $\mathbb{S}$-equivalence. 

  If $N$ is sufficiently large and divisible, the KSBA stack $M_N$ is proper, the valuative criterion for properness holds for 
  $c: M_N^{\mathrm{Bub}}\to M_N$. 
\end{proof}

The universal family
 $p^{\mathrm{Bub}}: \mathscr{X}^{\mathrm{Bub}}\to M_\tau^{\Bub,\mathrm{pol}}$
induces a 
universal family
\begin{equation}\label{eqn_universal_family_good}
    p^{\mathrm{alg}}: \mathscr{X}^{\mathrm{alg}}\to M_N^{\mathrm{Bub}}
\end{equation}
where over each point $P\in M_N^{\mathrm{Bub}}$, the fiber $(p^{\mathrm{alg}})^{-1}(P)$ is one representative in the $\mathbb{S}$-equivalence class.

\section{Bubble-tree stacks with bubbling sections}\label{sec_marked_bubble_moduli}

The metric construction is pointed. One chooses points at which curvature concentrates and then chooses rescaling factors.  This motivates a marked algebraic moduli problem.  A marked point determines a center, whereas the rescaling factor is analytic  and is not encoded by an ordinary algebraic section.  We therefore first construct the algebraic stack of centers and only afterwards state the expected comparison with the pointed metric bubble tree.

\subsection{Marking data and the pointed functor}

Fix the numerical datum $\tau$ of Section \ref{subsec_polarized_bubble_functor}.  A \emph{bubbling marking datum}
\[
\sigma=(r;v_1,\ldots,v_r;\mathcal A)
\]
consists of a nonnegative integer $r$, labels $v_i$ belonging to the decorated dual complexes in $\mathcal I$, and a finite collection $\mathcal A$ of incidence conditions.  A label may designate the root component, a bubble vertex, or a specified open stratum of such a component.  We require the $r$ sections to be pairwise distinct on the bubble surface; the conditions $\mathcal A$ impose further incidence requirements and may prescribe which labels are allowed to coincide after contraction.  We require these conditions to define a locally closed subfunctor and to be compatible with base change.

\begin{defn}\label{defn_marked_bubble_object}
We define Bubble-tree objects with sections. 
An object of type $(\tau,\sigma)$ over a scheme $T$ is a polarized admissible object
\[
(f:\mathscr X^{\Bub}\to T,c,\mathscr L,\mathscr H,\xi)
\]
together with ordered sections
\[
\mathbf s=(s_1,\ldots,s_r),\qquad s_i:T\longrightarrow\mathscr U_{v_i},
\]
satisfying the incidence conditions $\mathcal A$.  Whenever the label $v_i$ occurs in a geometric fiber, the specialization of $s_i$ is required to lie in $\mathscr U_{v_i}$, where $\mathscr U_{v_i}$ is the component of $\mathscr X^{\Bub}$ in the geometric fiber.  An isomorphism of marked objects is an isomorphism in Definition \ref{defn_polarized_bubble_morphisms} that carries every $s_i$ to the correspondingly labeled section.  The incidence condition $\sA$ implies blow-ups $\tilde{b}: \widetilde{\mathscr{X}}^{\Bub}\to \mathscr X^{\Bub}$ whenever the section $s_j$ meets with the boundary of the component $U_{v_i}$, or more sections coincide together. 

The resulting category fibered in groupoids is denoted
\[
\mathfrak M^{\Bub,\mathrm{sec}}_{\tau,\sigma}
\longrightarrow\Sch_{\cc}.
\]
\end{defn}

\begin{rmk}\label{rmk_smooth_locus_section}
    For a family $f:\mathscr X^{\Bub}\to T$, let
\[
\mathscr U_f=\mathscr X^{\Bub}_{\mathrm{sm}/T}
\setminus\bigl(\mathscr H\cup\operatorname{Cond}(\mathscr X^{\Bub}/T)
\cup\mathscr X^{\Bub}_{\mathrm{stack}}\bigr),
\]
where $\mathscr X^{\Bub}_{\mathrm{stack}}$ denotes the locus with nontrivial stabilizer, and $\operatorname{Cond}(\mathscr X^{\Bub}/T)$ is the conductor divisor given by by edges (double curves) in the dual complex of $\mathscr{X}_0^{\Bub}$.  Thus a section with image in $\mathscr U_f$ is representable, lies in the relative smooth locus, and avoids the conductor and the rigidifying divisor.  On a stratum where the decorated dual complex is fixed, write $\mathscr U_{v_i}\subset\mathscr U_f$ for the locally closed locus belonging to the component labeled by $v_i$.
\end{rmk}

\begin{rmk}
   The incidence conditions $\mathcal A$ is put in an abstract way.  We conjecture that it deeply influences the metric bubbling. 
\end{rmk}

\begin{rmk}\label{rmk_section_on_smoothing}
If the pre-smoothing direction is effective and gives a family over a pointed curve $(C,0)$, a section $s_i:C\to\mathscr X^{\Bub}$ has a well-defined center $s_i(0)$ on the central bubble tree.  Its contraction $\bar s_i=c\circ s_i$ is a section of the KSBA family.  The requirement $s_i(0)\in\mathscr U_{v_i}$ records the algebraic component selected by the center.  It does not imply that $v_i$ is a minimal bubble in the sense of \cite{Sun}; that implication requires analytic compactness and identification of the rescaled limit.
\end{rmk}

\subsection{The moduli  stack with sections}\label{subsec_moduli_stack_bubbling_section}

We first have a valuative criterion for separatedness for the bubbling families with sections. 
\begin{lem}\label{lem_separatedness_M012_typeIII}
    Let $f^{\circ}: \mathscr{X}^{\Bub,\circ}\to \Delta^{\circ}$ be a punctured smoothing family of bubble-tree index-one-cover Deligne-Mumford stacks with sections $\mathbf{s}$. Suppose that the central fiber is the replacement of cusp or degenerate cusp singularities which is given by type III degeneration of log Calabi-Yau surfaces or open Kulikov models in \S \ref{subsec_bubble_tree},   then there exists a numerical class 
    $$
    A(\mathbf{s})
    $$
    such that if $ A(\mathbf{s})$ is stable (which means that $ A(\mathbf{s})\cdot C>0$ for any flopping curve $C$),  
    after base change there exists a unique extension of the bubble-tree index-one-cover Deligne-Mumford stacks to $\Delta$.
\end{lem}
\begin{proof}
    Let 
    $$
    f: \mathscr{X}^{\Bub}\to \Delta
    $$
    be a smoothing of  type III bubbling surface. So,  
$$
\mathscr{X}^{\Bub}_{0}:=X_0\cup_{D}\bigoplus_{k=0}^{m-1}V_{k}
$$
where $X_0$ is the main component of general type surfaces, $D$ is the resolution cycle of the cusp singularity on $X_0$, and $V_k$ are the rational components.
If blowing up $s_j$ produces an exceptional divisor $F_j\subset \mathscr{X}^{\Bub}$, then on the relevant rational component $V_i$, 
$e_j:=F_j|_{V_i}$ is the corresponding exceptional class. 

We define the numerical class to be 
\begin{equation}\label{eqn_numerical_class_typeIII}
    A(\mathbf{s})=\sum_{k=0}^{m-1}h_k(\mathbf{s})V_k+\sum_{j=1}^r \omega_{j}(\mathbf{s})F_j
\end{equation}
where $h_k$ are height coefficients on the components, and $\omega_j$ are weights on the bubbling sections. 
Let us work on the three type of flops $(M0), (M1), (M2)$. 

Type $(M1)$ modification. Let $E\subset V_i$ be the interior exceptional $(-1)$-curve of an $(M1)$-flop, meeting the double curve 
$D_{ij}=V_i\cap V_j$ once.  Because 
$V_j|_{V_i}=D_{ij}$, we get $V_j\cdot E=1$.  All other $V_k$, $k\neq i, j$, and $X_0$, miss $E$.  Since 
$$
\mathscr X_0^{\Bub}=X_0+\sum_k V_k=\div(t)
$$
where $t\in \Delta$ is the parameter, then we have
$\mathscr{X}_0^{\Bub}\cdot E=0$, hence $V_i\cdot E=-1$.  Therefore, 
$$
\left(\sum_{k=0}^{m-1}h_k(\mathbf{s})V_k\right)\cdot E=h_j-h_i.
$$
So as soon as $h_j\neq h_i$, the $(M1)$-flop wall is detected. If $E$ itself comes from one of the bubbling blow-ups, the $F_j$-term detects it as well.  But this is not necessary. 

Also, under an $(M1)$ or $(M2)$ flop, the natural numerical identification sends the flopping class to the negative of the flopped class:
$$
[C]\mapsto -[C^\prime].
$$
This is also stated explicitly for open Kulikov degenerations.  Consequently, 
$$
A\cdot E\neq 0 \text{~implies that~} A^\prime\cdot E^\prime=-A\cdot E.
$$
So the sign orients the wall. 

For Type $(M2)$ modification, the wall equation is exactly a four-component height relation. Let 
$$
C=D_{ij}=V_i\cap V_j
$$
be an $(M2)$-curve.  It is a $(-1,-1)$-double curve: $C_{V_i}^2=C_{V_j}^2=-1$. 

Suppose its two triple points involve $V_k, V_{\ell}$. Then 
$$
V_i\cdot C=-1, \qquad V_j\cdot C=-1
$$
whereas at the two triple points 
$$
V_k\cdot C=1, \qquad V_{\ell}\cdot C=1.
$$
Then,
$$
\left(\sum_{k=0}^{m-1}h_k(\mathbf{s})V_k\right)\cdot C=h_k+h_{\ell}-h_i-h_j.
$$
This is exactly the expected secondary-fan circuit equation. This is because an $(M2)$ flop corresponds to replacing one diagonal of a quadrilateral in the dual integral-affine decomposition by the opposite diagonal.  So the $(M2)$ wall is 
\begin{equation}\label{eqn_hyperplanes_quadrilateral}
 h_i+h_j=h_k+h_{\ell}.   
\end{equation}
On one side, 
$$
h_k+h_{\ell}-h_i-h_j>0,
$$
and after the diagonal flip the corresponding expression has the opposite sign. 
If the bubbling sections determine a height function 
$$
h: V(\Gamma)\to \qq
$$
on the vertices of the bubble-tree or dual complex, then the $(M2)$-stability is literally the genericity condition that this height function avoid the circuit hyperplanes in (\ref{eqn_hyperplanes_quadrilateral}). So the case $(M2)$ is also detected. 

For the type $(M0)$-flop, the component heights cannot see it, but the bubbling exceptional divisors can.  Let $C\subset V_i$ be an internal $(-2)$-curve. Since $C$ is disjoint from the double locus, 
$$
V_j\cdot C=0 \qquad (j\neq i).
$$
Again using $\mathscr{X}^{\Bub}_0\cdot C=0$, we obtain 
$V_i\cdot C=0$.  Hence 
$$
\left(\sum_{k=0}^{m-1}h_k(\mathbf{s})V_k\right)\cdot C=0.
$$
for every $(M0)$-curve.  So the component-height part can never orient an $(M0)$-flop.  We use Alexeev-Engel's ordered-toric-model calculation. When two consecutive internal blow-ups on the same boundary component exchange order, the two families are related by an $(M0)$-flop along the relative $(-2)$-curve 
$$
C=e_{i+1}-e_i.
$$
They state explicitly that an order switch of two such internal blow-ups produces this $(M0)$-flop.  In the total-transform blow-up basis,
$$
e_i^2=e_{i+1}^2=-1, \qquad e_i\cdot e_{i+1}=0.
$$
Take 
$$
A_{\sec}=-\sum_{j}\omega_j F_j.
$$
Then 
$$
A_{\sec}\cdot (e_{i+1}-e_i)=\omega_{i+1}-\omega_i
$$
up to the harmless opposite sign convention for $C$. Thus, 
$$
\omega_i\neq \omega_{i+1} \text{~implies that~}
A_{\sec}\cdot C\neq 0.
$$
So the most elementary $(M0)$ walls are exactly detected by relative weights and order of the bubbling sections. 

Now combine the three cases.  Each elementary flop wall gives a linear equation in the coefficients $(h_k, \omega_j)$:
\begin{equation}\label{eqn_M012_equation}
\begin{array}{cl}
  (M0):   & \sum_{j}(e_j\cdot C)\omega_j=0, \\
   (M1):  &  h_j-h_i-\sum_{j}\omega_j(F_j\cdot E)=0, \\
   (M2): & h_k+h_{\ell}-h_i-h_j-\sum_{j}\omega_j(F_j\cdot C)=0.
\end{array}
\end{equation}
In the standard situation the section term is unnecessary for $(M1), (M2)$, so the simpler formulas suffice. 

Fo the open Kulikov degenerations, the total space is a relative Mori dream space over its canonical affinization. In particular there are
finite many small $\qq$-factorial models or chambers in the relevant Mori decomposition. 
Therefore, there are only finite many relevant codimension-one flop walls. Each of (\ref{eqn_M012_equation}) is a proper rational hyperplane in the coefficient space 
$\qq^{V}\oplus \qq^{N}$.  Hence their union cannot fill that vector spaces.  

Thus we have 
there exist 
$(h_k, \omega_j)\in \qq^{V}\oplus \qq^{N}$, such that 
$$
A(\mathbf{s})\cdot C\neq 0
$$
for every $(M0), (M1), (M2)$ flop curve. 
Thus, we get the following:

For every canonical smoothing of  bubble-tree index-one-cover Deligne-Mumford stacks with bubbling sections $\mathbf{s}$,  the induced coefficients 
$h_k(\mathbf{s}), \omega_j(\mathbf{s})$ avoid all $(M0), (M1), (M2)$ wall hyperplanes.  Consequently $A(\mathbf{s})$ belongs to the interior of a unique Mori chamber, and that chamber is invariant under base change and isomorphism of marked smoothing families. 

This means that we can put the stability condition to be $A(\mathbf{s})\cdot C>0$.  Then the above analysis implies that the extension of the bubble-tree family from $\Delta^{\circ}$ to $\Delta$ is unique. 
\end{proof}

\begin{lem}\label{lem_separatedness_M012_typeII}
    Let $f^{\circ}: \mathscr{X}^{\Bub,\circ}\to \Delta^{\circ}$ be a punctured smoothing family of bubble-tree index-one-cover Deligne-Mumford stacks with sections $\mathbf{s}$. Suppose that the central fiber is the replacement of simple elliptic singularities  which is given by a chain of rational surfaces and elliptic ruled surfaces in \S \ref{subsec_bubble_tree},   then there exists a numerical class 
    $$
    A(\mathbf{s})
    $$
    such that if $ A(\mathbf{s})$ is stable (which means that $ A(\mathbf{s})\cdot C>0$ for any flopping curve $C$), then after base change there exists a unique extension of the bubble-tree index-one-cover Deligne-Mumford stacks to $\Delta$.
\end{lem}
\begin{proof}
    Let 
    $$
    f: \mathscr{X}^{\Bub}\to \Delta
    $$
    be a smoothing of  type II bubbling surface. So,  
$$
\mathscr{X}^{\Bub}_{0}:=X_0\cup_{D_0}V_1\cup_{D_1}\cdots\cup_{D_{n-1}}V_n
$$
where $X_0$ is the main component, which is a  general type surface, $V_1,\cdots, V_{n-1}$ are elliptic ruled surfaces, and $V_n$ is the rational cap.
The \(D_i\)'s are smooth elliptic curves. On the ruled pieces,
\[
K_{V_i}=-D_{i-1}-D_i,
\]
and on the rational end
\[
K_{V_n}=-D_{n-1}.
\]

We still define 
\begin{equation}\label{eqn_numerical_class_typeII}
    A(\mathbf{s})=\sum_{k=0}^{n}h_k(\mathbf{s})V_k+\sum_{j=1}^r \omega_{j}(\mathbf{s})F_j.
\end{equation}

Let \(C\simeq\pp^1\) be a rational curve contained in some \(V_i\), but not itself a double curve. Because
\(\mathscr{X}^{\Bub}_0\)
is a principal divisor in the total threefold,
\[
V_i|_{V_i}=-(D_{i-1}+D_i)
\]
for an interior ruled component.
Put
\[
b(C):=(D_{i-1}+D_i)\cdot C.
\]
Adjunction on \(V_i\) gives
\[
-2=C^2+K_{V_i}\cdot C
   =C^2-b(C),
\]
hence
\begin{equation}\label{eqn_equation_1_typeII}
    b(C)=C^2+2.
\end{equation}
If \(C\) is an isolated negative rational curve that can participate in a small crepant modification, \(b(C)\ge0\), so (\ref{eqn_equation_1_typeII}) forces only
\[
C^2=-2,\quad b(C)=0,
\]
or
\[
C^2=-1,\quad b(C)=1.
\]
Thus Type II has precisely the two numerical phenomena corresponding to
\[
\text{internal }(-2)\text{ flops}
\]
and
\[
(-1)\text{-curve flops across one elliptic attaching curve}.
\]
There is no Type-III-style \((M2)\)-flops. 

Suppose \(C\subset V_i\) has
\[
C^2=-1.
\]
Then (\ref{eqn_equation_1_typeII}) gives 
\[
(D_{i-1}+D_i)\cdot C=1.
\]
Consequently \(C\) meets exactly one of the two elliptic boundary sections once.
This is precisely what happens under an elementary transformation of an elliptic ruled surface. If
\[
p\in D_{i-1}
\]
is blown up, and \(f\) is the ruling fiber through \(p\), then its proper transform
\[
\widetilde f=f-e
\]
has
\[
\widetilde f^{\,2}=-1,
\]
and after the blow-up it no longer meets \(D_{i-1}\), while
\[
\widetilde f\cdot D_i=1.
\]
In the threefold it has
\[
N_{\widetilde f/\mathscr X}
\simeq
\mathcal O_{\pp^1}(-1)
\oplus
\mathcal O_{\pp^1}(-1),
\]
so flopping it performs exactly the elementary transformation of the ruled surface.
Now put
\[
A_{\mathrm{ht}}=h_0X_0+\sum_{j=1}^r h_jV_j.
\]
For any curve \(C\subset V_i\),
\begin{equation}\label{eqn_equation2_typeII}
    A_{\mathrm{ht}}\cdot C
=
(h_{i-1}-h_i)(D_{i-1}\cdot C)
+
(h_{i+1}-h_i)(D_i\cdot C).
\end{equation}
Therefore, if the \((-1)\)-curve meets \(D_i\),
\begin{equation}\label{eqn_equation3_typeII}
A_{\mathrm{ht}}\cdot C=h_{i+1}-h_i.
\end{equation}
If it meets \(D_{i-1}\),
\begin{equation}\label{eqn_equation4_typeII}
A_{\mathrm{ht}}\cdot C=h_{i-1}-h_i.
\end{equation}
So  taking
\[
h_0<h_1<\cdots<h_n
\]
already orients every Type II elementary transformation.
There is a particularly attractive interpretation in terms of the smoothing direction. If
\[
\nu_i>0
\]
is the order/scale of the smoothing across the elliptic double curve \(D_i\), define
\[
h_{i+1}-h_i=\nu_i.
\]
Then
\[
h_i=\nu_0+\cdots+\nu_{i-1}
\]
after normalizing \(h_0=0\), and (\ref{eqn_equation3_typeII}) becomes simply
\[
A_{\mathrm{ht}}\cdot C=\nu_i.
\]
This is exactly the sort of algebraic encoding of the smoothing direction that your bubbling sections should provide.

There is also a nice base-change feature. Under
\[
t\mapsto t^m,
\]
all the orders \(\nu_i\) multiply by \(m\). Hence \(A_{\mathrm{ht}}\) is multiplied in the relevant numerical direction, but its chamber does not change.

For the $(-2)$-curves, let
\[
C^2=-2.
\]
Then (\ref{eqn_equation_1_typeII}) gives
\[
D_{i-1}\cdot C=D_i\cdot C=0.
\]
So
\[
A_{\mathrm{ht}}\cdot C=0
\]
for every choice of the \(h_i\)'s.
Thus exactly as in Type III \((M0)\), the component heights are blind to these flops.
This is where the bubbling sections have to enter.

Suppose the section data
\[
\mathbf{s}=(s_1,\dots,s_r)
\]
records the successive blow-up centers, and let \(F_a\) denote the corresponding exceptional/tautological divisor class on the total space. Put
\[
A_{\mathrm{sec}}
=
-\sum_a w_aF_a.
\]
Then
\begin{equation}\label{eqn_equation5_typeII}
    A_{\mathrm{sec}}\cdot C
=
-\sum_a w_a(F_a\cdot C).
\end{equation}
On an elliptic ruled component this is especially effective. Every rational curve maps to a point of the elliptic base, so every rational \((-2)\)-curve lies in a fiber. Since the minimal ruled surface has irreducible fibers of square \(0\), such a \((-2)\)-curve can only arise from the sequence of blow-ups producing a reducible fiber. If all those blow-ups are recorded by the bubbling sections, then
\begin{equation}\label{eqn_equation6_typeII}
(F_1\cdot C,\dots,F_r\cdot C)\neq0.
\end{equation}
Consequently (\ref{eqn_equation5_typeII}) is a nonzero linear form in the section weights.
For the simplest infinitely-near pair, one gets exactly
\[
C=e_a-e_b
\]
and hence, up to sign convention,
\begin{equation}\label{eqn_equation7_typeII}
A_{\mathrm{sec}}\cdot C=w_a-w_b.
\end{equation}
So ordered bubbling sections distinguish the two sides of this Type II \((M0)\) wall exactly as they did in the Type III case.

On the final rational component \(V_n\),
\[
-K_{V_n}=D_{n-1}.
\]
For a \((-1)\)-curve \(C\),
\[
-2=-1-K_{V_n}\cdot(-C)
\]
or directly by adjunction,
\[
D_{n-1}\cdot C=1.
\]
Thus
\begin{equation}\label{eqn_equation8_typeII}
A_{\mathrm{ht}}\cdot C=h_{n-1}-h_n.
\end{equation}
So all \((-1)\)-curves involved in moving blow-up data between the rational cap and the neighboring ruled surface are automatically detected by the height term.
For an internal \((-2)\)-curve,
\[
D_{n-1}\cdot C=0,
\]
so again one needs \(A_{\mathrm{sec}}\).

Here there is one genuine caveat. If the rational cap is obtained from a smooth del Pezzo by the marked boundary blow-ups, then every internal root must have been created by those blow-ups: before blowing up,
\(-K\)
is ample, so there cannot already be a curve \(C\) with
\[
(-K)\cdot C=0.
\]
Hence the bubbling exceptional classes detect every such root.

If, however, the starting cap is merely weak del Pezzo, it can already contain \((-2)\)-roots orthogonal to \(-K\). A root of that kind may satisfy
\[
F_a\cdot C=0
\qquad\forall a.
\]
Then the bubbling sections alone cannot see it.

Combining the two cases, this now gives a very concrete candidate:
\begin{equation}\label{eqn_equation10_typeII}
A(\mathbf{s})
=
A_{\mathrm{ht}}+\epsilon A_{\mathrm{sec}}
=
\sum_i h_i(\mathbf{s})V_i
-
\epsilon\sum_a w_a(\mathbf{s})F_a,
\qquad
0<\epsilon\lll 1.
\end{equation}

Choose the scales hierarchically.
For every \((-1)\)-flop,
\[
A_{\mathrm{ht}}\cdot C\neq0,
\]
so for sufficiently small \(\epsilon\),
\[
\operatorname{sgn}(A\cdot C)
=
\operatorname{sgn}(A_{\mathrm{ht}}\cdot C).
\]
For every internal \((-2)\)-flop,
\[
A_{\mathrm{ht}}\cdot C=0,
\]
while completeness of the bubbling marking gives
\[
A_{\mathrm{sec}}\cdot C\neq 0.
\]
Hence
\[
A\cdot C\neq 0.
\]
Therefore
\begin{equation}\label{eqn_equation11_typeII}
    A(\mathbf{s})\cdot C\ne0
\quad
\text{for every Type II flopping curve }C.
\end{equation}
And across the flop,
\[
[C]\longleftrightarrow-[C^+],
\]
so
\[
A\cdot C>0
\quad\Longleftrightarrow\quad
A^+\cdot C^+<0.
\]
Thus exactly one adjacent model lies in the \(A\)-positive chamber.

The section weights can also be made lexicographic. If the bubbling sections have an intrinsic ordering
\[
s_1\prec s_2\prec\cdots\prec s_N,
\]
take formally
\[
w_1\gg w_2\gg\cdots\gg w_N>0,
\]
or numerically
\[
w_a=\delta^a,\qquad0<\delta\lll 1.
\]
For a fixed finite collection of flop walls, every nonzero expression
\[
\sum_a m_aw_a
\]
then has a well-defined leading term, so it cannot vanish.

We need to get a canonical choice, rather than numerical existence. If two bubbling sections are completely indistinguishable and the moduli problem only remembers them as an unordered set, any invariant construction may force
\[
w_a=w_b,
\]
and then the wall
\[
C=e_a-e_b
\]
is invisible. So to get a separated marked moduli problem you really want the bubbling sections to retain their intrinsic level/order/valuation coming from the smoothing. If they do, the Type II calculation is actually stronger than the Type III one: the rooted interval gives a canonical orientation, and there is no \((M2)\) ambiguity.
For the index-one-cover Deligne-Mumford stack, the same calculation should be done upstairs on the Gorenstein cover. 

Therefore, there are only finite many relevant codimension-one flop walls, and we get the following:

For every canonical smoothing of  bubble-tree index-one-cover Deligne-Mumford stacks with bubbling sections $\mathbf{s}$, in the type II case,   $A(\mathbf{s})$ belongs to the interior of a unique Mori chamber, and that chamber is invariant under base change and isomorphism of marked smoothing families. 

This means that we can put the stability condition to be $A(\mathbf{s})\cdot C>0$.  Then the above analysis implies that the extension of the bubble-tree family from $\Delta^{\circ}$ to $\Delta$ is unique. 
\end{proof}

Lemma \ref{lem_separatedness_M012_typeIII} and Lemma \ref{lem_separatedness_M012_typeII} imply that we need introduce one numerical class $A_{\mathbf{s}}$ in 
Definition \ref{defn_marked_bubble_object}.   Let 
\[
\mathfrak M^{\Bub,\mathrm{sec}}_{\tau,\sigma, A_{\mathbf{s}}}
\longrightarrow\Sch_{\cc}.
\]
be the resulting category fibered in groupoids in Definition \ref{defn_marked_bubble_object} by adding another numerical class $A_{\mathbf{s}}$ as part of the polarization. 
Thus, an object of type $(\tau,\sigma, A_{\mathbf{s}})$ over a scheme $T$ is a polarized admissible object
\[
(f:\mathscr X^{\Bub}\to T,c,\mathscr L,\mathscr H,\xi)
\]
together with ordered sections
\[
\mathbf s=(s_1,\ldots,s_r),\qquad s_i:T\longrightarrow\mathscr U_{v_i},
\]
satisfying the incidence conditions $\mathcal A$.  Here the polarization is $\mathscr{L}+A_{\mathbf{s}}$. 

\begin{thm}\label{prop_marked_stack_construction}
The fiber category $\mathfrak{M}^{\Bub,\mathrm{sec}}_{\tau,\sigma, A_{\mathbf{s}}}$ is a separated finite-type Deligne--Mumford stack.  The morphism forgetting the sections,
\[
F:\mathfrak{M}^{\Bub,\mathrm{sec}}_{\tau,\sigma, A_{\mathbf{s}}}
\longrightarrow M^{\Bub,\mathrm{pol}}_\tau,
\]
is representable.  If the entire labeled incidence condition is open in the relative smooth locus over the chosen base stratum, then $F$ is smooth of relative dimension $2r$.

Moreover, there exists a proper morphism 
$$
\mathbf{c}_{\mathrm{sec}}: \mathfrak{M}^{\Bub,\mathrm{sec}}_{\tau,\sigma, A_{\mathbf{s}}}\to M_N
$$
if $N$ is sufficiently large and divisible. 
\end{thm}
\begin{proof}
Similar proof as in Theorem \ref{thm_moduli_stack} shows that $\mathfrak{M}^{\Bub,\mathrm{sec}}_{\tau,\sigma, A_{\mathbf{s}}}$ is a finite type Deligne-Mumford stack.
It is worth mentioning that 
the automorphism group of a marked object is the subgroup of the automorphism group of its underlying polarized object that fixes every labeled section.  It is therefore finite.

The forgetting morphism $F$ is representable since adding marking sections reduces the automorphism group of the bubble-tree index-one-cover Deligne-Mumford stacks.
Since the universal family $p: \mathscr{X}^{\Bub}\to M^{\Bub,\mathrm{pol}}_\tau$ is smooth of relative dimension two, the open configuration morphism $F$ is smooth of relative dimension $2r$.

If $N$ is sufficiently large and divisible, the KSBA moduli stack $M_N$ is proper. 
The contraction  morphism $\mathbf{c}:  \mathfrak{M}^{\Bub,\mathrm{sec}}_{\tau,\sigma, A_{\mathbf{s}}}\to M_N$ is proper from the valuative criterion for properness which we address now. 
    Consider the following diagram 
    \[
\xymatrix{
\spec(K^\prime)\ar[r]\ar[d]&\spec(K)\ar[r]\ar[d]&  \mathfrak{M}^{\Bub,\mathrm{sec}}_{\tau,\sigma, A_{\mathbf{s}}}\ar[d]^{f}\\
\spec(R^\prime)\ar[r]\ar@{-->}[urr]^{\phi}&\spec(R)\ar@{-->}[ur]^{\phi}\ar[r]^{\varphi}& M_N
}
\]
where $R$ is a valuation ring with field of fractions $K$, and residue field $\cc$. We can take $R=\cc[\![t]\!]$ and $K=\cc(\!(t)\!)$. 
The morphism $\varphi: \spec(R)\to M_N$  corresponds a flat $\qq$-Gorenstein family  $f: \sX\to \spec(R)$  of s.l.c. surfaces.
It also represents a flat family $f: \XX\to \spec(R)$ of index-one-cover Deligne-Mumford stacks. 

Suppose the family \( f: \mathcal{X} \to \spec(R) \) is a smoothing whose central fiber \( \mathcal{X}_0 = f^{-1}(0) \) contains no simple elliptic, cusp, or degenerate cusp singularities of embedded dimension at least 6. Then, by construction, the lifting \( f: \mathscr{X}^{\Bub} \to \spec(R) \) coincides with the original family \( f \), plus the labeled sections,  yielding a morphism \( \phi: \spec(R) \to  \mathfrak{M}^{\Bub,\mathrm{sec}}_{\tau,\sigma, A_{\mathbf{s}}}\).

Now suppose the central fiber \( \mathcal{X}_0 \) does contain such singularities (simple elliptic, cusp, or degenerate cusp) of embedded dimension \( \geq 6 \). Following the constructions in \S\ref{defn_bubble_tree_surface} and \S\ref{subsec_bubble_tree_smoothing}, we may build the bubble-tree index-one-cover Deligne-Mumford stack \( \mathscr{X}_0^{\Bub} \) from \( \mathcal{X}_0 \). By Theorem \ref{thm_bubble_tree_DM_smoothing}, there exists a smoothing \( f: \mathscr{X}^{\Bub} \to \spec(R) \) with central fiber \( f^{-1}(0) = \mathscr{X}_0^{\Bub} \). The tree \( \sT \) used to construct \( \mathscr{X}_0^{\Bub} \) is determined by the asymptotic behavior of the bubbling metric tangent cones in K\"ahler geometry, as shown in \cite{Sun}. 
Further, there exists a global polarization $\mathscr{L}$ over the family satisfying the conditions in the construction of the moduli stack. 
This construction again produces a morphism \( \overline{\phi}: \spec(R) \to M_\tau^{\Bub,\mathrm{pol}} \).   Adding the labeled sections, it gives a morphism 
\( \phi: \spec(R) \to \mathfrak{M}^{\Bub,\mathrm{sec}}_{\tau,\sigma, A_{\mathbf{s}}} \).
We prove that the morphism $\phi$ is unique. 
Suppose that there are  two families
\[
(f_1:\mathscr X_1^{\Bub}\to T,c_1,\mathscr L_1,\mathscr H_1,\xi_1; \mathbf s_1=(s_1,\ldots,s_r))
\]
and 
\[
(f_2:\mathscr X_2^{\Bub}\to T,c_2,\mathscr L_2,\mathscr H_2,\xi_2; \mathbf s_2=(s_1,\ldots,s_r)),
\]
so that $f_1, f_2$ are in one $\mathbb{S}$-equivalence class, i.e., they are  related by a type $(M0)$ $(M1)$ or $(M2)$-flop and the generic fibers are the same as the former family 
$f: \mathscr{X}^{\Bub}\to \spec(R)$,  then from Lemma \ref{lem_separatedness_M012_typeIII}, and Lemma \ref{lem_separatedness_M012_typeII},  we can choose different sections $\mathbf{s}_1$ and $\mathbf{s}_2$ to distinct these two families so that the stability condition 
$$A_{\mathbf{s}}\cdot C>0$$ 
for the flopping curves $C$  make only one family satisfying the stability condition.  Hence  
the morphism $\phi: \spec(R) \to \mathfrak{M}^{\Bub,\mathrm{sec}}_{\tau,\sigma, A_{\mathbf{s}}}$ is unique.

Suppose the family $f: \mathcal{X} \to \spec(R)$ is a deformation of slc surfaces containing bad singularities of type simple elliptic, cusp, or degenerate cusp, with embedded dimension at least $6$. In this situation, we use the deformation of the bubble-tree index-one-cover Deligne-Mumford stack $\mathscr{X}_t^{\Bub}$ for any $t \in \spec(R)$.

Such deformations exist; see \cite{DJ} for a study of finite group equivariant deformations of simple elliptic and cusp singularities.

Consequently, after a possible base change $\spec(R') \to \spec(R)$, there exists a morphism $\spec(R') \to \mathfrak{M}^{\Bub,\mathrm{sec}}_{\tau,\sigma, A_{\mathbf{s}}}$. This implies that the morphism $\mathbf{c}_{\mathrm{sec}}: \mathfrak{M}^{\Bub,\mathrm{sec}}_{\tau,\sigma, A_{\mathbf{s}}} \to M_N$ is proper.
\end{proof}

\begin{rmk}
The $(M0), (M1), (M2)$ flops cause trouble for the separatedness of the moduli stack $M^{\Bub,\mathrm{pol}}_\tau$, since they produce double points.  After adding sections, families in the $\mathbb{S}$-equivalence class correspond to different points in the moduli stack
$\mathfrak{M}^{\Bub,\mathrm{sec}}_{\tau,\sigma, A_{\mathbf{s}}}$.
A similar idea can be applied to the moduli compactification of K3 surfaces by using bubbling sections on type II and III degenerations. 
\end{rmk}

\subsection{Deformation theory of the sections}

Let $x=(\mathscr X^{\Bub},\mathbf s)$ be a geometric marked point, suppressing the polarization, contraction, and pre-smoothing direction from the notation.  On the open configuration locus, with no additional closed incidence constraints, moving the markings is unobstructed relative to the unmarked object and has tangent space
\[
T_{F,x}=\bigoplus_{i=1}^r T_{\mathscr X^{\Bub},s_i}.
\]
Write $T^j_\tau$ for the deformation spaces of the underlying object with the homotopy-fiber conditions imposed by $(c,\mathscr L,\mathscr H,\xi)$.  The low-degree part of the deformation sequence is
\[
0\to\operatorname{InfAut}(x)\to T^0_\tau
\longrightarrow\bigoplus_iT_{\mathscr X^{\Bub},s_i}
\longrightarrow\operatorname{Def}(x)
\longrightarrow T^1_\tau
\to0,
\]
and the obstruction is the pullback of the class in $T^2_\tau$.  Thus sections add tangent directions and reduce infinitesimal automorphisms, but introduce no new relative obstruction as long as they remain in $\mathscr U$.

\begin{prop}[Induced obstruction theory]\label{prop_marked_pot}
Suppose $E^\bullet\to \ll^{\bullet}_{\mathfrak M^{\Bub,\mathrm{pol}}_\tau}$ is the perfect obstruction theory described in Section \ref{subsec_pot_tautological_intro}, suppose $F$ is smooth, and assume that this theory admits a compatible lift across the cotangent transitivity triangle for $F$.  Then
\[
E^\bullet_{\mathrm{sec}}
\longrightarrow \ll^{\bullet}_{\mathfrak{M}^{\Bub,\mathrm{sec}}_{\tau,\sigma, A_{\mathbf{s}}}}
\]
is a compatible absolute perfect obstruction theory, where $E^\bullet_{\mathrm{sec}}$ fits into a distinguished triangle
\[
F^*E^\bullet\longrightarrow E^\bullet_{\mathrm{sec}}
\longrightarrow L_F\longrightarrow F^*E^\bullet[1].
\]
Locally, after choosing a splitting of the smooth cotangent triangle, it is represented by $F^*E^\bullet\oplus\Omega_F$.  Its virtual dimension is
\[
\operatorname{vdim}\mathfrak{M}^{\Bub,\mathrm{sec}}_{\tau,\sigma, A_{\mathbf{s}}}
=\operatorname{vdim} M^{\Bub,\mathrm{pol}}_\tau+2r.
\]
\end{prop}
\begin{proof}
The relative cotangent complex of the smooth morphism $F$ is the vector bundle $\Omega_F$ in degree zero.  The assumed compatibility with the cotangent transitivity triangle gives $E^\bullet_{\mathrm{sec}}$.  Its virtual rank is the virtual rank of $E^\bullet$ plus the rank $2r$ of $L_F$.
\end{proof}

\begin{rmk}
Smoothness of $F$ makes the relative deformation theory unobstructed, but it does not by itself choose a global splitting of the cotangent triangle.  This is why Proposition \ref{prop_marked_pot} states the compatibility requirement explicitly rather than claiming that the displayed local direct sum is canonical.
\end{rmk}

\subsection{Comparison with the former moduli space and with metric bubbles}

The three relevant morphisms are
\[
\xymatrix{
\mathfrak{M}^{\Bub,\mathrm{sec}}_{\tau,\sigma, A_{\mathbf{s}}}
\ar[r]^-F\ar[dr]_{\mathbf c_{\mathrm{sec}}}&
M^{\Bub,\mathrm{pol}}_\tau\ar[d]^{\mathbf c}\\
&M_N,
}
\qquad \mathbf c_{\mathrm{sec}}=\mathbf c\circ F.
\]
The former stack remembers the bubble components, polarization, rigidifying divisor, and pre-smoothing direction.  The new stack remembers in addition ordered centers on specified bubble components. 
There is also a morphism to the  KSBA moduli stack
\[
(\mathscr X^{\Bub},\mathbf s)\longmapsto
(\mathscr X,c\circ\mathbf s).
\]
 Distinct sections on different bubble components may have the same image on $\mathscr X$; hence this morphism forgets precisely the information about which algebraic bubble carries the center.

\begin{con}[Comparison with Song Sun's pointed bubble tree]\label{con_marked_metric_comparison}
Let $(\mathcal X\to\Delta,\bar s_1,\ldots,\bar s_r)$ be a pointed $\qq$-Gorenstein smoothing, and assume that its negative K\"ahler--Einstein metrics admit a finite pointed bubble tree at the marked sections.  Assume also that every non-conical vertex has a compatible algebraic compactification and that the neck incidence agrees with an admissible decorated dual complex in $\mathcal I$.  After choosing the auxiliary polarization, divisor, first-order deformation, and compatible lifted sections, finite base change and semistable reduction are expected to produce a point of $\mathfrak{M}^{\Bub,\mathrm{sec}}_{\tau,\sigma, A_{\mathbf{s}}}$.  Conversely, an effective smoothing of a marked algebraic object determines the centers and incidence labels of the metric bubble tree whenever the required analytic compactness and uniqueness statements hold.
\end{con}

\begin{rmk}
The conjecture compares centers, vertex labels, and incidence relations.  To recover the complete pointed metric object one must enrich the algebraic datum by scale parameters, or, if a comparison theorem relates them, by valuation data recording degeneration rates.  Such an enrichment is closer to a logarithmic or non-Archimedean moduli problem and is not part of $\mathfrak{M}^{\Bub,\mathrm{sec}}_{\tau,\sigma, A_{\mathbf{s}}}$.
\end{rmk}

\section{A differential-geometric framework for cusp bubbling}\label{sec_discussion}

This section explains the differential-geometric motivation for the bubble-tree construction and, in particular, its relation to cusp singularities. We separate the established analytic input from the additional compactness and identification statements that remain conjectural.

\subsection{Bubble limits}
Let $(X_j,L_j,\omega_j,p_j)$ be a sequence of polarized K\"ahler--Einstein manifolds, with $\omega_j\in 2\pi c_1(L_j)$ and
\[
\operatorname{Ric}(\omega_j)=c_j\omega_j,\qquad c_j\in\{-1,0,1\}.
\]
When $c_j=0$, we also assume that $K_{X_j}$ is trivial. A degeneration arising from a one-parameter smoothing of a stable surface supplies such a sequence on the smooth fibers. The converse, however, is not automatic: an arbitrary metric sequence need not arise from an algebraic one-parameter family.

Suppose that the sequence satisfies the volume non-collapsing condition that there exists a $\kappa>0$ such that
$$
\text{Vol}(B(p_j,1))\ge \kappa
$$
for all $j$. This is the non-collapsing K\"ahler--Einstein setting established for polarized K\"ahler--Einstein Fano manifolds, corresponding to $c_j=1$. After passing to a subsequence, Donaldson--Sun \cite{DS1} proved that the Gromov--Hausdorff limit $(X_{\infty}, d_{\infty}, p_{\infty})$ is a normal projective variety with log-terminal singularities, where $d_{\infty}$ is induced by a weak K\"ahler--Einstein metric $\omega_{\infty}$ on $X_{\infty}$. They also proved in \cite{DS2} that there is a unique metric tangent cone $\sC_0(X_{\infty})$ at $p_{\infty}$. This metric tangent cone is a normal affine variety endowed with a weak Calabi--Yau cone metric.
The algebraic structure of $\sC_0(X_{\infty})$ can be obtained from a $2$-step degeneration from the local ring $\sO_{(X_{\infty}, p_{\infty})}$, and \cite{DS2} conjectured, which was proved in \cite{LWX}, that the cone is only determined by the germ $(X_{\infty}, p_{\infty})$.

From \cite{Sun}, tangent cones can be regarded as only the {\em first order} approximation
of the singular behavior of solutions to an elliptic nonlinear PDE, and $\sC_0(X_{\infty})$ does not encode the information on how the singularity is formed. Then \cite{Sun} considered the rescaled spaces $(X_j, \lambda_j\omega_j, p_j)$ for all possible sequences of integers $\lambda_j\ge 1$. When $\lambda_j\to \infty$, then passing to a subsequence we obtain a non-compact pointed Gromov-Hausdorff limit $(Z, d_Z, p_Z)$, which is naturally an affine variety with a weak Calabi-Yau metric $\omega_Z$, which is asymptotic to a unique metric cone $\sC_{\infty}(Z)$ at infinity. There is also a unique metric tangent cone $\sC_0(Z)$ at $p_Z$.
Hence we get the definition of bubble limit in \cite{Sun}.

\begin{defn}\label{defn_bubble_limit}
 A {\em bubble limit} $(Z,\omega_Z, p_Z)$ is a subsequential pointed Gromov-Hausdorff limit of the sequence $(X_j, \lambda_j\omega_j, p_j)$ for some $\lambda_j\ge 1$.
\end{defn}

A bubble limit is endowed with both an algebraic structure and a metric structure. Let
$\mathfrak{B}$ be the class of all bubble limits modulo equivalence relation by holomorphic isometries preserving the base points. We have
$X_{\infty}\in \mathfrak{B}$ is the only compact bubble limit. Also
$\sC_0(X_{\infty})\in \mathfrak{B}$. For any non-compact $Z\in \mathfrak{B}$, both the tangent cone $\sC_0(Z)$ and the asymptotic cone $\sC_{\infty}(Z)$ are in $\mathfrak{B}$.

A cone limit is a bubble limit $(\sC,O)$ such that $\sC$ is a metric cone with vertex $O$.
A {\em minimal bubble} is a bubble limit which is not a cone limit but is asymptotic to $\sC_0(X_{\infty})$. \cite[Theorem 1.11]{Sun} proved the existence of minimal bubbles, and
\cite[Corollary 1.13]{Sun} showed that, after fixing a subsequence and passing to a further subsequence, the bubble limits of $(X_j, \omega_j, p_j)$ modulo metric scaling form a finite set $\{(Z_{\alpha}, q_{\alpha})\}_{\alpha=1}^{2m}$ with $\text{Vol}_0(Z_{\alpha})$ strictly increasing, where
$(Z_1,q_1)=(X_{\infty},p_{\infty})$, $(Z_{2m}, q_{2m})=(\cc^n, 0)$, and for all $\alpha$, $(Z_{2\alpha}, q_{2\alpha})$ is a cone limit, and
$\sC_0(Z_{2\alpha-1})=Z_{2\alpha}=\sC_{\infty}(Z_{2\alpha+1})$. The integer $m$ is uniformly bounded and only depends on the dimension $n$ and $\Xi\ge 0$ for $\text{diam}(X_j, \omega_j)\le \Xi$.

\subsection{General type surfaces}\label{subsec_general_type_Bubb}

Smooth surfaces with ample canonical bundle admit negative K\"ahler--Einstein metrics. Fix a projective smoothing $f:\sX\to\Delta$ with smooth canonically polarized general fiber, and choose points and a sequence of parameters tending to zero. A general sequence in moduli need not itself lie on one algebraic curve; the one-parameter family is part of the hypothesis.
Suppose that the central fiber $\sX_0$ is an slc stable surface which contains simple elliptic, cusp and degenerate cusp singularities. Simple elliptic and cusp singularities are log-canonical surface singularities which are normal singularities, while degenerate cusp singularities are non-normal singularities.

The algebraic family converges to the stable surface $X_\infty=X_0$. Analytically, the negative K\"ahler--Einstein metrics converge locally smoothly on compact subsets of the regular locus, under the standard hypotheses in \cite{BG}. Near a log-canonical, non-klt point, the limiting metric is complete and places the point at infinite distance. Consequently, one should regard the metric completion and the algebraic KSBA limit as related but distinct objects; it is generally inaccurate to identify the pointed Gromov--Hausdorff limit with the projective surface without specifying the normalization and base points.
Note that \cite{BG} constructed K\"ahler--Einstein metrics on the regular locus of an slc variety $X$; and \cite{DFS} constructed K\"ahler--Einstein metrics
near an isolated log-canonical singularity.

If curvature concentrates near the cusp, a suitable blow-up may produce a complete non-compact Ricci-flat K\"ahler four-manifold. Its universal cover is then hyperk\"ahler, and, under the usual finite-energy and asymptotic assumptions, it is natural to compare it with a gravitational instanton. Sun--Zhang establish the corresponding bubble-tree picture in the hyperk\"ahler K3 setting \cite{SZ}. Their theorem does not directly imply the analogous statement for negative K\"ahler--Einstein surfaces.

In the KSBA setting, the metric develops a complete end near a simple elliptic or cusp singularity. Determining whether curvature blow-ups exist and classifying their limits are separate analytic problems. For simple elliptic singularities, \cite{FHJ} constructs an $ALH^*$ Ricci-flat bubble whose rotated complex structure is $Y\setminus E$, where $Y$ is a smooth del Pezzo surface and $E\in|-K_Y|$ is a smooth elliptic curve. This result motivates the del Pezzo components in the construction of \S\ref{subsec_bubble_tree}.

Cusp singularities are subtler. By Looijenga's conjecture, a cusp is smoothable precisely when the exceptional cycle of the dual cusp occurs as an anticanonical cycle on a rational surface. This algebro-geometric criterion motivates the following conditional differential-geometric statement. It is deliberately formulated so that the analytic steps not presently supplied by the algebraic construction are visible in the hypotheses.

\begin{prop}\label{prop_conditional_cusp_bubbling}
Let $g_j$ be negative K\"ahler--Einstein metrics on smooth canonically polarized surfaces, with $\operatorname{Ric}(g_j)=-g_j$. Suppose $\lambda_j\to\infty$ and $(X_j,\lambda_jg_j,p_j)$ converge smoothly in the pointed Cheeger--Gromov sense to a complete smooth four-manifold $(Z,g,p)$. Then $g$ is Ricci-flat and K\"ahler, after passing to a subsequence for the parallel complex structures. Its universal cover is hyperk\"ahler. If $|\operatorname{Rm}(\lambda_jg_j)|(p_j)=1$, the limit is nonflat.
\end{prop}
\begin{proof}
The Einstein constant of $\lambda_jg_j$ is $-\lambda_j^{-1}$, so smooth convergence implies $\operatorname{Ric}(g)=0$. In the convergence charts, the complex structures are parallel orthogonal tensors. Their values at the base point have a convergent subsequence; parallelness and the smooth convergence of the connections then give a limiting parallel complex structure. In real dimension four the restricted holonomy of a Ricci-flat K\"ahler metric lies in $SU(2)$, and the simply connected cover is hyperk\"ahler. The curvature normalization passes to the limit.
\end{proof}

This proposition proves the differential-geometric conclusion available directly from rescaling. It does not establish the convergence hypothesis, finite curvature energy, global hyperk\"ahler structure on $Z$ itself, or a rational compactification. The earlier version assumed the desired compactification and consequently did not prove the essential Type III identification. The following subsections give a more precise approach.

\subsection{What a Type III model can determine}\label{subsec_cusp_reconstruction}

Fix a smoothing component of a normal cusp and a compatible Type III anticanonical model as in \cite{Engel,FE}. One must distinguish the rational pair on a general fiber, the exceptional cycle of the dual cusp, and the rational components in the Type III central fiber. In the classical construction the distinguished component is an Inoue--Hirzebruch surface; the entire analytic central fiber is not a projective union of rational surfaces. Any replacement of that distinguished component by a resolved global KSBA component requires a separate gluing and contraction argument.

The Type III data records component pairs $(V_i,D_i)$, double curves, triple points, and an integral-affine structure on the dual complex. The dual complex is two-dimensional and may be a sphere; a metric bubble tree is instead organized by nested concentration scales. These are different combinatorial objects. Recovering the latter requires choosing regions of the former and proving how their sizes compare in the degenerating metrics.

\begin{prop}\label{prop_cusp_hyperbolic}
For a cusp with resolution cycle of length $n\geq3$ and self-intersections $(-b_1,\ldots,-b_n)$, where $b_i\geq2$ and at least one $b_i>2$, choose an orientation convention and put
\[
A_C=\prod_{i=1}^{n}\begin{pmatrix}b_i&-1\\1&0\end{pmatrix}.
\]
The link is a torus bundle with this monodromy, up to the usual changes of basis and inversion from orientation conventions. In particular $\operatorname{tr}(A_C)>2$: the link has hyperbolic torus monodromy and Sol geometry, not unipotent torus monodromy.
\end{prop}
\begin{proof}
Plumbing the normal disk bundles of the resolution cycle gives the displayed change of basis for successive torus cross-sections; see \cite{KasuyaCusp}. For the trace claim, conjugation by $P=\left(\begin{smallmatrix}1&1\\1&0\end{smallmatrix}\right)$ gives
\[
P^{-1}\begin{pmatrix}b&-1\\1&0\end{pmatrix}P
=\begin{pmatrix}1&1\\b-2&b-1\end{pmatrix}.
\]
These matrices have nonnegative entries, dominate $U=\left(\begin{smallmatrix}1&1\\0&1\end{smallmatrix}\right)$ entrywise, and at least one has positive lower-left entry. Replacing all other factors by $U$ shows the trace of the product is strictly greater than two. The torus-bundle description and its Sol geometry follow. Cycles of length one or two require their corresponding nodal or multiple-intersection plumbing conventions.
\end{proof}

The term ``parabolic subgroup'' in a Hilbert modular cusp description must not be confused with the conjugacy class of this $SL(2,\mathbb Z)$ monodromy. In particular, a single Kodaira $I_b$ end, whose monodromy is unipotent, cannot be identified with the whole cusp link as a torus bundle. An individual bubble may have a different end because its neck encloses only part of the degenerating region.

Suppose a region in a proposed affine base is a disk with finitely many singular points and carries a rank-two integral local system away from them. Choose based loops and parallel-transport paths. Then its outer monodromy factors as
\[
A_{\mathrm{out}}=A_1\cdots A_k.
\]
This is the image under the monodromy representation of the fundamental-group boundary relation. It is a necessary compatibility test for a proposed decomposition into candidate bubble regions, once one has identified this local system with the relevant torus fibration. That identification is an extra assumption; affine monodromy and the link monodromy should not be equated without it.

The factorization alone does not determine the factors or the scales. For example
\[
\begin{pmatrix}1&a\\0&1\end{pmatrix}
\begin{pmatrix}1&0\\b&1\end{pmatrix}
=\begin{pmatrix}1+ab&a\\b&1\end{pmatrix},\qquad a,b>0,
\]
has hyperbolic trace $2+ab$ although each factor is unipotent. This matrix identity illustrates the possibility of different local and outer monodromies; it is not a realization theorem for a positive geometric fibration.

\begin{prop}\label{prop_component_charge}
Let $V$ be a smooth rational surface and $D=D_1+\cdots+D_n\in|-K_V|$ an SNC cycle with $n\geq3$. Then the charge 
\[
Q(V,D):=e(V\setminus D)=12-K_V^2-n
=12-3n-\sum_iD_i^2.
\]
\end{prop}
\begin{proof}
For a rational surface Noether's formula gives $e(V)=12-K_V^2$. The cycle has $e(D)=2n-n=n$, and $K_V^2=D^2=\sum D_i^2+2n$. Euler characteristic with compact supports is additive and agrees here with the ordinary Euler characteristic, giving the formula.
\end{proof}
This charge is a topological quantity, not a curvature-energy formula. A Gauss--Bonnet calculation on a noncompact bubble includes boundary terms. Thus neither $Q=0$ nor a prescribed charge distribution alone proves that a Type III component is flat or that it carries one nonflat metric vertex.

\subsection{A precise conditional analytic route}\label{subsec_cusp_analytic_route}

Here is a usable sufficient criterion for identifying a candidate suggested by a Type III model. It isolates estimates on metrics, rather than assuming in advance that an unknown bubble has the desired compactification.

\begin{thm}\label{thm_cusp_estimate_identification}
Let $(Z,g_*,J_*,p)$ be a complete Ricci-flat K\"ahler surface chosen as a candidate from the algebraic data. Let $g_j$ be negative K\"ahler--Einstein metrics on the smooth fibers of a fixed smoothing, and let $\lambda_j\to\infty$. Suppose there are embeddings $\Phi_j:U_j\to X_j$ from an exhaustion $U_j\nearrow Z$, taking $p$ to $p_j$, such that:
\begin{enumerate}
\item $\Phi_j^*J_j\to J_*$ smoothly on compact subsets. The pulled-back metrics $h_j=\Phi_j^*(\lambda_jg_j)$ are uniformly elliptic and bounded in $C^k$ on each compact subset for every $k$.
\item The distance in $(X_j,\lambda_jg_j)$ from $p_j$ to the complement of $\Phi_j(U_j)$ tends to infinity.
\item Every subsequential local limit belongs to a fixed class $\mathcal C$ of complete metrics on $Z$, specified by asymptotics, periods, and normalization. Within $\mathcal C$, $g_*$ is the unique Ricci-flat K\"ahler metric, up to the specified pointed holomorphic equivalences.
\end{enumerate}
Then $(X_j,\lambda_jg_j,p_j)$ converges to $(Z,g_*,p)$ in pointed smooth topology, modulo those equivalences. If the limiting normalized curvature at $p$ is nonzero, this is a nonflat bubble.
\end{thm}
\begin{proof}
The compact-set bounds and ellipticity give smooth subsequential convergence by diagonal extraction. The limiting complex structure is $J_*$ and the limiting metric is K\"ahler. Rescaling the Einstein equation gives Ricci-flatness. The second hypothesis ensures that bounded pointed balls are captured by the charts and prevents an artificial finite-distance boundary in the limit. The third hypothesis supplies completeness and places the limit in the uniqueness class. Every subsequential limit is therefore $g_*$ up to the stated equivalences, giving convergence of the full sequence modulo them. Nonvanishing curvature passes to the limit.
\end{proof}

For a Type III application, the work is to construct these charts and prove the hypotheses. In a normal-crossing chart $xy=t^\ell u$, with $u$ a unit, the integer $\ell$ and the section's orders of vanishing give logarithmic degeneration data. They do not yet give the factor $\lambda_j$. One needs estimates relating these orders and the K\"ahler periods to distances and curvature. Locally one can write a pulled-back K\"ahler form as
\[
\omega_j^{\mathrm{scaled}}=\omega_{*,j}
+\sqrt{-1}\partial\bar\partial u_j.
\]
Uniform potential, ellipticity, and derivative estimates, together with control on expanding necks, are the analytic input to Theorem \ref{thm_cusp_estimate_identification}. Period and asymptotic matching must also be checked before applying a uniqueness theorem; the $ALH^*$ Torelli theorem \cite{CJL_Torelli} illustrates the need for such additional data.

To recover an entire bubble tree one must repeat this analysis at separated scales, prove that the regions are exhaustive, and control residual neck energy. Under a uniform total $L^2$ curvature bound $E$ and a positive lower bound $\varepsilon_0$ for the energy of each disjoint extracted nonflat region, at most $E/\varepsilon_0$ such regions occur. This follows by additivity on disjoint regions and scale invariance of $\int|\operatorname{Rm}|^2$ in real dimension four. Establishing the energy gap, disjointness, and the absence of further concentration is additional work in a collapsing cusp degeneration.

We have a conjecture for Type III data as a constraint on cusp bubbles.
\begin{con}\label{con_cusp_reconstruction_revised}
Fix a projective smoothing of a cusp, its K\"ahler classes, a section, and a compatible Type III model with period and gluing data. When a finite bubble decomposition exists, its noncollapsed regions should be recoverable from suitable regions of an integral-affine refinement together with the relative metric scales. The monodromy and topology of these regions must satisfy the compatibility tests above. Identification of a candidate metric vertex is governed by Theorem \ref{thm_cusp_estimate_identification}; no assertion is made that every rational component contributes one vertex or that the undecorated central fiber determines a unique metric.
\end{con}

The answer to the reconstruction question is therefore partial but concrete: Type III geometry provides candidates and tests that can rule out proposed bubble pictures. Turning this into a theorem for a specified smoothing requires the metric estimates and matching conditions above. A local cusp smoothing alone does not specify negative K\"ahler--Einstein metrics without a global canonical polarization or appropriate boundary data.

We propose a question.

\begin{question}
 If $p_{\infty}$ is a cusp singularity of high embedding dimension, does a non-collapsed Ricci-flat blow-up limit exist, and is it of type $ALG$ or $ALG^*$? How does its end monodromy encode the duality between a smoothing of the cusp and the existence and deformation theory of a Looijenga pair for the dual cusp cycle?

 Moreover, in view of Proposition \ref{prop_minus-one_form_Looijenga}, are the compactifications of these bubble limits related to Type III degenerations in minus-one form, with every rational component not adjacent to the main index-one-cover Deligne--Mumford stack $\XX_0$ in minus-one form?
\end{question}

\begin{rmk}
After choosing an appropriate complex structure in the hyperk\"ahler family, many $ALG$ and $ALG^*$ gravitational instantons admit compactifications related to rational elliptic surfaces. It is therefore natural to ask whether the cusp bubbles above admit such compactifications and how these are related to Type III degenerations of log Calabi--Yau surfaces. We do not assume that every $ALG$ or $ALG^*$ space is literally the complement of a fiber in the original complex structure.
\end{rmk}

Suppose that an algebraic degeneration has a non-normal central fiber $X_\infty$ with a degenerate cusp $p_\infty$. The smoothing and semistable-reduction construction in \S\ref{sec_deg_cusp} suggests the following metric counterpart.
Here is a conjecture for the degenerate-cusp bubbles.

\begin{con}\label{con_deg_cusp} 
Fix a $\qq$-Gorenstein smoothing of a degenerate cusp and a section through the cusp. Assume that the associated negative K\"ahler--Einstein metrics admit a finite pointed bubble tree and that every non-conical Ricci-flat vertex admits a compatible hyperk\"ahler rotation and algebraic compactification. Then the union of these compactified vertices, with the incidence data induced by the neck regions, is $\mathbb{S}$-equivalent to the relative minimal open Kulikov model constructed in \S\ref{sec_deg_cusp}. A converse identification would require period, asymptotic, and scale data together with uniqueness for the corresponding metrics.
\end{con}

\subsection{Comparison with the K3 surface case}

\subsubsection{K3 surfaces}
There is an interesting comparison with the moduli space of K3 surfaces. Let $\overline{\MM}$ be the full Gromov-Hausdorff compactification of K3 surfaces, and $\MM$ the partial compactification of moduli of smooth K3 surfaces by adding the non-collapsing Gromov-Hausdorff limits. In algebraic geometry $\MM$ parametrizes K3 surfaces with klt singularities.

From \cite[Theorem 1.1]{SZ}, any collapsed limit in $\overline{\MM}\setminus \MM$ is one of the following three cases:
\begin{enumerate}
 \item in real dimension $3$, a flat orbifold $\ttt^3/\mu_2$;
 \item in real dimension $2$, a singular special K\"ahler metric on $\mathbb{S}^2$ with local integral monodromy;
 \item in real dimension $1$, a unit interval.
\end{enumerate}

\cite[Theorem 1.2]{SZ} proved that all the non-flat gravitational instantons are the six $ALE, ALF, ALG, ALG^*, ALH, ALH^*$-spaces. In algebraic geometry, \cite{AE23} constructed the moduli space of polarized K3 surfaces. There are three types of degeneration of K3 surfaces.

Type I: smooth K3 surfaces degenerate to K3 surfaces with ADE singularities. This is the non-collapsing case discussed above.

Type II: smooth K3 surfaces, or K3 surfaces with ADE singularities, degenerate to a chain of surfaces whose end components are rational and whose intermediate components are elliptic ruled. This is Case (3), where the collapsing limit is a unit interval. It is analogous to the construction of bubble-tree surfaces for simple elliptic singularities in \S\ref{subsec_bubble_tree}.

Type III: smooth K3 surfaces, or K3 surfaces with ADE singularities after simultaneous resolution, degenerate to normal-crossing unions of rational surfaces whose dual complex is a sphere $\mathbb{S}^2$.
This is Case (2) above for the collapsing limit in \cite{SZ}.

Case (1) does not happen for the moduli space of polarized K3 surfaces in \cite{AE23}.

From \cite[Theorem (0.4)]{FS} (Proposition \ref{prop_minus-one_form_Looijenga}), for a Type III degeneration
$$
f: \mathscr{X}\to \Delta
$$
of K3 surfaces, we can always apply type $(M1)$ and $(M2)$ birational modifications so that in $\mathscr{X}_0$ each double curve $D_{ij}$ has self-intersection number $D_{ij}^2=-1$ if it is a smooth rational curve $\pp^1$; while $D_{ij}^2=1$ on the smooth component side, and $-1$ on another side if $D_{ij}$ is a rational nodal curve.

The relation between a Type III algebraic degeneration and a two-dimensional collapsed metric limit requires the polarization and metric normalization to be specified. The Type III label alone does not imply that every nonflat bubble is ALG or $ALG^*$. Hein constructs gravitational instantons from complements of suitable singular fibers in rational elliptic surfaces \cite{Hein}. After a hyperk\"ahler rotation, one may therefore expect a compactification by a Looijenga pair in minus-one form. Establishing such a compactification, and proving its compatibility with the incidence data of the complete bubble tree, are precisely the missing steps recorded in Conjecture \ref{con_uniqueness}.

\begin{example}[The degree-six del Pezzo pair and an $I_6$ $ALH^*$ model]\label{ex_dp6_I6}
Let
\[
Y=\operatorname{Bl}_{q_1,q_2,q_3}\pp^2,
\]
where $q_1,q_2,q_3$ are the three torus-fixed points. This is the toric del Pezzo surface of degree six. Its toric boundary
\[
D=D_1+\cdots+D_6\in|-K_Y|
\]
is a cycle of six smooth rational curves, and $D_i^2=-1$ for every $i$. Thus $(Y,D)$ is a Looijenga pair in minus-one form.

Choose a smooth member $E\in|-K_Y|$ meeting $D$ transversely away from its nodes. Adjunction gives
\[
E\cdot D_i=(-K_Y)\cdot D_i=2+D_i^2=1,
\]
so $E\cap D$ consists of six points $p_i$, one on each $D_i$. Blow up these points:
\[
\rho:V=\operatorname{Bl}_{p_1,\ldots,p_6}Y\longrightarrow Y.
\]
The pencil generated by $D$ and $E$ becomes base-point free on $V$ and defines a genus-one fibration $\pi:V\to\pp^1$. Each exceptional curve $S_i$ is a section and satisfies $S_i^2=-1$. The strict transform
\[
\widetilde D=\widetilde D_1+\cdots+\widetilde D_6
\]
is an anticanonical fiber of Kodaira type $I_6$, because $\widetilde D_i^2=D_i^2-1=-2$. Conversely, contracting the six disjoint sections $S_i$ recovers the pair $(Y,D)$ and changes the boundary self-intersections from $(-2)^6$ to $(-1)^6$.

The open surface $V\setminus\widetilde D$ is an $I_6$ rational-elliptic complement. Hein's construction \cite{Hein} equips it with metrics of volume growth $R^{4/3}$, of type $ALH^*$ in current terminology, rather than ALG. An ALG model instead arises from suitable finite-monodromy fibers; $I_b^*$ fibers give $ALG^*$ models. On the Type III side, $(Y,D)$ has exactly the boundary data required of a component in minus-one form. Indeed, each $D_i$ contains two triple points in a Type III central fiber, and the triple-point formula forces the self-intersection on the adjacent component to be
\[
D_i^2|_{Y'}=-2-D_i^2|_Y=-1.
\]
The triple-point formula is a necessary local compatibility condition, not a proof that this pair occurs in a globally smoothable Type III K3 central fiber. Such an occurrence requires neighboring components, compatible gluings, global $d$-semistability, and a smoothing theorem.

There is also a crucial distinction between the two open surfaces. The contraction $V\to Y$ contracts curves $S_i\setminus\widetilde D\cong\mathbb A^1$ that lie in the open surface. It therefore does not identify $V\setminus\widetilde D$ with $Y\setminus D$. Indeed
\[
e(Y\setminus D)=12-6-6=0,\qquad
e(V\setminus\widetilde D)=12-0-6=6.
\]
The former is the toric open surface $(\mathbb C^*)^2$. The six interior boundary blow-ups change its topology and supply the nonzero charge of the elliptic complement. This example is useful precisely because it shows that a minus-one hexagon by itself does not determine the topology or metric type of a bubble. The pencil and its base points are indispensable extra data. Moreover $D^2=6$, so the hexagon is not a negative-definite cusp resolution cycle; it can only be a local component boundary in a larger Type III construction.
\end{example}
\begin{rmk}
In differential geometry, it is natural to do hyperk\"ahler rotations for $ALG$ and $ALG^*$ spaces.
 It would be interesting to determine the relevant torsor class in the Tate--Shafarevich group of the Jacobian fibration and to relate possible finite-order twists to $ALG^*$ spaces and hyperk\"ahler rotation. No assertion about the order of that group is needed for Example \ref{ex_dp6_I6}.
\end{rmk}

\subsubsection{Case of bubble surfaces}
In the bubble-tree surface $X^{\Bub}$, a cusp $(X,x)$ is replaced by a Type III degeneration of log Calabi--Yau surfaces as in \S\ref{subsec_bubble_tree}. The spherical dual complex resembles the two-dimensional collapsed limit in \cite[Theorem 1.1]{SZ}, but this analogy is not an identification of the two metric limits. From the mirror-symmetry viewpoint, a smoothing of a cusp is associated with a vanishing Lagrangian sphere; see \cite{Aroux}. The corresponding singular integral-affine sphere is the combinatorial input for the Type III construction in \cite{FE}. Clarifying whether the neck decomposition of a negative K\"ahler--Einstein degeneration recovers this integral-affine structure is part of Conjecture \ref{con_uniqueness}.


\begin{thebibliography}{s}

\bibitem{AH} D. Abramovich and B. Hassett, \newblock Stable varieties with a twist, {\em Classification of algebraic varieties, EMS Ser. Congr. Rep.}, Eur. Math. Soc., Z\"urich, 2011, 1-38.
\bibitem{Alexeev2} Alexeev, \newblock Boundedness and $K^2$ for log surfaces, {\em International Journal of Mathematics}, Vol. 05, No. 06, 779-810 (1994).
\bibitem{AM04}V. Alexeev and S. Mori, \newblock Bounding singular surfaces of general type, In: Christensen C., Sathaye A., Sundaram G., Bajaj C. (eds) Algebra, Arithmetic and Geometry with Applications. Springer, Berlin, Heidelberg.
\bibitem{Alexeev_Ann}V. Alexeev, \newblock Complete moduli in the presence of semiabelian group action, {\em Ann. Math.},
 Volume 155 (2002), Issue 3, 611-708.
\bibitem{Alexeev_2023} Alexeev, \newblock Kappa classes on KSBA spaces, {\em Moduli}, 2, e4 (2025) 1-14, arXiv:2309.14842.
\bibitem{AE23} V. Alexeev and P. Engel, \newblock Compact moduli of K3 surfaces, {\em Annals of Math}, 198 , no.2, 727-789 (2023).
\bibitem{AET23}V. Alexeev, P. Engel and A. Thompson, \newblock
Stable pair compactification of moduli of K3 surfaces of degree 2, {\em J. Reine Angew. Math.}, 799(2023), 1-56.
%\bibitem{AT21} V. Alexeev and A. Thompson, \newblock ADE surfaces and their moduli, {\em J. Algebraic Geom.} 30 (2021), no. 2, 331-405.
\bibitem{AAB2024}V. Alexeev, H. Arg\"uz and P. Bousseau, \newblock The KSBA moduli space of stable log Calabi-Yau surfaces, arXiv:2402.15117.

\bibitem{AHL-H}J. Alper, D. Halpern-Leistner, and J. Heinloth, \newblock Existence of moduli spaces for algebraic stacks, {\em Invent. Math.},
Vol. 234 (2023), No. 3, pp. 949-1038.
\bibitem{Aroux}D. Aroux, \newblock
Mirror symmetry and T-duality in the complement of
an anticanonical divisor, {\em Journal of G\"okova Geometry Topology},
Vol. 1 (2007) 51-91.

\bibitem{Artin} M. Artin, \newblock Versal deformations and algebraic stacks, {\em Invent. Math.}, 27, 165- 189 (1974).

\bibitem{ABBDILW}K. Ascher, D. Bejleri, H. Blum, K. DeVleming, G. Inchiostro, Y. Liu, and X. Wang, \newblock Moduli of boundary polarized Calabi-Yau pairs, arXiv:2307.06522.

\bibitem{ABJ}K. Ascher, D. Bejleri and Y. Jiang, \newblock Moduli space of twisted stable maps as moduli of lci covers for fibered surfaces, in preparation.

\bibitem{BF} K. Behrend and B. Fantechi, \newblock The intrinsic normal cone,
 alg-geom/9601010, {\em Invent. Math}. 128 (1997), no. 1, 45-88.
\bibitem{Brieskorn} E. Brieskorn, \newblock Rationale Singularitaten komplexer Fl\"achen, {\em Invent. Math.}, 4, 336-358 (1968).
\bibitem{BG}R. Berman and H. Guenancia, \newblock K\"ahler--Einstein metrics on stable varieties and log canonical pairs, {\em Geom. Funct. Anal.}, Vol. 24 (2014) 1683-1730.

\bibitem{CJL_Torelli} T. Collins, A. Jacob, and Y.-S. Lin, \newblock The Torelli theorem for $ALH^*$ gravitational instantons, arXiv:2111.09260. \url{https://arxiv.org/abs/2111.09260}.
\bibitem{DJ}S. Das and Y. Jiang, \newblock Equivariant deformation of minimally elliptic singularities, arXiv:2510.10253 [math.AG].
\bibitem{DFS}V. Datar, X. Fu, and J. Song, \newblock K\"ahler--Einstein metrics near isolated log-canonical singularities, arXiv:2106.05486.
\bibitem{dBS24} M. de Borbon, C. Spotti, \newblock Some models for bubbling of (log) K\"ahler--Einstein metrics, {\em Annali
Dell’Universita’di Ferrara}, (2024), pp. 1-32.

\bibitem{Donaldson} S. K. Donaldson, \newblock Fredholm topology and enumerative geometry: reflections on some words of Michael Atiyah, {\em Proceedings of 26th G\"okova Geometry-Topology Conference}, 1-31.

\bibitem{DS1}S. Donaldson and S. Sun, \newblock Gromov-Hausdorff limits of K\"ahler manifolds and algebraic geometry, {\em Acta Mathematica}, 213 (1):63-106, 2014.

\bibitem{DS2}S. Donaldson and S. Sun, \newblock Gromov-Hausdorff limits of K\"ahler manifolds and algebraic geometry II, {\em J. Diff. Geom.}, 107(2):327-371, 2017.
\bibitem{Engel} P. Engel, \newblock A proof of Looijenga's conjecture via integral-affine geometry,
{\em Journal of Differential Geometry}, 109(3): 467-495, 2018.

\bibitem{FE} P. Engel and R. Friedman, \newblock Smoothings and rational double point adjacencies for cusp singularities, {\em Journal of Diff. Geom.},118(1): 23-100 (May 2021).
\bibitem{Feehan} P.M.N. Feehan, \newblock Geometry of the ends of the moduli space of anti-self-dual connections, {\em J. Differential Geom.}, 1995, 42(3), 465-553.

\bibitem{Friedman}R. Friedman, Global smoothing of varieties with normal crossings, {\em Annals Math.}, 118 (1983), 75-114.
\bibitem{FS}R. Friedman and F. Scattone, \newblock Type III degeneration of K3 surfaces, {\em Invent. Math.} 83, 1-39, (1986).
\bibitem{FM}R. Friedman and R. Miranda \newblock Smoothing cusp singularities of small length, {\em Math. Ann.} 263, 185-212 (1983).

\bibitem{FHJ}X. Fu, H. Hein and X. Jiang, \newblock A continuous cusp closing process for negative K\"ahler--Einstein metrics, arXiv:2401.11468.

\bibitem{GHK15} M. Gross, P. Hacking and S. Keel, \newblock Mirror symmetry for log Calabi-Yau surfaces I, {\em Publ. Math. Inst. Hautes Etudes Sci.}, 122:65-168, (2015).

\bibitem{KasuyaCusp} N. Kasuya, \newblock The contact structure on the link of a cusp singularity, arXiv:1202.2198. \url{https://arxiv.org/abs/1202.2198}.

\bibitem{Hacking}P. Hacking, \newblock Compact moduli of plane curves, {\em Duke Math. Journal}, Vol. 124, No. 2 (2004), 213-257.

\bibitem{Hein}H.-J. Hein, \newblock Gravitational instantons from rational elliptic surfaces, {\em J. Amer. Math. Soc.} 25 (2012), no. 2, 355--393.


\bibitem{HMX14} C.D. Hacon, J. MacKernan, and C. Xu, \newblock Boundedness of moduli of varieties of general type,
 J. Euro. Math. Soc. 20 (2018), Issue 4, 865-901, arXiv:1412.1186.

 \bibitem{Hirzebruch} F. Hirzebruch, \newblock Hilbert modular surfaces, {\em Ens. Math.} 19 183-281 (1973).

\bibitem{Illusie2} L. Illusie, \newblock Cotangent complex and deformations of torsors and group schemes, in Toposes, algebraic geometry and logics, LNM 274, Springer, 1972, p. 159-189.

\bibitem{Inoue}M. Inoue, \newblock New surfaces with no meromorphic functions, II, In {\em Complex analysis and algebraic geometry}, 91-106, Iwanami Shoten, Tokyo, 1977.
% \bibitem{Jiang2008}Y. Jiang, \newblock The orbifold cohomology of simplicial toric stack bundles, {\em Illinois Journal of Mathematics},Vol. 52, No.2 (2008), 493-514, math.AG/0504563.
\bibitem{Jiang_2021}Y. Jiang, \newblock A note on higher obstruction spaces for surface singularities, preprint, arXiv:2112.10679.
\bibitem{Jiang_2022}Y. Jiang, \newblock The virtual fundamental class for the moduli space of surfaces of general type, arXiv:2206.00575.
\bibitem{Jiang_cusp}Y. Jiang, \newblock Equivariant smoothing for cusp singularities, preprint, arXiv:2302.00637.
\bibitem{Jiang_2023}Y. Jiang, \newblock Smoothing of surface singularities via equivariant smoothing of lci covers, arXiv:2309.16562.

\bibitem{Jiang_K3}Y. Jiang, \newblock Bubbling compactification of K3 surfaces,  in praparation. 
%\bibitem{Kollar}J. Koll\'ar, \newblock Moduli of varieties of general type, {\em Handbook of Moduli}: Volume II, Adv. Lect. Math. (ALM), vol. 24, Int. Press, Somerville, MA, 2013, pp. 115-130, arXiv:1008.0621.

\bibitem{Kollar-Shepherd-Barron} J. Koll\'ar and N. I. Shepherd-Barron, \newblock Threefolds and deformations of surface singularities,
{\em Invent. Math.}, 91, 299-338 (1988).
\bibitem{Kollar-Book}J. Koll\'ar, \newblock {\em Families of varieties of general type}, Cambridge Tracts in Mathematics, 231, Cambridge University Press,
April 2023,
ISBN:
9781009346115, DOI: https://doi.org/10.1017/9781009346115.

\bibitem{KP15} S. J. Kov\'acs and Z. Patakfalvi, \newblock Projectivity of the moduli space of stable log-varieties and subadditivity of log-Kodaira dimension, J. Amer. Math. Soc. 30 (2017), 959-1021, arXiv:1503.02952.


\bibitem{Laufer} H. Laufer, \newblock On minimally elliptic singularities, {\em Amer. J. Math.}, Vol. 99, No. 6 (1977), 1257-1295.

\bibitem{LWX}C. Li, X. Wang and C. Xu, \newblock Algebraicity of the metric tangent cone and equivariant K-stability, {\em Journal of the American Mathematical Society}, 34(4):1175-1214, 2021.

\bibitem{LT} J. Li and G. Tian, \newblock Virtual moduli cycles and Gromov-Witten invariants of algebraic varieties, {\em J. Amer. Math. Soc.}, 11, 119-174, 1998, math.AG/9602007.

\bibitem{Looijenga_Top}E. Looijenga, \newblock On the semi-universal deformation of a simple elliptic hypersurface singularity Part II: the discriminant, {\em Topology} Vol. 17, pp. 23-40, (1978).
\bibitem{Looijenga}E. Looijenga, \newblock Rational surfaces with an anticanonical cycle, {\em Ann. of Math.} (2), 114 (2):267-322, 1981.
\bibitem{MTT}D. Markushevich, A. Tikhomirov and G. Trautmann, \newblock Bubble-tree compactification of moduli spaces of vector bundles on surfaces, {\em Cent. Eur. J. Math.}, 10(4), (2012), 1331-1355.

\bibitem{MP}R. Miranda and U. Persson, \newblock On Extremal Rational Elliptic Surfaces, {\em Math. Z.}, 193, 537-558 (1986).


\bibitem{Nakamura} I. Nakamura, \newblock Inoue-Hirzebruch surfaces and a duality of hyperbolic unimodular singularities. I. {\em Math. Ann.} 252, 221-235 (1980).

\bibitem{Odaka}Y. Odaka, \newblock Algebraic geometry of bubbling K\"ahler metrics, arXiv:2406.14518.

\bibitem{Odaka2}Y. Odaka, \newblock Polystable log Calabi-Yau varieties and gravitational instantons, {J. Math. Sci. Univ. Tokyo}, 29, (2022) 1-50, arXiv:2009.13876.

\bibitem{PXu15} Z. Patakfalvi and C. Xu, \newblock Ampleness of CM line bundle on the moduli space of canonically polarized varieties,
{\em Algebraic Geometry}, 4(1), 29-39, (2017),
arXiv:1503.08668.
\bibitem{Pinkham} H. C. Pinkham, \newblock Automorphisms of cusps and Inoue-Hirzebruch surfaces, {\em Compositio Math.} 52, No.3 (1984), 299-313.

% \bibitem{PS} Y. Prokhorov and C. Shramov, \newblock Automorphism groups of Inoue and Kodaira surfaces, {\em Asian Journal of Mathematics}, Vol. 24 (2020)
%No. 2, 355-368.
\bibitem{Reid}M. Reid, \newblock Canonical 3-folds, in Algebraic Geometry Angers 1979, 273-310.

\bibitem{Stack_Project} Stack Project: https://stacks.math.columbia.edu/download/stacks-morphisms.pdf.

\bibitem{Shepherd-Barron} N. I. Shepherd-Barron, \newblock Degenerations with numerically effective canonical divisor, {\em The birational geometry of degenerations}, Progress in Mathematics, Vol. 29, 33-84.



\bibitem{Sun}S. Sun, \newblock Bubbling of Kähler-Einstein metrics, arXiv:2303.11309.
\bibitem{SZ}S. Sun and R. Zhang, \newblock Collapsing geometry of hyperk\"ahler 4-manifolds and applications,
{\em Acta Math.}, 232(2) 325-424 (2024).
\bibitem{Taubes}C. H. Taubes, \newblock A framework for Morse theory for the Yang--Mills functional, {\em Invent. Math.} 94 (1988), no. 2, 327--402.

\bibitem{Uhlenbeck} K. K. Uhlenbeck, \newblock Removable singularities in Yang-Mills fields, {\em Comm. Math. Phys.}, 1982, 83(1), 11-29.

\end{thebibliography}
\end{document}